\documentclass[twoside,a4paper]{article}
\usepackage[top=4cm, bottom=4cm, inner=3cm,outer=3cm]{geometry}
\usepackage[utf8]{inputenc}
\usepackage{amsthm}
\usepackage{mathtools}
\usepackage{csquotes}
\usepackage{amsmath}
\usepackage{graphicx}
\usepackage{dsfont}
\usepackage{soul}
\usepackage{amssymb}
\usepackage{multicol}
\usepackage{courier}
\usepackage{listings}
\usepackage{eurosym}
\usepackage{subfiles}
\usepackage{multicol}
\usepackage{pdfpages}
\usepackage{footnote}
\usepackage{color}
\definecolor{myblue}{rgb}{.25, .25, .9}
\definecolor{myred}{rgb}{.6, .4, .4}
\definecolor{myred2}{rgb}{.9, .1, .1}
\definecolor{mygreen}{rgb}{.25, .6, .5}
\usepackage{adjustbox}
\usepackage[colorlinks=true,
            linkcolor=myred2,
            urlcolor=myred,
            citecolor=myred]{hyperref}
\usepackage{pdfpages}
\usepackage[english]{babel}
\usepackage{array}
\usepackage{enumerate}
\usepackage{gensymb}
\usepackage[nottoc]{tocbibind}
\usepackage{lettrine}
\usepackage[normalem]{ulem}
\usepackage[font=small,labelfont=bf]{caption}
\usepackage{cancel}
\usepackage{nomencl}
\usepackage{mathrsfs}
\usepackage{changepage}
\usepackage{amsmath,stmaryrd,graphicx,float}
\numberwithin{equation}{section}

\usepackage{tikz-cd}
\usetikzlibrary{arrows}

\usepackage{subcaption} 

\usepackage{titlesec}

\newtheoremstyle{mystyle}
  {}
  {}
  {\itshape}
  {}
  {\bfseries}
  {.}
  { }
  {}

\theoremstyle{mystyle}

\newtheorem{theorem}{Theorem}[section]

\newtheorem{lemma}[theorem]{Lemma}
\newtheorem{proposition}[theorem]{Proposition}

\newtheorem{example}[theorem]{Example}
 
\newtheorem{algorithm}[theorem]{Algorithm} 
\newtheorem{remark}[theorem]{Remark} 
\newtheorem{assumption}[theorem]{Assumption}

\usepackage{etoolbox}

\DeclareMathAlphabet{\mathbbold}{U}{bbold}{m}{n}

\titleformat{\section}
{\large\filcenter\bfseries}{\llap{\thesection\hskip 9pt}}{0pt}{}
\titlespacing*{\section}
{0pt}{3.25ex plus 1ex minus .2ex}{1.5ex plus .2ex}

\titleformat{\subsection}[runin]
{\bfseries}{\llap{\thesubsection\hskip 9pt}}{0pt}{}
\titlespacing*{\subsection}
{0pt}{3.25ex plus 1ex minus .2ex}{1.5ex plus .2ex}

\titleformat{\subsubsection}[runin]
{\bfseries}{\llap{\thesubsubsection\hskip 9pt}}{0pt}{}
\titlespacing*{\subsubsection}
{0pt}{3.25ex plus 1ex minus .2ex}{1.5ex plus .2ex}

\titleformat{\paragraph}[runin]
{\bfseries}{\llap{\theparagraph\hskip 9pt}}{0pt}{}
\titlespacing*{\paragraph}
{0pt}{3.25ex plus 1ex minus .2ex}{1.5ex plus .2ex}

\usepackage[titles]{tocloft}

\usepackage{fancyhdr}
\usepackage{lastpage}
\fancypagestyle{plain}{
\fancyhf{}
\fancyhead[LE]{\thepage\quad\textbf{\nouppercase\leftmark}}
\fancyhead[RO]{\textbf{\nouppercase\rightmark}\quad\thepage}
}

\fancypagestyle{basicstyle}{%
\fancyhf{}
\fancyhead[LE]{\thepage}
\fancyhead[RO]{\thepage}
}

\DeclareMathOperator*{\subjectto}{\text{subject to}}
\DeclareMathOperator*{\minimize}{\text{minimize}}

\DeclareMathOperator*{\argmin}{argmin}

\usepackage{algorithm}
\usepackage{algorithmic}[1]

\usepackage{chngcntr}
\counterwithin{figure}{section}

\usepackage[backend=bibtex,
maxbibnames=99,
style=alphabetic,
bibencoding=ascii,
giveninits=true, 
]{biblatex}
\renewbibmacro{in:}{}

\usepackage{rotating}
\usepackage{mathtools}

\DeclarePairedDelimiter\floor{\lfloor}{\rfloor}

\title{{\textbf{\Large Imaginary Zeroth-Order Optimization}}}
\author{{Wouter Jongeneel}\footnote{The author is with the Risk Analytics and Optimization Chair, \'Ecole Polytechnique F\'ed\'erale de Lausanne (EPFL) and is supported by the Swiss National Science Foundation under the NCCR \textit{Automation}, grant agreement~51NF40\_180545. The author is grateful to RAO colleagues and ICCOPT 2022 participants, their comments greatly improved the work. Contact: \texttt{wouter.jongeneel@epfl.ch}, \url{wjongeneel.nl}}
}
\date{\small{upload: December 14, 2021, last update: August 14, 2022}}

\begin{document}
\maketitle
\thispagestyle{empty}
\begin{abstract}
Zeroth-order optimization methods are developed to overcome the practical hurdle of having knowledge of explicit derivatives. Instead, these schemes work with merely access to noisy functions evaluations. One of the predominant approaches is to mimic first-order methods by means of some gradient estimator. The theoretical limitations are well-understood, yet, as most of these methods rely on finite-differencing for shrinking differences, numerical cancellation can be catastrophic. The numerical community developed an efficient method to overcome this by passing to the complex domain. This approach has been recently adopted by the optimization community and in this work we analyze the practically relevant setting of dealing with computational noise. To exemplify the possibilities we focus on the strongly-convex optimization setting and provide a variety of non-asymptotic results, corroborated by numerical experiments, and end with local non-convex optimization.  
\end{abstract}

{\footnotesize{
\noindent\textbf{\textit{Keywords}}|zeroth-order optimization, derivative-free optimization, complex-step derivative, gradient estimation, numerical optimization.\\
\textbf{\textit{AMS Subject Classification (2020)}}|65D25, 65G50, 65K05, 65Y04, 65Y20, 90C56.
}}


\section{Introduction}
``\textit{La voie la plus courte et la meilleure entre deux vérités du domaine réel passe souvent par le domaine imaginaire.}''|J. Hadamard\footnote{See \url{http://homepage.math.uiowa.edu/~jorgen/hadamardquotesource.html}.}\\\\
From the Fourier transformation, quantum mechanics to the Nyquist stability criterion, the complex numbers grew out to be quintessential mathematical machinery. 

Building upon the work by~\cite{ref:kiefer1952stochastic,Ref:Lyness:Moler,nemirovsky1983problem,ref:SquireTrapp,ref:Flaxman,ref:nesterov2017random}, it is shown in~\cite{ref:JongeneelYueKuhnZO2021} that (randomized) zeroth-order optimization also benefits from passing to the complex domain as one can derive an inherently numerically stable method, which is in sharp contrast to common finite-difference methods. This work departs from~\cite{ref:JongeneelYueKuhnZO2021} by introducing an indispensable layer of realism; noise. 

We are interested in numerically solving optimization problems of the form 
\begin{equation*}
\label{equ:opt:main}
    \minimize_{x\in \mathcal{X}}\quad f(x),
\end{equation*}
where $f:\mathcal{D}\to \mathbb{R}$ is a smooth objective function defined on an open set~$\mathcal{D}\subseteq \mathbb{R}^n$, and~$\mathcal{X}\subseteq \mathcal{D}$ is a non-empty closed feasible set. Optimizers, which based on the context could be \textit{globally} or \textit{locally} optimal, are denoted by $x^{\star}$. We extend~\cite{ref:JongeneelYueKuhnZO2021} and assume that the objective function~$f$ can only be accessed through a \textit{\textbf{zeroth-order oracle}} that outputs \textit{corrupted} function evaluations at prescribed test points, that is, \textit{with} noise. As we only have access to such a zeroth-order oracle, our work belongs to the field of \textit{zeroth-order optimization}, \textit{derivative-free optimization} or more generally \textit{black-box optimization}~\cite{ref:Conn,ref:audet2017derivative}. 

We start by highlighting two important assumptions made throughout this work. 
\begin{assumption}[Smoothness]
\label{ass:Cw}
The objective function $f$ is real-analytic over $\mathcal{D}\subseteq \mathbb{R}^n$, denoted $f\in C^{\omega}(\mathcal{D})$. 
\end{assumption}

Recall, a function is \textit{\textbf{real-analytic}} when it can be locally expressed by a convergent power series, which is stronger than smoothness,~\textit{i.e.}, $C^{\omega}(\mathcal{D})\subset C^{\infty}(\mathcal{D})$. A complex-analytic function is called \textit{\textbf{holomorphic}}\footnote{More formally, a complex differentiable function is called \textit{holomorphic}, but as it turns out, complex differentiability coincides with complex analyticity~\cite{ref:Krantz:Complex}.}.  
With few exceptions~\cite{absil2005convergence}, Assumption~\ref{ass:Cw} does not appear often explicitly in the optimization literature. However, by means of the results in~\cite{polderman2} it does appear indirectly in for example the context of reinforcement-learning~\cite{ref:Fazel_18,ref:malik2019derivative}.
As in~\cite{ref:JongeneelYueKuhnZO2021}, Assumption~\ref{ass:Cw} is again mainly there to allow for the next assumption. As will be explained below, having access to $\Im(f(z))$ for some $z\in \mathbb{C}^n$ is at the core of the approach. In contrast to~\cite{ref:JongeneelYueKuhnZO2021} we allow for the presence of (computational) noise.

\begin{assumption}[Stochastic complex oracle]
\label{ass:complex:oracle}
Consider some \textit{unknown} function $f\in C^{\omega}(\mathcal{D})$ which admits a holomorphic extension to $\Omega\subseteq \mathbb{C}^n$. We assume to have access to an oracle which can output $\Re(f(z))+\xi$ and 
$\Im(f(z))+\xi$ for any $z\in \Omega$ with $\xi$ a zero-mean random variable supported on $\Xi\subseteq \mathbb{R}$ with $\mathbb{E}[\xi^2]\leq \sigma_{\xi}$ for some $\sigma_{\xi}>0$.
\end{assumption}

Assumption~\ref{ass:complex:oracle} is particularly important in the simulation-based context. As there the evaluation of $f(z)$ might pertain to millions of floating-point operations, chopping and round-off errors are easily introduced. The set $\Omega$ will be specified later on. We will make no further assumptions regarding the distribution of $\xi$.   

\subsection{Related work}
Arguably the first algorithm that uses noisy finite-differences to approximate gradient algorithms is the Kiefer-Wolfowitz algorithm~\cite{ref:kiefer1952stochastic},~\cite[Section~2.3.5]{ref:kushner1978stochastic}. \textcite{nemirovsky1983problem} contributed the first single-point gradient estimator and perhaps more importantly, the need for lower bounds. A large fraction of the work on zeroth-order optimization entails mimicking first-order algorithms via some approximation of the gradient. These types of algorithms are generally scalable\footnote{See however the discussion in~\cite{ref:scheinberg2022finite} to put this in the correct perspective.}, easy to implement and as they mimic first-order methods, they usually come with guarantees. A common gradient estimator is of the form
\begin{equation}
    \label{equ:g:example}
    \widehat{g}_{\delta}(x) = \frac{n}{\delta }\left(f(x+\delta u)-f(x)\right)u,
\end{equation}
for some choice of the smoothing parameter $\delta>0$ and some appropriately chosen random variable $u$. The $\delta$ is sometimes referred to as the \textit{exploration} parameter. See that~\eqref{equ:g:example} requires two function evaluations, as such we speak of a \textit{\textbf{multi-point method}}. Using estimators of the form~\eqref{equ:g:example} was popularized in the bandit-context~\cite{ref:Flaxman}, although for a single-point estimator, and relates largely to work on stochastic approximation algorithms~\cite{ref:kushner2003stochastic,ref:spall2005introduction} and to some extent to inexact/biased first-order methods~\cite{luo1993error,ref:Aspremont:2008smooth,ref:devolder2014first,ref:ajalloeian2021analysis}, \cite[Section~4]{ref:SDOT}. 

Compared to first-order methods, zeroth-order methods are commonly $O(n)$ times slower in the deterministic setting~\cite{ref:nesterov2017random}. When noise is involved the balance between bias and variance requires a more careful selection of the smoothing parameter $\delta$. 
Let $K\in \mathbb{N}$ denote the length of the sequence $x_1,x_2,\dots,x_K$ designed with the aim of converging in some sense to (some) $\argmin_{x\in \mathcal{X}} f(x)$. Let $\bar{x}_K=K^{-1}\sum^K_{k=1}x_k$ be a \textbf{\textit{uniformly-averaged iterate}}, we will be mostly interested in quantifying how fast the optimization error
\begin{equation*}
\label{equ:opt:error}
    \mathbb{E}[f(\bar{x}_K)-f(x^{\star})]
\end{equation*}
decays.
Here, the expectation $\mathbb{E}[\cdot]$ is over the oracle noise and the deliberate randomization within the proposed algorithms. That is, we can define an abstract probability space $(\Omega,\mathcal{F},\mathbb{P})$ and define $\mathbb{E}[\cdot]$ to be the expectation with respect to $\mathbb{P}$. 
In~\cite{ref:JamiesonNR12} the authors consider $\tau$-strongly convex functions with $L$-Lipschitz gradients and show that the expected optimization error decays like $\Omega(\sqrt{n/K})$ when using noisy single-point oracles. \textcite{shamir2013complexity} shows that in the quadratic case the result can be improved.  
If the objective $f$ is $r$-times continuously differentiable, Chen shows that a  rate of the order $O(K^{-(r-1)/2r})$ is optimal~\cite{ref:chen1988lower}. If $f$ is strongly convex the optimal rate becomes $O(K^{-(r-1)/r})$~\cite{polyak1990optimal}. 
See also~\cite{ref:RakhlinSS12} for more on optimal rates in the stochastic setting. 

In~\cite{duchi2015optimal} the authors show the information-theoretic optimality of multi-point (two-point) methods, yet,~in~\cite{ref:JongeneelYueKuhnZO2021} the authors show the numerical superiority of \textit{single}-point schemes. This work sets out to show to what extent this observation prevails when noise is present.
As highlighted throughout the recent survey article by~\textcite{ref:larson2019derivative}, it is not clear if there is a \textit{single}-point method which is as fast as multi-point methods. 
This observation motivates~\textcite{ref:zhang2020residualfeedback} to use some form of memory such that their estimator only demands a single new point each call. Nevertheless, in the end their method is reminiscent of a multi-point method. 
Another recent work observes how the continuous-time notion of~\textit{extremum seeking} can be translated to a zeroth-order optimization algorithm~\cite{ref:chen2021improve}. Their method turns out to be a combination of the aforementioned residual-feedback and momentum and achieves an optimization error of the order $O(n/K^{2/3})$, for a restrictive class of problems and a deterministic oracle. To the best of our knowledge, we will provide the first real single-point method which is capable of achieving an optimal rate.   

We focus on one particular approach to zeroth-order optimization. Different and successful lines of attack relate to \textit{model-based} (trust-region)~\cite{ref:Conn}, \textit{Bayesian}~\cite{ref:mockus2012bayesian} and more broadly \textit{black-box} optimization~\cite{ref:audet2017derivative}. 

\paragraph*{Contribution}
We show that catastrophic numerical cancellation errors are also inevitable in the widely used noisy multi-point case.
 We will show that this non-deterministic setting also benefits from the imaginary gradient estimator as proposed in~\cite{ref:JongeneelYueKuhnZO2021}.
 Using this \textit{single}-point estimator and building upon~\cite{ref:Bartlett2008,ref:akhavan2020exploiting}, we provide the non-asymptotic analysis for a variety of algorithms. 
Specifically, we consider for strongly convex functions the \textit{unconstrained}, \textit{constrained}, \textit{online} and \textit{quadratic} cases. In the last setting we can show that the algorithm is rate-optimal. To comply with zeroth-order knowledge we also propose an estimation scheme for the strong-convexity parameter. As an outlook we provide a local result in the nonconvex case and showcase PDE-constrained optimization as an area of application. Besides, we generalize some results from~\cite{ref:JongeneelYueKuhnZO2021} and we hope that a secondary contribution of this work is to bring numerical intricacies further to the attention. 

\paragraph*{Structure}
We start in Section~\ref{sec:num} by detailing numerical problems in zeroth-order optimization. In Section~\ref{sec:im:grad} we highlight the imaginary gradient estimator as proposed in~\cite{ref:JongeneelYueKuhnZO2021} to overcome the aforementioned obstacles. 
Section~\ref{sec:strong:cvx:opt} and Section~\ref{sec:non:convex} provide all algorithms, corresponding convergence rates and a few numerical experiments. Section~\ref{sec:smooth} briefly comments on merely smooth non-analytic functions and we conclude the work in Section~\ref{sec:conclusion}. Some auxiliary results can be found in the appendix.

\paragraph*{Notation}
The real and imaginary parts of a complex number $z=a+ib\in \mathbb{C}$ are denoted by $\Re(z)=a$ and $\Im(z)=b$, while $\mathbb{B}^n=\{x\in \mathbb{R}^n : \|x\|_2\leq 1\}$ is the Euclidean $n$-ball and $\mathbb{S}^{n-1}=\partial \mathbb{B}^n$ denotes the Euclidean $(n-1)$-sphere. Let $\mathsf{Y}\subset \mathbb{R}^n$ be a Borel measurable set such that $\partial \mathsf{Y}$ is an orientable compact differentiable manifold. We write $y\sim \mathsf{Y}$ to declare that~$y$ is a random vector following the uniform distribution on~$\mathsf{Y}$, and for any Borel measurable function $g:\mathsf{Y}\subset \mathbb R^n \rightarrow \mathbb R$ we denote by
\begin{equation*}
    \mathbb{E}_{y\sim \mathsf{Y}}[g(y)]= \frac{1}{\mathrm{vol}(\mathsf{Y})}\int_{\mathsf{Y}}g(y) \mathrm{d}V(y)
\end{equation*}
the expected value of~$g(y)$, where $\mathrm{d}V$ represents the Borel measure induced by the \textit{\textbf{volume form}} on $\mathsf{Y}$, and $\mathrm{vol}(\mathsf{Y})$ represents the volume of~$\mathsf{Y}$. The set of all $r$ times continuously differentiable real-valued functions on the open set~$\mathcal D\subseteq\mathbb R^n$ is denoted by $C^r(\mathcal{D})$.
Non-negative constants are denoted by $C_1,C_2,\dots,C_y,C_z$. Their values can change from line to line. Regarding complexity notation, $\Omega(\cdot)$, $\Theta(\cdot)$ and $O(\cdot)$ have their usual meaning with $\widetilde{O}(\cdot)$ hiding logarithmic factors.
The proof contain explicit errors, whenever possible. 

Using the notation from~\cite{nesterov2003introductory} a function $f$ is said to be $C^{k,r}_L(\mathcal{D})$-\textbf{\textit{smooth}} when $f$ is $k$ times continuously differentiable with additionally having its $r^{\mathrm{th}}$-derivative being $L$-Lipschitz over some open set $\mathcal{D}\subseteq \mathbb{R}^n$. Here, $k$ is an element of $\mathbb{N}_{\geq 0}\cup \{\infty\}\cup\{\omega\}$. That is, if $f\in C^{1,1}_{L_1(f)}(\mathcal{D})$, then, $f$ has a \textbf{\textit{Lipschitz gradient}},~\textit{i.e.}, 
\begin{equation}
    \label{equ:grad:Lipschitz}
    \|\nabla f(x) - \nabla f(y) \|_2 \leq L_1(f) \|x-y\|_2, \quad \forall x,y\in \mathcal{D}.
\end{equation}
Similarly, if $f\in C^{2,2}_{L_2(f)}(\mathcal{D})$, then, $f$ has a \textbf{\textit{Lipschitz Hessian}},~\textit{i.e.}, 
\begin{equation}
    \label{equ:Hess:Lip}
    \|\nabla^2 f(x)- \nabla^2 f(y)\|_2 \leq L_2(f)\|x-y\|_2 \quad \forall x,y \in \mathcal{D}. 
\end{equation}
Instead of the $\ell_2$-norm $\|\cdot\|_2$ one can generalize the above to any norm $\|\cdot\|$ and its dual $\|\cdot\|_{\star}$. 
Note that when $f\in C^{\omega}(\mathcal{D})$, then, the existence of $L_1(f)$ and $L_2(f)$ for $f$ restricted to compact subsets of $\mathcal{D}$ is trivial. Yet, to aid the reader, we will always indicate when we work with these constants.


\section{Numerical stability in zeroth-order optimization}
\label{sec:num}

Multi-point finite-difference estimators dominate the zeroth-optimization literature,~\textit{e.g.}, see~\cite{ref:hazan2014bandit,duchi2015optimal,ref:nesterov2017random,ref:gasnikov2017stochastic,ref:shamir2017optimal,ref:akhavan2020exploiting,ref:lam2021minimax,ref:novitskii2021improved} or the recent survey articles~\cite{ref:larson2019derivative,ref:IEEEsurvey}. The motivation largely follows from the observation that the initial single-point schemes as proposed in~\cite{nemirovsky1983problem,ref:Flaxman} have an unbounded variance, even when the function evaluations come without noise. The multi-point schemes avoid this by constructing estimators akin to numerical directional derivatives~\cite{ref:agarwal2010optimal,ref:nesterov2017random}.

Nevertheless, as pointed out in~\cite{ref:JongeneelYueKuhnZO2021}, multi-point schemes \textit{do} suffer from catastrophic numerical cancellation. See also~\cite{ref:shi2021numerical} for an extensive numerical study on the numerical performance of finite-difference methods in the context of optimization. 

\subsection{Numerical cancellation}
 The smallest $\epsilon_M\in \mathbb{Q}_{>0}$ such that on a particular machine $1+\epsilon_M>1$ is called the \textit{\textbf{machine precision}}. Nowadays, the number $\epsilon_M$ is commonly of the order $10^{-16}$, which is the number we will use. So in general, for a continuous function $f:\mathbb{R}\to \mathbb{R}$, when $x,y\in \mathbb{R}$ are chosen such that $f(x)-f(y)\leq \epsilon_M$ the numerical evaluation of $f(x)-f(y)$ can be problematic. 
Now for zeroth-order gradient estimators, given some $f\in C^{r+1}(\mathbb{R})$ with $r\geq 1$, then in the approximation 
\begin{equation}
\label{equ:part:f}
    \partial_x f(x) = \frac{f(x+\delta)-f(x)}{\delta}+O(\delta) 
\end{equation}
one cannot make $\delta>0$ arbitrarily small and expect to recover $\partial_x f(x)$. For a sufficiently small $\delta$ the evaluations $f(x+\delta)$ and $f(x)$ will be numerically indistinguishable and cancellation errors appear, see~\cite[Chapter 11]{ref:overton2001numerical}. Running into these machine-precision problems is inherent to finite-difference (multi-point) optimization methods as one looks for (at) the flattest part of $f$. 

A celebrated work-around in the numerical community is the so-called \textit{\textbf{complex-step method}}. This approach was introduced in~\cite{Ref:Lyness:Moler} with the first concrete \textit{complex-step} approach appearing in~\cite{ref:SquireTrapp} and with later elaborations to higher-order derivatives, matrices and Lie groups in~\cite{ref:Martins,ref:Higham,ref:Abreu:geo,ref:Abreu:wave:2018,ref:cossette2020complex}. In short, via the Cauchy-Riemann equations one can show that for a holomorphic function $f$, one has
\begin{equation}
\label{equ:cs:deriv}
    \partial_x f(x) = \frac{\Im(f(x+i\delta))}{\delta}+O(\delta^2).
\end{equation}
Not only is numerical cancellation impossible, the error term improved compared to~\eqref{equ:part:f}. 
This approach recently surfaced in the optimization community~\cite{ref:complexwithnoise,ref:Hare2020} with the first complete deterministic non-asymptotic analysis appearing in~\cite{ref:JongeneelYueKuhnZO2021}.
The first applications of the complex-step derivative to Reinforcement Learning appeared in~\cite{ref:wang2021improved,ref:wang2021model}. Of course, as complex arithmetic is more expensive than real arithmetic, numerical stability does not come for free\footnote{For example, to compute the multiplication of 2 complex numbers $(a_1+ib_1)(a_2+ib_2)$ one needs 3 real multiplications~,\textit{i.e.}, $(a_1+b_1)(a_2+b_2)$, $a_1a_2$ and $b_1b_2$, see also~\cite{ref:alt1981complexity}.}. 

To visualize the power of the complex-step approach we provide a short example. 
\begin{example}[Numerical estimator stability]
\label{ex:num:est:stab}
\upshape{
We showcase the forward-difference ($\mathsf{fd}$), central-difference ($\mathsf{cd}$) and complex-step ($\mathsf{cs}$) for $f(x)=\log(x)$ at $x=1$, that is, we compare
\begin{align*}
       f_{\mathsf{fd}}(x,\delta) &= \frac{f(x+\delta)-f(x)}{\delta},\\  f_{\mathsf{cd}}(x,\delta) &= \frac{f(x+\delta)-f(x-\delta)}{2\delta},\\  f_{\mathsf{cs}}(x,\delta) &= \frac{\Im\left( f(x+i\delta) \right)}{\delta}
\end{align*}
for $\delta\downarrow 0$, see~Figure~\ref{fig:estlog}. 
\begin{figure*}[t!]
    \centering
    \begin{subfigure}[b]{0.3\textwidth}
        \includegraphics[width=\textwidth]{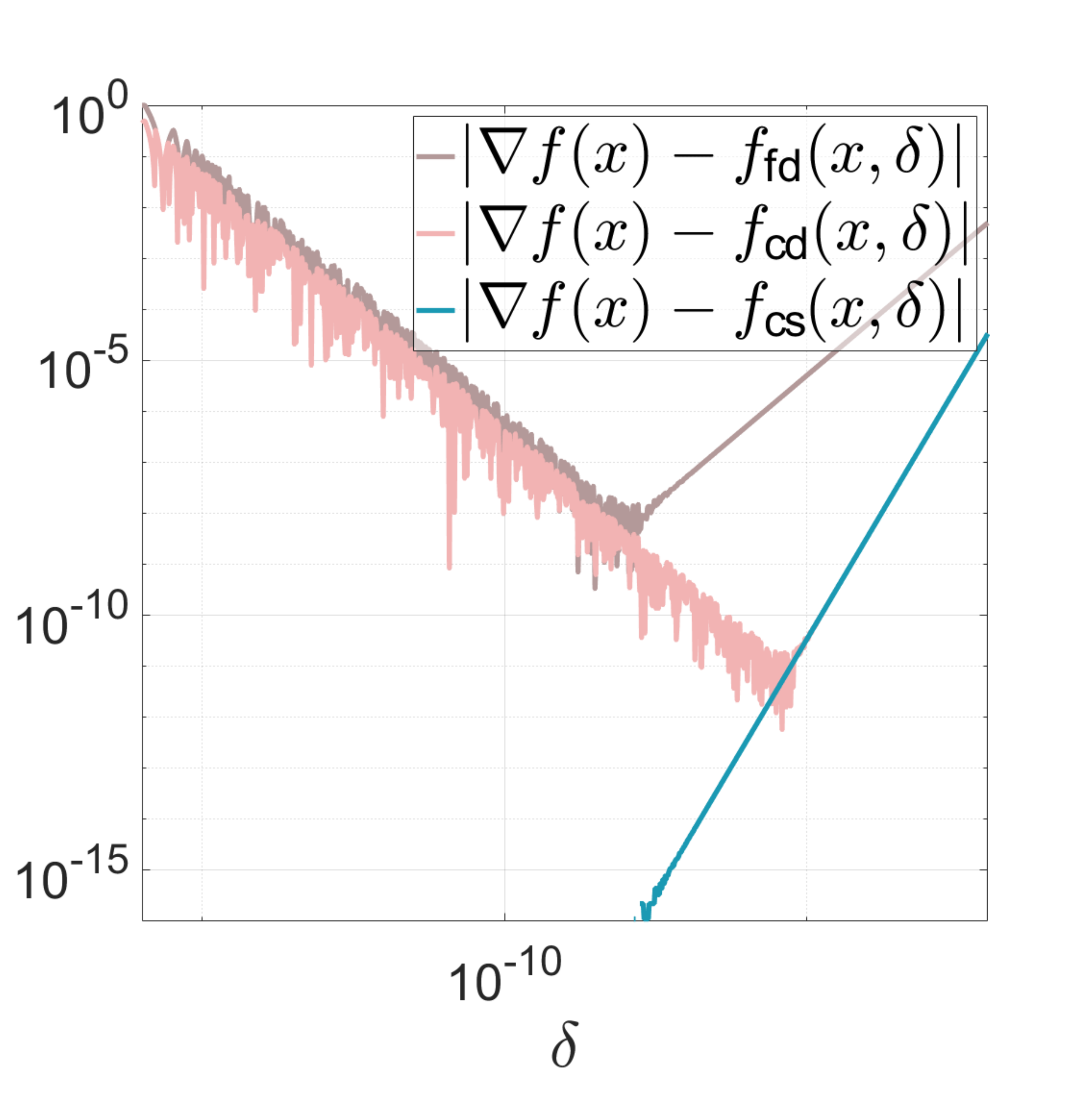}
        \caption{Example~\ref{ex:num:est:stab}, gradient estimator comparisons. See also~\cite[Example~2.6]{ref:JongeneelYueKuhnZO2021}.}
        \label{fig:estlog}
    \end{subfigure}\quad
    \begin{subfigure}[b]{0.3\textwidth}
        \includegraphics[width=\textwidth]{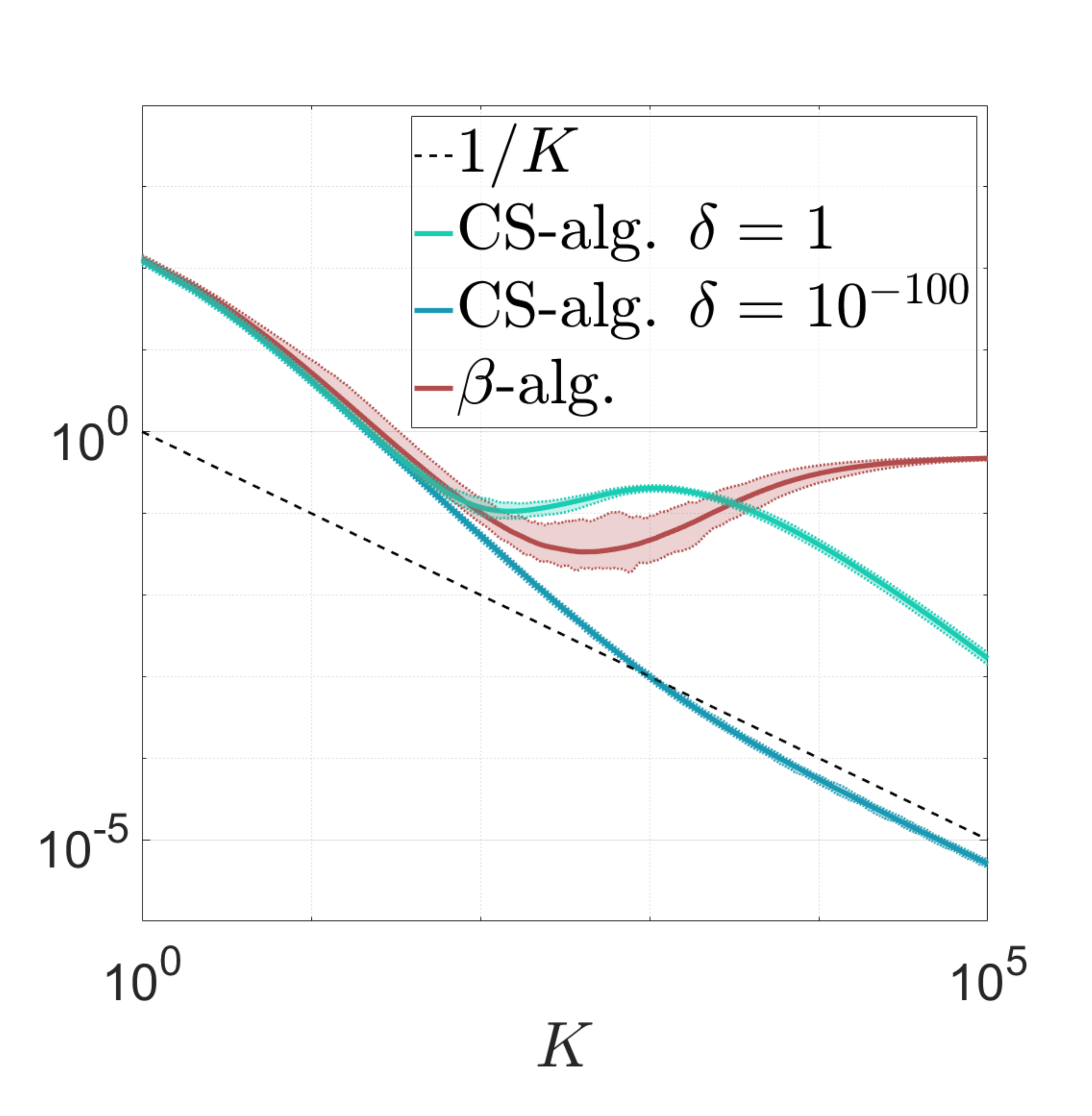}
        \caption{Suboptimality gap for Example~\ref{ex:strongcvx}, Theorem~\ref{thm:strong:convex:noise:quadratic} (CS) vs.~\cite[Theorem~5.1]{ref:akhavan2020exploiting} ($\beta$).}
        \label{fig:strongcvx}
    \end{subfigure}\quad 
    \begin{subfigure}[b]{0.3\textwidth}
        \includegraphics[width=\textwidth]{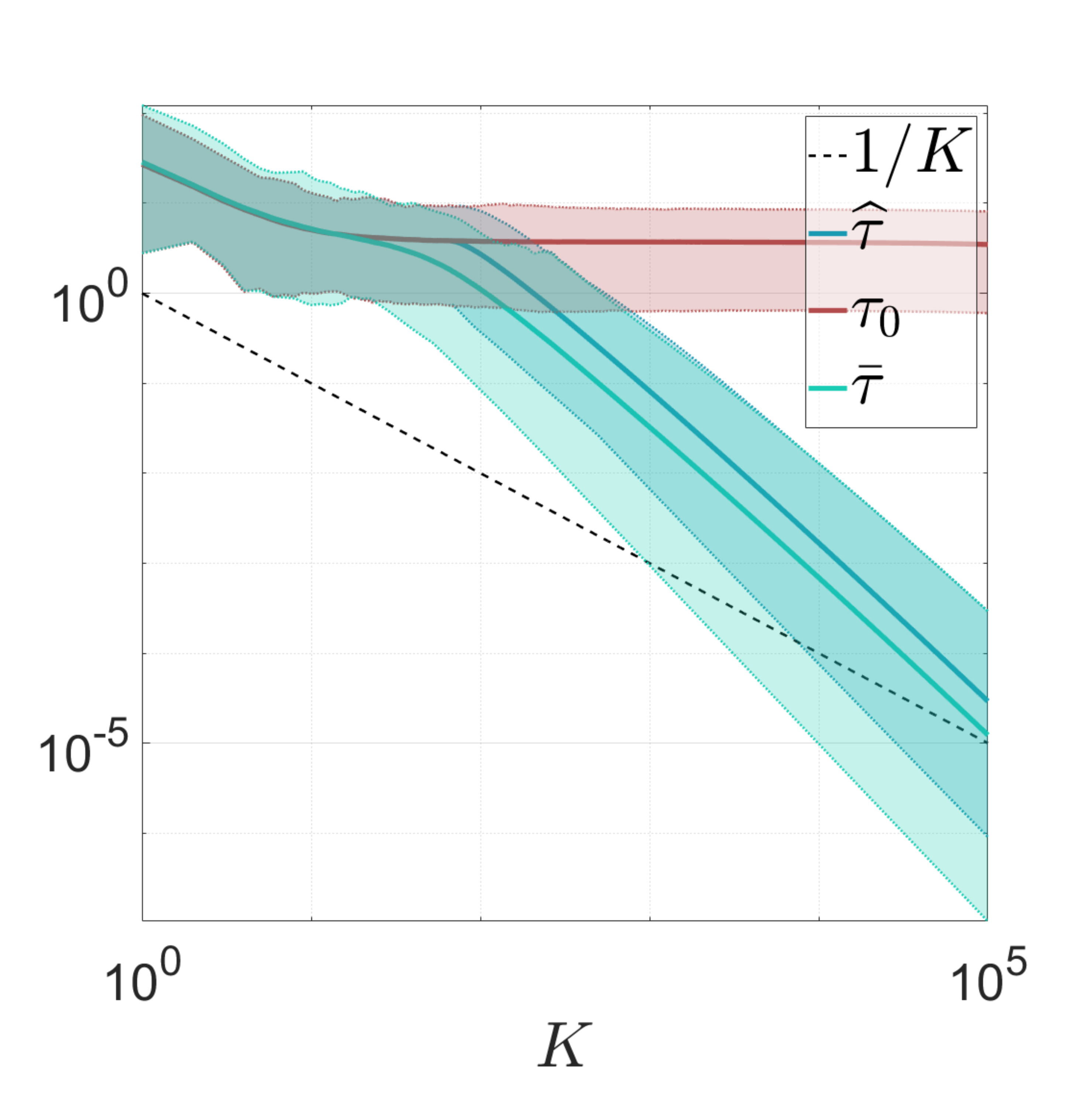}
        \caption{Suboptimality gap for Example~\ref{ex:reg}, Algorithm~\ref{alg:convex_unconstrained}~(b) combined with the estimation scheme~\eqref{equ:SDP:tau}.}
        \label{fig:cvxreg}
    \end{subfigure}
    \caption[]{Numerical experiments. Each figure displays all available data.}
    \label{fig:sim}
\end{figure*}
Only the complex-step estimator can reach machine precision, yet the other two methods are used frequently in zeroth-order optimization under the assumption that one can select $\delta$ arbitrarily close to $0$. As such, these methods leave something to be desired, numerically. 
}
\end{example}

\begin{figure*}[t!]
    \centering
    \begin{subfigure}[b]{0.23\textwidth}
        \includegraphics[width=\textwidth]{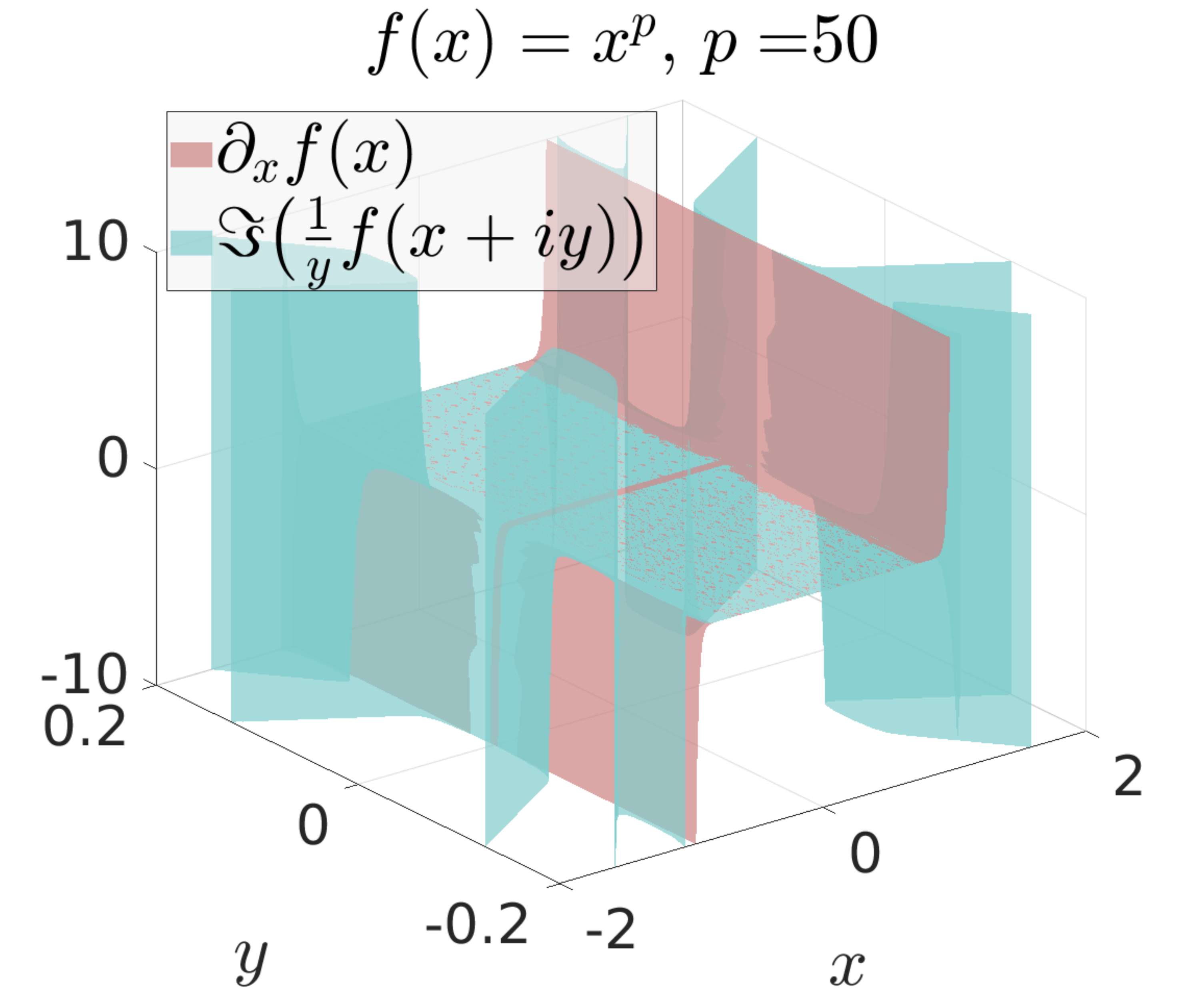}
        \caption{$p=50$}
        \label{fig:p50}
    \end{subfigure}\quad
    \begin{subfigure}[b]{0.23\textwidth}
        \includegraphics[width=\textwidth]{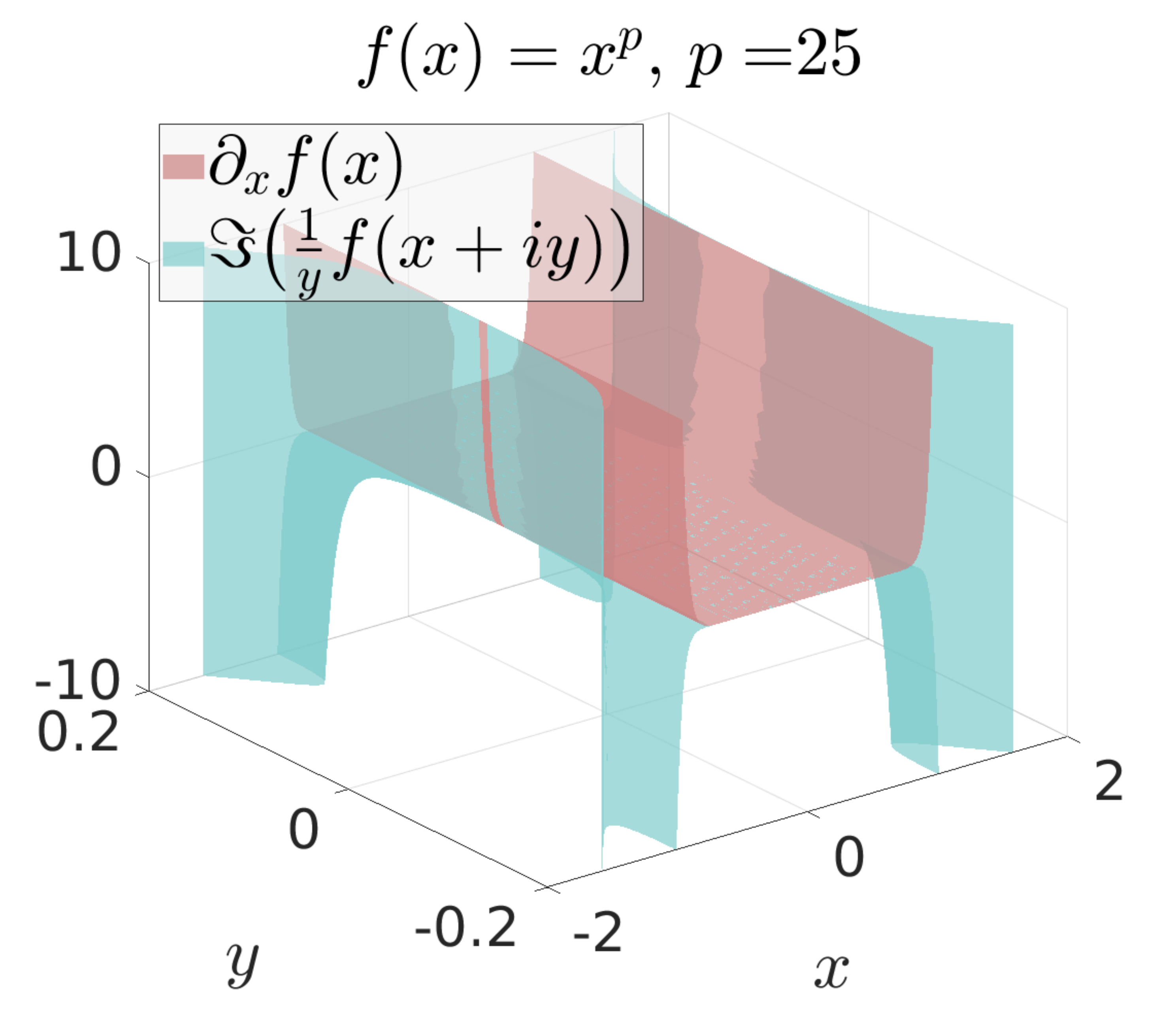}
        \caption{$p=25$}
        \label{fig:p25}
    \end{subfigure}\quad    \begin{subfigure}[b]{0.23\textwidth}
        \includegraphics[width=\textwidth]{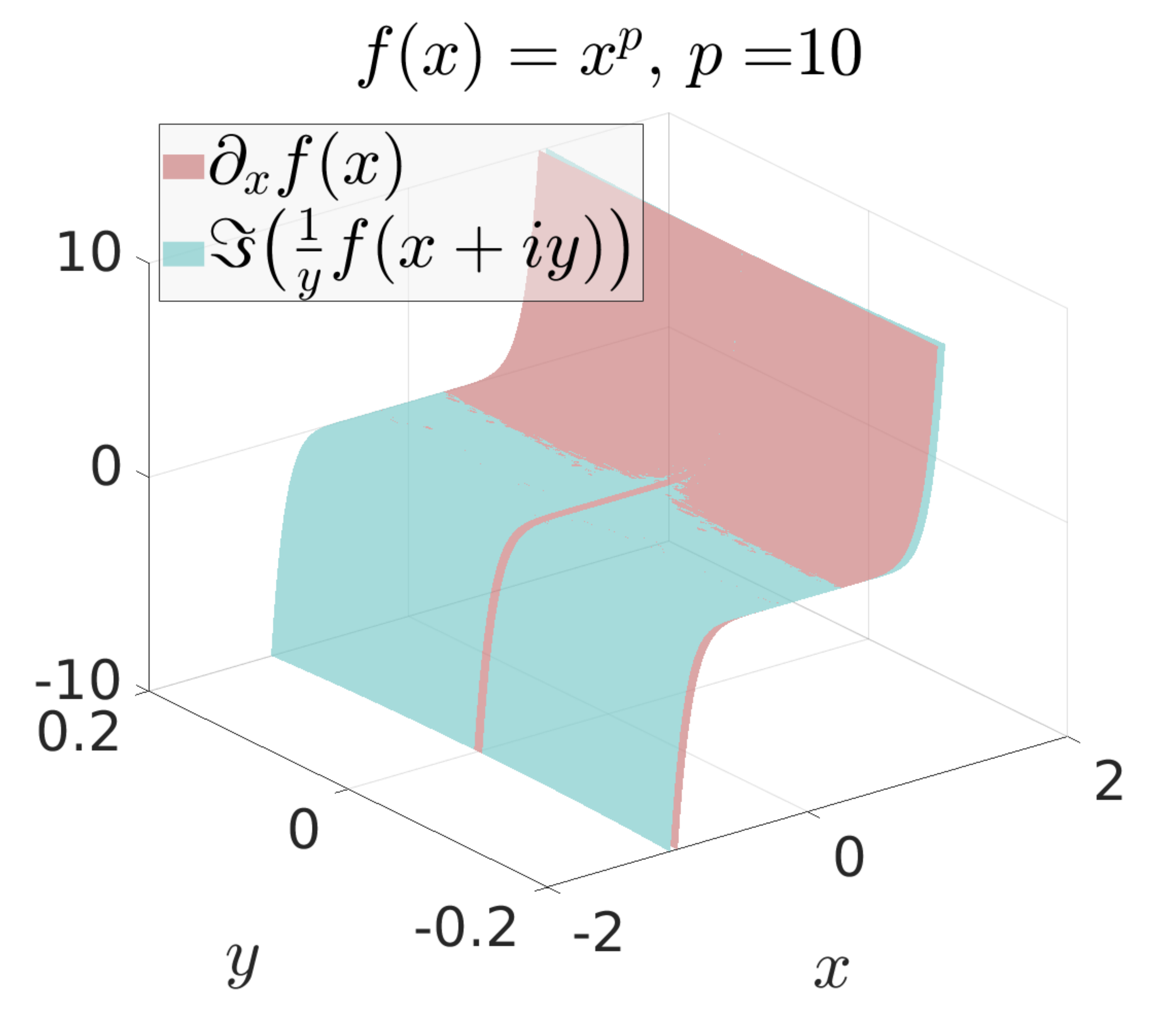}
        \caption{$p=10$}
        \label{fig:p10}
    \end{subfigure}\quad    \begin{subfigure}[b]{0.23\textwidth}
        \includegraphics[width=\textwidth]{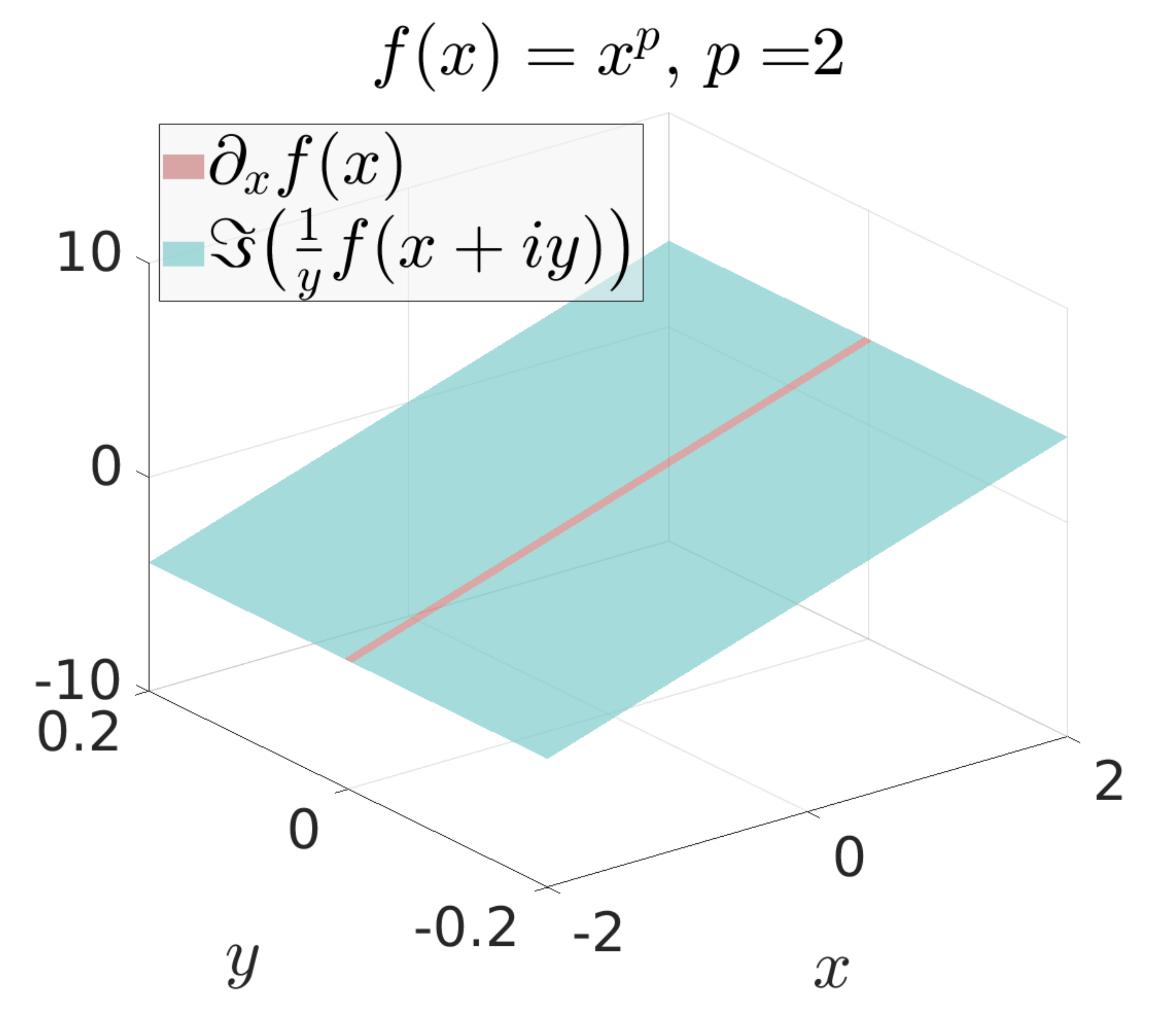}
        \caption{$p=2$}
        \label{fig:p2}
    \end{subfigure}
    \caption[]{Further visualizations of~\eqref{equ:cs:deriv},~\textit{i.e.}, $\partial_x f(x)\approx \Im f(x+iy)/y$.}
    \label{fig:sim:Im}
\end{figure*}

At last we elaborate on Example~\ref{ex:num:est:stab} and visualize the imaginary lifting of $f(x)$. That is, for $f(x)=x^p$, with $x\in \mathbb{R}$ and $p\in \mathbb{N}$, we show $\Im(f(x+iy)/y)$. Indeed, for sufficiently small $y$ we see in Figure~\ref{fig:sim:Im} that this number converges to $\partial_x f(x)$ for $p\downarrow 2$\footnote{See \url{http://wjongeneel.nl/ZO.gif} for an animated version of Figure~\ref{fig:sim:Im}.}. 


\section{Imaginary gradient estimation}
\label{sec:im:grad}

In this section we summarize the main tool as set forth by~\cite{ref:JongeneelYueKuhnZO2021}.
Motivated by Example~\ref{ex:num:est:stab}, we consider the \textbf{\textit{imaginary}} $\delta$-\textbf{\textit{smoothed}} version of $f$ as proposed in~\cite{ref:JongeneelYueKuhnZO2021}, that is
\begin{equation}
\label{equ:fw:Re}
    f_{\delta}(x) = \mathbb{E}_{v\sim \mathbb{B}^n} \left[\Re\big(f(x+i\delta v)\big) \right]. 
\end{equation}
Here, the parameter $\delta\in \mathbb{R}_{> 0}$ is the tuneable smoothing parameter and relates to the radius of the ball we average over. As mentioned before, the offset $\delta v$ in~\eqref{equ:fw:Re} relates to exploration\footnote{This notion of exploration \textit{could} be a benefit of these randomized approaches~\cite{ref:scheinberg2022finite}.}, due to our limited amount of information on the objective, each direction is potentially worthwhile exploring and as such we consider a perfectly symmetric shape; the ball $\mathbb{B}^n$. See~\cite{ref:hazan2014bandit} and Lemma~\ref{lem:gen:CR:grad} for comments and results beyond $\mathbb{B}^n$.

To make sure $f_{\delta}$ is well-defined, $f(x+i\delta v)$ needs to be well-defined and as such we assume the following. 
\begin{assumption}[Holomorphic extension]
\label{ass:holomorphic}
The function $f:\mathcal{D}\subseteq \mathbb{R}^n\to \mathbb{R}$ is real-analytic over the open set $\mathcal{D}$ and admits a holomorphic extension to $\mathcal{D}\times i\cdot (-\bar{\delta},\bar{\delta})^n\subset \mathbb{C}^n$ for some $\bar{\delta}\in (0,1)$. 
\end{assumption}
See~\cite[Section~2.1]{ref:JongeneelYueKuhnZO2021} for more on the existence of such an extension. Note, the interval $(0,1)$ is merely a convenient choice for the exposition.  

Next we highlight the approximation quality of $f_{\delta}$. 
\begin{lemma}[Approximation quality of the complex-step function~{\cite[Proposition~3.2]{ref:JongeneelYueKuhnZO2021}}]
\label{lem:approx:function}
Let $f\in C^{\omega,1}_{L_1(f)}(\mathcal{D})$ with $L_1(f)\geq 0$ satisfy Assumption~\ref{ass:holomorphic} for some $\bar{\delta}\in (0,1)$. Then, for $f_{\delta}$ as in~\eqref{equ:fw:Re} and any fixed $x\in \mathcal{D}$ and $\kappa\in (0,1)$ there exists some constant $C_0\geq 0$, vanishing with $L_1(f)$, such that
\begin{align}
\label{ineq:f_delta_deviation}
    \left|f_\delta (x) - f(x) \right| &\leq {C_0 \delta^2  }\quad \forall\,\delta\in (0,\kappa \bar{\delta}]. 
\end{align}
\end{lemma}

It is imperative to remark that convexity of $f$ does not always carry over to $f_{\delta}$,~\textit{e.g.,} see~\cite[Example~3.6]{ref:JongeneelYueKuhnZO2021}. 

Now we state one of the key contributions of~\cite{ref:JongeneelYueKuhnZO2021}, which is the integral representation of $\nabla f_{\delta}$. This result is the complex-step version of the approach as proposed in~\cite[Section 9.3]{nemirovsky1983problem} and popularized by~{\cite[Lemma 1]{ref:Flaxman}}.

\begin{lemma}[The gradient of the complex-step function~{\cite[Proposition~3.3]{ref:JongeneelYueKuhnZO2021}}]
\label{lem:CR:grad}
Let $f\in C^{\omega}(\mathcal{D})$ satisfy Assumption~\ref{ass:holomorphic} for some $\bar{\delta}\in (0,1)$, then, $f_{\delta}$ as in~\eqref{equ:fw:Re} is differentiable and for any $x\in \mathcal{D}$ we have for any $\delta\in (0,\bar{\delta})$
\begin{equation}
\label{equ:unif:grad}
    \nabla f_{\delta}(x) = \frac{n}{\delta}\cdot \mathbb{E}_{u\sim \mathbb{S}^{n-1}}\left[\Im\left(f(x+i\delta u) \right)u \right]. 
\end{equation}
\end{lemma}

The following result allows for showing consistency,~\textit{i.e.}, $\lim_{\delta\downarrow 0}\nabla f_{\delta}(x)=\nabla f(x)$. 
\begin{lemma}[Integration over the $(n-1)$-sphere]
\label{lem:int:sphere}
Given any $x\in \mathbb{R}^n$, then
\begin{equation}
\label{equ:int:sphere}
    \frac{n}{\mathrm{vol}(\mathbb{S}^{n-1})}\cdot \int_{\mathbb{S}^{n-1}}\langle x, u \rangle u \mathrm{d}V(u) = x.
\end{equation}
\end{lemma}
Although this result is well-known, for completeness we also provide the proof.
\begin{proof}
First, rewrite~\eqref{equ:int:sphere} as $n\cdot \int_{\mathbb{S}^{n-1}}uu^{\mathsf{T}}du\,x$ and recall that $uu^{\mathsf{T}} x = \langle x, u \rangle u$. Now we would like to show that $n\cdot \int_{\mathbb{S}^{n-1}}uu^{\mathsf{T}} du = \mathrm{vol}(\mathbb{S}^{n-1})\cdot I_n$.
To that end, use the geometric tracing identity
$n\cdot \int_{\mathbb{S}^{n-1}} \langle Xu, u \rangle du = \mathrm{Tr}(X)\cdot \mathrm{vol}(\mathbb{S}^{n-1})$~\cite[Lemma~3.100]{ref:gallot2004riemannian}, 
differentiating both sides with respect to $X$ yields $n\cdot\int_{\mathbb{S}^{n-1}}uu^{\mathsf{T}}du = \mathrm{vol}(\mathbb{S}^{n-1})\cdot I_n$ indeed, which concludes the proof.
\end{proof}

Since $f$ is real-analytic, the directional derivative at $x\in \mathcal{D}$ in the direction $u\in \mathbb{S}^{n-1}$ is well-defined and given by $\langle \nabla f(x), u \rangle$.
Then, observe from~\eqref{equ:unif:grad} that the approximation is asymptotically consistent, that is, by appealing to the dominated convergence theorem we have  
\begin{equation}
\label{equ:unif:consistent}
    \lim_{\delta\downarrow 0}\nabla f_{\delta} (x) = \frac{n}{\mathrm{vol}(\mathbb{S}^{n-1})}\cdot \int_{\mathbb{S}^{n-1}} \langle \nabla f(x), u \rangle u \mathrm{d}V(u) \overset{\eqref{equ:int:sphere}}{=} \nabla f(x). 
\end{equation}
Showing consistency of this type, albeit for the estimator, was one of the key observations in~\cite{ref:agarwal2010optimal,ref:nesterov2017random} to reduce gradient estimator variance.
Such an observation does not hold for other known single-point estimators \textit{cf.}~\cite[Section 1.1]{ref:Flaxman}.  
 
Lemma~\ref{lem:CR:grad} provides us immediately with a (noisy) single-point estimator of $\nabla f_{\delta}(x)$, namely 
\begin{equation}
\label{equ:grad:est:g}
    g_{\delta}(x) = \frac{n}{\delta}\Im\left(f(x + i \delta u)\right)u+\frac{n}{\delta}\xi u,\quad u\sim \mathbb{S}^{n-1} 
\end{equation}
for some noise term $\xi\in \Xi$. 
In contrast to the noise-free setting in~\cite{ref:JongeneelYueKuhnZO2021}, equation~\eqref{equ:grad:est:g} immediately reveals the delicacy in selecting $\delta\in \mathbb{R}_{>0}$. Note, the term $n/\delta$ follows from our choice to average over $\mathbb{B}^n$,~\textit{i.e.}, by~\eqref{equ:fw:Re}. Below we will clarify that this term, and thereby the offset due to the noise, cannot be decreased by any other choice of solid. In that sense, $\mathbb{B}^n$ is \textit{geometrically optimal}. We will use~\eqref{equ:grad:est:g} in gradient descent algorithms of the form $x_{k+1}=x_k-\mu_k g_{\delta_k}(x_k)$, as detailed in Algorithm~\ref{alg:convex_unconstrained}~(a) and Algorithm~\ref{alg:convex_unconstrained}~(b), for $\mu_k\in \mathbb{R}$ a stepsize and $\delta_k\in \mathbb{R}_{>0}$ the smoothing parameter. 

The next assumption on the (computational) noise will be assumed throughout. 
\begin{assumption}[Independence]
\label{ass:xi:indep}
The random variable $\xi$ is drawn independently of $u\sim \mathbb{S}^{n-1}$.
\end{assumption}

\begin{proposition}[Gradient approximation quality~{\cite[Proposition~3.4]{ref:JongeneelYueKuhnZO2021}}]
\label{prop:grad:approx}
Let $f\in  C^{\omega,2}_{L_2(f)}(\mathcal{D})$ with $L_2(f)\geq 0$ satisfy Assumption~\ref{ass:holomorphic} for some $\bar{\delta}\in (0,1)$. Then, for any fixed $x\in \mathcal{D}$ and $\kappa\in (0,1)$ there is a constant $C_1\geq 0$, vanishing with $L_2(f)$, such that
\begin{align}
\label{equ:grad:delta:approx:error}
    \|\nabla f_{\delta}(x) - \nabla f(x) \|_2 &\leq {C_1 n  \delta^2 }\quad \forall \delta\in (0,\kappa\bar{\delta}]. 
\end{align}
\end{proposition}

We see that the simple \textit{single}-point approach allows for an error of the form $O(\delta^2)$ which is what can be commonly achieved using central-difference \textit{multi}-points methods~\textit{cf.}~\cite{ref:nesterov2017random}. 

From~\eqref{equ:grad:delta:approx:error} it appears that~\eqref{equ:grad:est:g} is potentially a biased gradient estimator. Consider the special case of $f$ being quadratic (see Figure~\ref{fig:sim:Im} for a visualization). In that case, $\nabla f_{\delta}=\nabla f$, that is, the estimator is unbiased: $\mathbb{E}[g_{\delta}]=\nabla f$. This property will be exploited in Section~\ref{sec:strong:convex:opt:opt}. 

In general, however, there will be a bias, controlled in part by selecting the sequence $\{\delta_k\}_{k\geq 1}$ and unfortunately, a fixed bias prohibits (local) convergence in general~\cite{ref:ajalloeian2021analysis}. However, by looking at~\eqref{equ:grad:est:g}, it can be shown that to overcome this, a selection of $\{\mu_k\}_{k\geq 1}$ and $\{\delta_k\}_{k\geq 1}$ should satisfy the following;  
\begin{enumerate}[(i)]
    \item\label{pt:1} As $\mu_k=\Theta(k^{-1})$~\cite{ref:RakhlinSS12}, for fixed $\delta_k= \delta>0$ a bias term prevails of the form $\sum^K_{k=1}\mu_k\delta = O(\log(K)+1)$. This can be avoided by selecting $\delta_k$ to be asymptotically vanishing.
    \item However, as the data is noisy, a term of the form $\mu_k/\delta_k$ also accumulates. As such, by~\eqref{pt:1} $\delta_k \to 0$, but slower than $\mu_k\to 0$. 
\end{enumerate}
With this in mind we see that when $\mathbb{E}[g_{\delta}]\neq \nabla f$ zeroth-order optimization algorithms resort to selecting the smoothing-parameter sequence $\{\delta_k\}_{k\geq 1}$ such that $\delta_1$ converges to $0$ sufficiently slow,~\textit{cf.}~\cite[Theorem~1]{ref:novitskii2021improved},~\cite[Theorem~3]{ref:balasubramanian2021zeroth}. See also~\cite{ref:Fabian1971optimizing},~\cite[Chapter~6]{ref:spall2005introduction},~\cite[Assumption~1]{ref:wang2021model} for similar assumptions from the stochastic approximation viewpoint. Motivated by the observation that $\delta_k\to 0$ is necessary for an abundance of algorithms, this work provides a framework that can handle this requirement numerically. That means, a framework where $\delta_k$ can be made arbitrarily small\footnote{Up to what the machine at hand can produce, usually $2^{-1023}\approx 10^{-308}$.}. 

At last, to characterize the effectiveness of our algorithms, we need to bound the second moment of the estimator~\eqref{equ:grad:est:g}. 
We observe the same attractive property as highlighted in~\cite{ref:nesterov2017random}, there is no need to assume boundedness of the second moment of our stochastic estimator, \textit{cf.}~\cite{ref:RakhlinSS12}. 
As we allow for computational noise, the bound will differ slightly from the result in~\cite{ref:JongeneelYueKuhnZO2021}. 

\begin{lemma}[Estimator second moment]
\label{lem:im:lip}
Let $f\in C^{\omega,2}_{{L_2}(f)}(\mathcal{D})$ satisfy Assumption~\ref{ass:holomorphic} for some $\bar{\delta}\in (0,1)$ and ${L_2}(f)\geq 0$. Then, for any fixed $x\in \mathcal{D}$, $\kappa\in (0,1)$ and $g_{\delta}(x)$ as in~\eqref{equ:grad:est:g} there are constants $C_a,C_b\geq 0$, vanishing with $L_2(f)$, such that for any $\delta\in (0,\kappa \bar{\delta}]$ one has 
\begin{align}
\label{equ:g:var}
\mathbb{E}_{u\sim \mathbb{S}^{n-1}}\left[ \|g_{\delta}(x)\|_2^2\right] &\leq {C_a n^2 \delta^4} + C_b n^2\delta^2 \|\nabla f(x)\|_2 + n\|\nabla f(x)\|_2^2+ \tfrac{n^2}{\delta^2}\sigma_{\xi}.
\end{align}
\end{lemma}
\begin{proof}
First, observe from Algorithm~\ref{alg:convex_unconstrained}, Assumption~\ref{ass:complex:oracle} and Assumption~\ref{ass:xi:indep} that 
\begin{equation*}
    \mathbb{E}_{u\sim \mathbb{S}^{n-1}}\left[ \|g_{\delta}(x)\|_2^2\right] = \tfrac{n^2}{\delta^2}\mathbb{E}_{u\sim \mathbb{S}^{n-1}}\left[\left(\Im\left(f(x+i\delta u) \right)\right)^2 \right]+\tfrac{n^2}{\delta^2}\mathbb{E}_{\xi}[\xi^2].
\end{equation*}
Then, the claim follows directly by the same reasoning as in~\cite[Corollary~3.5]{ref:JongeneelYueKuhnZO2021}.
\end{proof}

As with standard gradient-descent, the more isotropic the level sets of the objective are, the better. The common way to enforce this is by means of changing the underlying metric via the Hessian,~\textit{i.e.}, Newton's method.
With this in mind, averaging over some solid ellipsoid might appear more beneficial than over the ball. In the spirit of~\cite{ref:hazan2014bandit} and~\cite[Proposition~3,Lemma~4]{ref:hu2016bandit} we generalize Lemma~\ref{lem:CR:grad} to more generic solids and show|perhaps unsurprisingly| that spherical smoothing is optimal in the sense that it minimizes the offset due to noise in~\eqref{equ:g:var}. 

To be in line with Assumption~\ref{ass:holomorphic} we assume that this generic solid $\mathsf{M}$ is a subset of $(-1,1)^n$.

\begin{lemma}[The gradient of the complex-step function for generic solids]
\label{lem:gen:CR:grad}
Let $\mathsf{M}\subset (-1,1)^n \subset \mathbb{R}^n$ be diffeomorphic to $\mathbb{B}^n$. 
Let $f\in C^{\omega}(\mathcal{D})$ satisfy Assumption~\ref{ass:holomorphic} for some $\bar{\delta}\in (0,1)$, then, $f_{\delta,\mathsf{M}}$ as in 
\begin{subequations}
\begin{equation}
    \label{equ:gen:f:delta}
    f_{\delta,\mathsf{M}}(x) = \mathbb{E}_{v\sim \mathsf{M}}\left[\Re(f(x+i\delta v) \right]
\end{equation}
is differentiable and for any $x\in \mathcal{D}$ we have for any $\delta\in (0,\bar{\delta})$
\begin{equation}
\label{equ:gen:grad:f:delta}
    \nabla f_{\delta,\mathsf{M}}(x) = \frac{\mathrm{vol}(\delta\partial\mathsf{M})}{\mathrm{vol}(\delta \mathsf{M})}\cdot \mathbb{E}_{u\sim\partial \mathsf{M}}\left[\Im\left(f(x+i\delta u) \right)N(u) \right]. 
\end{equation}
\end{subequations}
for $N(u)$ a unit normal in $T^{\perp}_u\partial\mathsf{M}$.
\end{lemma}
\begin{proof}
As $\mathsf{M}\subset \mathbb{R}^n$ is a compact oriented manifold with boundary, we can appeal to the Divergence theorem~\cite[Theorem~16.32]{Lee2} (under the Euclidean metric), which states that for any smooth vector field $X$ on $\mathsf{M}$ one has 
\begin{equation}
\label{equ:div}
    \int_{\mathsf{M}}\mathrm{div}(X(v))dV(v) = \int_{\partial \mathsf{M}}\langle X(u),N(u)\rangle dV(u),
\end{equation}
for $N$ denoting the unit normal vector (field) along $\partial \mathsf{M}$. That is, $\mathbb{R}^n=T_p\partial\mathsf{M}\oplus T^{\perp}_p\partial\mathsf{M}$ for all $p\in \partial\mathsf{M}$ and $N(p)\in T^{\perp}_p\partial\mathsf{M}$.

Using the same reasoning as for example in~\cite{ref:JongeneelYueKuhnZO2021}, since one can select $X=f \cdot C$ for $C$ some constant vector field on $\mathsf{M}$, then, as $\mathrm{div}(C)=0$ and we can select $C$ to be aligned with any coordinate axis,~\eqref{equ:div} implies that
\begin{equation}
\label{equ:div2}
    \int_{\mathsf{M}}\nabla f(v) dV(v) = \int_{\partial \mathsf{M}}f(u)N(u)dV(u). 
\end{equation}
Note, $\nabla f(v)$ is well defined as $\mathsf{M}\subset \mathbb{R}^n$ is diffeomorphic to $\mathbb{B}^n$. 

Now we obtain the generalization of the result in~\cite{ref:JongeneelYueKuhnZO2021}, that is, by compactness, the Dominated Convergence theorem~\cite[Section~2.3]{ref:Folland}, the Divergence theorem~\eqref{equ:div} and the Cauchy-Riemann equations~\cite{ref:Krantz:Complex} we get
\begin{align*}
     \nabla_x \int_{\delta \mathsf{M}}\Re\left(f(x+iv)\right)\mathrm{d} V(v) {=}\int_{\delta \partial\mathsf{M}}\Im\left(f(x+iu)\right)N(u) \mathrm{d}V(u), 
\end{align*}
~\textit{e.g.}, see~\cite{ref:JongeneelYueKuhnZO2021} for more on this line of reasoning.
Then, due to the distributional assumption (uniformity), we write 
\begin{align*}
    f_{\delta,\mathsf{M}}(x) =& \mathbb{E}_{v\sim \mathsf{M}} \left[\Re\left(f(x+i\delta v)\right) \right] =\frac{1}{\mathrm{vol}(\delta \mathsf{M})}\int_{\delta \mathsf{M}}\Re\left(f(x+i v) \right) \mathrm{d}V(v),
\end{align*}
and similarly,
\begin{align*}
    \mathbb{E}_{u\sim \partial\mathsf{M}}\left[\Im\left(f(x+i\delta u) \right)N(u) \right] = \frac{1}{\mathrm{vol}(\delta \partial \mathsf{M})}\int_{\delta \partial \mathsf{M}} \Im\left(f(x+iu) \right)
    N(u)\mathrm{d}V(u). 
\end{align*}
Combining it all yields~\eqref{equ:gen:grad:f:delta}. 
\end{proof}

As $N(u)\in T^{\perp}_u\partial\mathsf{M}$ is a unit vector, the offset term in the variance~\eqref{equ:g:var} is minimized when we select $\mathsf{M}$ as
\begin{equation}
\label{equ:var:opt}
    \argmin_{\mathsf{M}\in \mathscr{M}}\frac{\mathrm{vol}(\delta\partial\mathsf{M})}{\mathrm{vol}(\delta \mathsf{M})},
\end{equation}
where $\mathscr{M}$ is the set of manifolds diffeomorphic to $\mathbb{B}^n$ and $\delta\in \mathbb{R}_{>0}$. 
To retrieve the optimizer, consider the isoperimetric inequality in $\mathbb{R}^n$~\cite{ref:osserman1978isoperimetric} which implies that $\mathsf{M}^{\star}=\mathbb{B}^n$ is optimal in the sense of~\eqref{equ:var:opt}. 

To get (the complex-step version of)~\cite[Corollary~6]{ref:hazan2014bandit} from Lemma~\ref{lem:gen:CR:grad}, let $\mathcal{E}^n_Q=\{x\in \mathbb{R}^n:\langle Q^{-1}x, x \rangle \leq 1\}$ for some $Q\in \mathcal{S}^n_{\succ 0}$. Now, $T_p\partial \mathcal{E}^n_Q=\{v\in \mathbb{R}^n:\langle Q^{-1}p,v\rangle =0\}$.
As $\mathcal{E}^n_Q=Q^{1/2}\mathbb{B}^n$ one can write
\begin{subequations}
\begin{equation}
 \label{equ:ellip:f}   
     f_{\delta,\mathcal{E}^n_Q}(x) = \mathbb{E}_{v\sim \mathcal{E}^n_Q} \left[f(x+i\delta v)\right] = \mathbb{E}_{v\sim \mathbb{B}^n} \left[f(x+i\delta Q^{1/2} v)\right]. 
\end{equation}
Via the rightmost term in~\eqref{equ:ellip:f} and the proof of Lemma~\ref{lem:gen:CR:grad} it follows immediately that
\begin{equation}
 \label{equ:ellip:grad}   
     \nabla f_{\delta,\mathcal{E}^n_Q}(x) = \mathbb{E}_{u\sim \mathbb{S}^{n-1}} \frac{n}{\delta}\left[f(x+i\delta Q^{1/2}u)Q^{-1/2}u\right].
\end{equation}
\end{subequations}
Equivalently, one can directly appeal to~\eqref{equ:gen:grad:f:delta}. However, here one needs to appeal to the isoperimetric ratio for ellipsoids~\cite{ref:rivin2007surface}. 

At last, we provide further comments on applicability. 
The complex-step derivative appears in a host of numerical applications, most notably, it is reported in~\cite[Page 44]{ref:cox2004software} that a value of $\delta=10^{-100}$ is successfully used in National Physical Laboratory software.
In the context of zeroth-order optimization, due to the complex-lifting,~\textit{i.e.}, we need $f(x+i\delta u)$, we cannot use immediately use physical measurement data, but we can work with any simulation-based system or data that admits a complex representation. A few areas of application are  
\begin{enumerate}[(i)]
    \item Simulation-based optimization,~\textit{e.g.,}~reinforcement learning and PDE-constrained optimization, see also~\cite{ref:wang2021improved} and Example~\ref{ex:PDE};
    \item Privacy-sensitive optimization,~\textit{e.g.,} the objective is known, but not to everybody; 
    \item Black-box objective,~\textit{e.g.,} $f(x)$ has been implemented in deprecated software, see also~\cite{ref:nesterov2017random}.
\end{enumerate}


\section{Strongly convex imaginary zeroth-order optimization}
\label{sec:strong:cvx:opt}

In this section we will utilize the imaginary gradient estimator $g_{\delta}$ as given by~\eqref{equ:grad:est:g} in the context of zeroth-order optimization algorithms. We will not focus on fully generic convex optimization problems as the flat parts of real-analytic convex functions must have measure zero~\cite{ref:Krantz:Complex,ref:JongeneelYueKuhnZO2021}. Hence, without too much loss of generality we omit convex functions which are not strongly convex\footnote{Future work will highlight the intimate relation between convex and strongly convex functions under the assumption that both are real analytic.}. See also~\cite{ref:kakadeduality} for more on strong-convexity in the context of generalization. 

In this section we relax some of the assumptions in~\cite{ref:JongeneelYueKuhnZO2021}, not only can we handle computational noise, the algorithms demand less knowledge of the problem compared to other work. This is possible by introducing a time-varying stepsize and a construction very much in line with~\cite{ref:akhavan2020exploiting}. In fact, recall from~\cite{ref:RakhlinSS12} that $\mu_k=\Theta(k^{-1})$ to allow for optimal rates. The edge our results have, however, over these existing works is that our sequence of smoothing parameters $\{\delta_k\}_{k\geq 1}$ is \textit{never} catastrophic.  

The generic algorithm for the unconstrained case is detailed in Algorithm~\ref{alg:convex_unconstrained}~(a),~\textit{i.e.}, for $\mathcal{X}= \mathcal{D}$.
\begin{algorithm}[t]
\begin{algorithmic}[1]
\STATE \textbf{Input:} initial iterate $x_1\in\mathcal{X}$, stepsizes $\{\mu_k \}_{k\ge 1}$, smoothing parameters $\{\delta_k \}_{k\ge 1}$.
  \FOR{$k = 1,2,\ldots, K$} 
  \STATE generate random $u_k \sim \mathbb{S}^{n-1}$
  
  \STATE obtain noisy estimate $g_{\delta_k}(x_k) = \frac{n}{\delta_k}\Im\left(f(x_k + i \delta_k u_k)\right)u_k+\frac{n}{\delta_k}\xi_k u_k$

  \STATE set $x_{k+1} = \Pi_{\mathcal{X}}\left( x_k - \mu_k\cdot g_{\delta_k}(x_k)\right)$
\ENDFOR
\end{algorithmic}
\caption{Imaginary zeroth-order optimization:\\ (a) \textit{unconstrained} $\mathcal{X}=\mathcal{D}$ and (b) \textit{constrained} $\mathcal{X}=\mathcal{K}\subset \mathcal{D}$.}
\label{alg:convex_unconstrained}
\end{algorithm}
Given a compact (possibly non-convex) set $\mathcal{K}\subset \mathcal{D}$, the algorithm for the constrained case is detailed in Algorithm~\ref{alg:convex_unconstrained}~(b)~\textit{i.e.}, for $\mathcal{X}= \mathcal{K}$. Here, $\Pi_{\mathcal{K}}:\mathcal{D}\subseteq \mathbb{R}^n\to \mathcal{K}$ denotes the \textbf{\textit{projection operator}}. 
  

Note, in our algorithms we will assume that we can sample in a small $\delta$-neighbourhood contained in $\mathcal{D}\setminus \mathcal{K}$. As a key application of the proposed framework is simulation-based optimization this is deemed justifiable. Having access to a projection operator $\Pi_{\mathcal{K}}$, we will assume nothing more than feasibility regarding the initial condition $x_1$. 

\subsection{Strong convexity}
In this part we consider the setting of $f\in C^{\omega}(\mathcal{D})$ being $\tau(f)$-\textit{\textbf{strongly convex}} over $\mathcal{D}$, \textit{i.e.,} there is some $\tau(f)>0$ such that 
\begin{equation}
    \label{equ:strong:convex}
    f(y)\geq f(x) + \langle \nabla f(x), y-x \rangle + \tfrac{1}{2}\tau(f)\|y-x\|_2^2,\quad \forall x,y\in \mathcal{D}. 
\end{equation}
In particular~\eqref{equ:strong:convex} implies that for $\mathcal{D}$ such that $x^{\star}\in \mathrm{int}(\mathcal{D})$ one has
\begin{equation}
\label{equ:strong:ineq}
    f(x)-f(x^{\star}) \geq \tfrac{1}{2}\tau(f)\|x-x^{\star}\|_2^2,\quad \forall x\in \mathcal{D}. 
\end{equation}
If additionally $f\in C^{\omega,1}_{L_1(f)}$, then by $\|\nabla f(x)\|_2^2\geq 2 \tau(f)(f(x)-f(x^{\star}))$ one has 
\begin{equation}
    \label{equ:Lip:strongly:convex}
    \tau(f)\|x-x^{\star}\|_2  \leq \|\nabla f(x)\|_2 \leq L_1(f)\|x-x^{\star}\|_2. 
\end{equation}

In contrast to~\cite{ref:JongeneelYueKuhnZO2021}, our algorithms ``\textit{only}'' demand knowledge of the strong-convexity parameter. In Section~\ref{sec:est:tau} we mention how one could estimate $\tau(f)$.

\subsection{Generic convergence rates}
\label{sec:strong:convex:generic}
As in~\cite{ref:akhavan2020exploiting}, we start with the constrained case.
\begin{theorem}[Convergence rate of Algorithm~\ref{alg:convex_unconstrained}~(b) with noise]
\label{thm:strong:convex:noise}
Let $f\in C^{\omega}(\mathcal{D})$ be a $\tau(f)$-strongly convex function satisfying Assumption~\ref{ass:holomorphic} for some $\bar{\delta}\in (0,1)$ and let $\mathcal{K}\subset \mathcal{D}$ be a compact convex set. Suppose that $f$ has a Lipschitz Hessian over $\mathcal{K}$, that is, ~\eqref{equ:Hess:Lip} holds for a non-zero constant $L_2(f)$.   
 Let $\{x_k\}_{k\ge 1}$ be the sequence of iterates generated by Algorithm~\ref{alg:convex_unconstrained}~(b) with stepsize $ \mu_k  = {2}/({\tau(f) k})$ and the sequence of smoothing parameters defined for all $k\geq 1$ by $\delta_k = \delta k^{-1/6}$ with $\delta\in (0,\kappa\bar{\delta}]$ for some $\kappa\in (0,1)$. 
Then, if the oracle satisfies Assumption~\ref{ass:complex:oracle}, the uniformly-averaged iterate $\bar{x}_K=K^{-1}\sum^K_{k=1}x_k$ achieves the optimization error 
\begin{equation*}
    \mathbb{E}[f(\bar{x}_K)-f(x^{\star})]
    \leq  \widetilde{O}\left(\frac{n^2}{\tau(f)} \delta^{-\tfrac{1}{3}}\sigma_{\xi} K^{-\tfrac{2}{3}}\right).
\end{equation*}
\end{theorem}
\begin{proof}
We mainly follow~\cite{ref:akhavan2020exploiting}. To that end, let $\sup_{x\in \mathcal{K}} \|\nabla f(x)\|_2\leq G$.
As $\mathcal{K}$ is convex and compact we have by the properties of the operator $\Pi_{\mathcal{K}}$ that $\|x_{k+1}-x^{\star}\|_2^2 \leq \|x_k - \mu_k g_{\delta_k}(x_k) - x^{\star}\|_2^2$. This can be written as conveniently as
\begin{equation}
\label{equ:g:inner:ineq:th1}
    \langle g_{\delta_k}(x_k), x_k - x^{\star}\rangle \leq \tfrac{1}{2\mu_k}\left( \|x_k - x^{\star}\|_2^2 - \|x_{k+1}-x^{\star}\|_2^2\right) + \tfrac{\mu_k}{2}\|g_{\delta_k}(x_k)\|_2^2.
\end{equation}
After reordering the standard strong $\tau(f)$-convexity expression, one obtains
\begin{equation}
\label{equ:f:ineq:tau:th1}
    f(x_k) -f(x^{\star}) \leq \langle \nabla f(x_k), x_k - x^{\star}\rangle - \tfrac{\tau(f)}{2}\|x_k - x^{\star}\|_2^2.
\end{equation}
Set $a_k=\|x_k-x^{\star}\|_2^2$, then, an application of the Cauchy-Schwarz inequality after combining~\eqref{equ:g:inner:ineq:th1} with~\eqref{equ:f:ineq:tau:th1} and taking the expectation over $u_k$ and $\xi_k$ conditioned on $x_k$ yields
\begin{align*}
\everymath={\displaystyle}
\begin{array}{lcl}
    \mathbb{E}[f(x_k)-f(x^{\star})|x_k] &\leq & \left\|\mathbb{E}[g_{\delta_k}(x_k)|x_k]-\nabla f(x_k) \right\|_2\|x_k-x^{\star}\|_2 + \tfrac{1}{2\mu_k}\mathbb{E}[a_k-a_{k+1}|x_k]\\
    && \tfrac{\mu_k}{2}\mathbb{E}[\|g_{\delta_k}(x_k)\|_2^2|x_k]-\tfrac{\tau(f)}{2}\mathbb{E}[a_k|x_k]\\
    &\overset{\eqref{equ:grad:delta:approx:error}}{\leq} & C_1 n\delta_k^2 \|x_k-x^{\star}\|_2 + \tfrac{1}{2\mu_k}\mathbb{E}[a_k-a_{k+1}|x_k]\\
    && \tfrac{\mu_k}{2}\mathbb{E}[\|g_{\delta_k}(x_k)\|_2^2|x_k]-\tfrac{\tau(f)}{2}\mathbb{E}[a_k|x_k],
    \end{array}
\end{align*}
for some $C_1>0$.
Now, use $ab\leq \tfrac{1}{2}(a^2+b^2)$, in particular $ab \leq \tfrac{1}{2}(\gamma a^2 + \gamma^{-1}b^2)$ for $\gamma\neq 0$, to construct 
\begin{equation*}
    n\delta_k^2 \|x_k-x^{\star}\|_2 \leq \tfrac{1}{2}\left(\tfrac{2C_1 }{\tau(f)} n^2 \delta_k^4 + \tfrac{\tau(f)}{2C_1}\|x_k-x^{\star}\|_2^2 \right).
\end{equation*}
Next, take the expectation over $x_k$ and let $r_k=\mathbb{E}[a_k]$ such that we can write
\begin{equation}
\label{equ:fk:th1}
    \mathbb{E}[f(x_k)-f(x^{\star})]\leq \tfrac{1}{2\mu_k}(r_k-r_{k+1})-\tfrac{\tau(f)}{4}r_k + \tfrac{1}{\tau(f)}C_1^2 n^2\delta_k^4 + \tfrac{\mu_k}{2}\mathbb{E}[\|g_{\delta_k}(x_k)\|_2^2]. 
\end{equation}
Summing~\eqref{equ:fk:th1} over $k$ yields
\begin{align*}
    \textstyle\sum^{K}_{k=1}\mathbb{E}[f(x_k)-f(x^{\star})]\leq \tfrac{1}{2}&\textstyle\sum^{K}_{k=1}\left( \tfrac{1}{\mu_k}(r_k-r_{k+1})-\tfrac{\tau(f)}{2}r_k\right)\\ 
    + &\textstyle\sum^{K}_{k=1}\left(\tfrac{1}{\tau(f)}C_1^2n^2\delta_k^4 + \tfrac{\mu_k}{2}\mathbb{E}[\|g_{\delta_k}(x_k)\|_2^2]\right). 
\end{align*}
As we selected $\mu_k = 2/(\tau(f) k)$ we can simplify the above by using the same reasoning as in~\cite{ref:akhavan2020exploiting}, that is
\begin{align*}
    \textstyle\sum^K_{k=1}\left( \tfrac{1}{\mu_k}(r_k-r_{k+1})-\tfrac{\tau(f)}{2}r_k\right) \leq r_1 \left(\tfrac{1}{\mu_1}-\tfrac{\tau(f)}{2}\right) + \textstyle\sum^K_{k=2}r_k \left(\tfrac{1}{\mu_k}-\tfrac{1}{\mu_{k-1}}-\tfrac{\tau(f)}{2}\right) =0. 
\end{align*}
Note that we rely on the $\tau(f)$-strong convexity. Using the observation from above and plugging in the stepsize $\mu_k$ elsewhere yields by~\eqref{equ:g:var}
\begin{align*}
    \textstyle\sum^{K}_{k=1}\mathbb{E}[f(x_k)-f(x^{\star})] &\leq  \tfrac{1}{\tau(f)}\textstyle\sum^{K}_{k=1}\left(C_1^2n^2\delta_k^4 + \tfrac{1}{k}\mathbb{E}[\|g_{\delta_k}(x_k)\|_2^2]\right)\\
    &\leq  \tfrac{n^2}{\tau(f)}\textstyle\sum^{K}_{k=1}\left(C_1^2\delta_k^4 + \tfrac{1}{k}\left[{C_2 \delta_k^4} + C_3\delta_k^2\|\nabla f(x_k)\|_2+ \tfrac{1}{n}\|\nabla f(x_k)\|_2^2+ \frac{1}{\delta_k^2}\sigma_{\xi}\right]\right),
\end{align*}
for some $C_2,C_3>0$. 
Now, minimizing over $\{\delta_k\}_k$ is possible but yields smoothing parameters as a function of unknown constants. Instead, we retain the ``\textit{optimal}'' root\footnote{Let $a,b\in \mathbb{R}_{>0}$, then, see that $(b/(2a))^{\tfrac{1}{6}}=\argmin_{\delta\in \mathbb{R}_{\geq 0}}\{a\delta^4+b\tfrac{1}{\delta^2}\}$.} and propose
\begin{equation*}
    \widetilde{\delta}_k = \left( \frac{\alpha \sigma_{\xi}}{k } \right)^{\tfrac{1}{6}},
\end{equation*}
for some $\alpha\in (0,1)$ to be specified. 
Using this smoothing parameter sequence, that is, $\delta_k= \widetilde{\delta}_k$, together with $\sum^K_{k=1}k^{-1}\leq 1+\log(K)$ (Lemma~\ref{lem:log}) yields  
\begin{align*}
\everymath={\displaystyle}
\begin{array}{lcl}
    \textstyle\sum^{K}_{k=1}\mathbb{E}[f(x_k)-f(x^{\star})]
    &\leq&  \tfrac{n^2}{\tau(f)}\textstyle\sum^{K}_{k=1}\left(C_1^2\left( \frac{\alpha \sigma_{\xi}}{k } \right)^{\tfrac{2}{3}} + \tfrac{1}{k}\left[{C_2 \left( \frac{\alpha \sigma_{\xi}}{k } \right)^{\tfrac{2}{3}}} + \left( \frac{\alpha \sigma_{\xi}}{k } \right)^{-\tfrac{1}{3}}\sigma_{\xi}\right]\right)\\
    && + \tfrac{n}{\tau(f)}G^2 (1+\log(K))+\textstyle\tfrac{n^2}{\tau(f)}G C_3\sum^K_{k=1}\tfrac{1}{k}\left(\frac{\alpha \sigma_{\xi}}{k}\right)^{\tfrac{1}{3}}\\
    &=&  \tfrac{n^2}{\tau(f)}\textstyle\sum^{K}_{k=1}\left(C_1^2\left( \frac{\alpha \sigma_{\xi}}{k} \right)^{\tfrac{2}{3}} + \tfrac{1}{k}\left[{C_2\left( \frac{\alpha \sigma_{\xi}}{k} \right)^{\tfrac{2}{3}}} + k^{\tfrac{1}{3}}\sigma_{\xi}^{\tfrac{2}{3}}\alpha^{-\tfrac{1}{3}}\right]\right)\\
    && + \tfrac{n}{\tau(f)}G^2 (1+\log(K))+\textstyle\tfrac{n^2}{\tau(f)}G C_3\sum^K_{k=1}k^{-\tfrac{2}{3}}(\alpha \sigma_{\xi})^{\tfrac{1}{3}}\\
    &\leq&  \tfrac{n^2}{\tau(f)}\textstyle\sum^{K}_{k=1}C_4 k^{-\tfrac{2}{3}} \sigma_{\xi}^{\tfrac{2}{3}}\alpha^{-\tfrac{1}{3}} \\
    && + \tfrac{n}{\tau(f)}G^2 (1+\log(K))+\textstyle\tfrac{n^2}{\tau(f)}G C_3\sum^K_{k=1}k^{-\tfrac{2}{3}}(\alpha \sigma_{\xi})^{\tfrac{1}{3}}. 
    \end{array}
\end{align*}
Now, as $\sum^K_{k=1}k^{-\tfrac{2}{3}}\leq 3 K^{\tfrac{1}{3}}$ (Lemma~\ref{lem:fractional}) we can continue and write 

\begin{align*}
\everymath={\displaystyle}
\begin{array}{lcl}
    \textstyle\sum^{K}_{k=1}\mathbb{E}[f(x_k)-f(x^{\star})]
    &\leq&  \tfrac{n^2}{\tau(f)}C_5 K^{\tfrac{1}{3}} \sigma_{\xi}^{\tfrac{2}{3}}\alpha^{-\tfrac{1}{3}}  + \tfrac{n}{\tau(f)}G^2 (1+\log(K))+\textstyle\tfrac{n^2}{\tau(f)}G C_6 K^{\tfrac{1}{3}}(\alpha \sigma_{\xi})^{\tfrac{1}{3}}. 
    \end{array}
\end{align*}
and as such we obtain the optimization error
\begin{align*}
\everymath={\displaystyle}
\begin{array}{lcl}
    \mathbb{E}[f(\bar{x}_K)-f(x^{\star})]
    &\leq&  \tfrac{n^2}{\tau(f)}C_5 K^{-\tfrac{2}{3}} \sigma_{\xi}^{\tfrac{2}{3}}\alpha^{-\tfrac{1}{3}}  + \tfrac{n}{\tau(f)}G^2K^{-1} (1+\log(K))+\textstyle\tfrac{n^2}{\tau(f)}G C_6 K^{-\tfrac{2}{3}}(\alpha \sigma_{\xi})^{\tfrac{1}{3}}.
    \end{array}
\end{align*}
As $\alpha\in (0,1)$ was arbitrary, we can set $\delta=\alpha \sigma_{\xi}$ such that $\delta_k=\delta k^{-\tfrac{1}{6}}$ for some $\delta\in (0,\bar{\delta})$.
\end{proof}
The edge Theorem~\ref{thm:strong:convex:noise} has over existing work is that the requested sequence $\{\delta_k\}_{k\geq 1}$ \textit{can} always be safely implemented. With respect to optimality, we highlight a general method to pass from $\widetilde{O}(\cdot)$ to $O(\cdot)$ complexities.

\begin{remark}[Removing the logarithmic term]
\label{rem:suffix:averaging}
\upshape{
One can appeal to \textit{$\alpha$-suffix averaging} as proposed in~\cite{ref:RakhlinSS12} to remove the logarithmic term. This is achieved by averaged estimates of the form $\widetilde{x}_K=\tfrac{2}{K}\sum^K_{k=K/2+1}x_k$ and follows from $\sum^T_{t=(1-\alpha)T+1}\tfrac{1}{t}\leq \log(1/(1-\alpha))$ for $\alpha\in (0,1)$ such that $\alpha T, (1-\alpha) T \in \mathbb{Z}$. As the implementation of $\widetilde{x}_K$ is not always easier or more efficient than $\bar{x}_K$, the uniformly-averaged estimator remains competitive despite the slower rate. 
}
\end{remark}

Next we consider the unconstrained case. Here, we cannot appeal to an uniform bound on $\nabla f(x)$. Instead, we use the idea from~\cite[Theorem~3.2]{ref:akhavan2020exploiting} and bound a subset of iterates before strong-convexity kicks in. In practise, when $\tau(f)$ is small, the first few stepsizes will be relatively large and can lead to overflow. In some sense one could interpret this as some restarting mechanism. 

\begin{theorem}[Convergence rate of Algorithm~\ref{alg:convex_unconstrained}~(a) with noise]
\label{thm:strong:convex::unconstrained:noise}
Let $f\in C^{\omega}(\mathcal{D})$ be a $\tau(f)$-strongly convex function satisfying Assumption~\ref{ass:holomorphic} for some $\bar{\delta}\in (0,1)$ with $x^{\star}\in \mathrm{int}(\mathcal{D})$. Suppose that $f$ has a Lipschitz gradient and Hessian, that is, ~\eqref{equ:grad:Lipschitz} and~\eqref{equ:Hess:Lip} hold, for non-zero constants $L_1(f)$ and $L_2(f)$, respectively.   
 Let $\{x_k\}_{k\ge 1}$ be the sequence of iterates generated by Algorithm~\ref{alg:convex_unconstrained}~(a) for 
 \begin{align*}
\everymath={\displaystyle}
 \begin{array}{llll}
 &\mu_k = \tfrac{1}{\tau(f) K},\quad  &\delta_k = \delta K^{-\tfrac{1}{6}},\quad &k=1,\dots,K_0,   \\
   &\mu_k = \tfrac{2}{\tau(f) k},\quad  &\delta_k = \delta k^{-\tfrac{1}{6}},\quad &k=K_0+1,\dots,K,
 \end{array}
 \end{align*}
 with $K_0 = \floor*{\tfrac{8 n^2L_1(f)^2}{\tau(f)^2}}$ and $\delta\in (0,\kappa\bar{\delta}]$ for some $\kappa\in (0,1)$. Then, if the oracle satisfies Assumption~\ref{ass:complex:oracle} and $K\geq 2K_0$ we incur for $\bar{x}_{K_0,K}=\tfrac{1}{K-K_0}\sum^K_{k=K_0+1}x_k$
 the optimization error 
\begin{equation}
\label{equ:unconstrained:rate}
\begin{aligned}
    \mathbb{E}[f(\bar{x}_{K_0,K})-f(x^{\star})]
    &\leq O\left(\frac{n^2L_1(f)^2}{\tau(f) }\|x_1-x^{\star}\|_2^2K^{-1}\right) + O\left( \frac{n^2\sigma_{\xi}}{\tau(f)\delta^2} K^{-\tfrac{2}{3}}\right). 
\end{aligned}
\end{equation}
\end{theorem}
\begin{proof}
The proof will be similar to that of~\cite[Theorem~3.2]{ref:akhavan2020exploiting}. 
Again, set $a_k=\|x_k-x^{\star}\|_2^2$, then, as in the proof of Theorem~\ref{thm:strong:convex:noise} 
\begin{align*}
\everymath={\displaystyle}
\begin{array}{lcl}
    \mathbb{E}[f(x_k)-f(x^{\star})|x_k] &\leq & \left\|\mathbb{E}[g_{\delta_k}(x_k)|x_k]-\nabla f(x_k) \right\|_2\|x_k-x^{\star}\|_2 + \tfrac{1}{2\mu_k}\mathbb{E}[a_k-a_{k+1}|x_k]\\
    && \tfrac{\mu_k}{2}\mathbb{E}[\|g_{\delta_k}(x_k)\|_2^2|x_k]-\tfrac{\tau(f)}{2}\mathbb{E}[a_k|x_k]\\
    &\overset{\eqref{equ:grad:delta:approx:error}}{\leq} & C_1 n\delta_k^2 \|x_k-x^{\star}\|_2 + \tfrac{1}{2\mu_k}\mathbb{E}[a_k-a_{k+1}|x_k]\\
    && \tfrac{\mu_k}{2}\mathbb{E}[\|g_{\delta_k}(x_k)\|_2^2|x_k]-\tfrac{\tau(f)}{2}\mathbb{E}[a_k|x_k].
    \end{array}
\end{align*}
Now, use $ab\leq \tfrac{1}{2}(a^2+b^2)$ together with $\tau(f)$-strong convexity,~\textit{i.e.},~\eqref{equ:strong:ineq}, to construct 
\begin{equation*}
    n\delta_k^2 \|x_k-x^{\star}\|_2 \leq \tfrac{1}{2}\left(\tfrac{2C_1 }{\tau(f)} n^2 \delta_k^4 + \tfrac{\tau(f)}{2C_1 }\|x_k-x^{\star}\|_2^2 \right)\leq \tfrac{C_1 }{\tau(f)} n^2 \delta_k^4 + \tfrac{1}{2C_1 }(f(x_k)-f(x^{\star})).
\end{equation*}
Next, let $r_k=\mathbb{E}[a_k]$ such that by $\|\nabla f(x_k)\|_2^2\leq L_1(f)^2 \|x_k-x^{\star}\|_2^2$ we can write
\begin{equation}
\label{equ:fk}
\begin{aligned}
\everymath={\displaystyle}
\begin{array}{lcl}
    \mathbb{E}[f(x_k)-f(x^{\star})]&\leq& \tfrac{1}{\mu_k}(r_k-r_{k+1})-\tau(f)r_k + \tfrac{2}{\tau(f)}C_1^2n^2\delta_k^4 + {\mu_k}\mathbb{E}[\|g_{\delta_k}(x_k)\|_2^2]\\
    &\overset{\eqref{equ:g:var}}{\leq}& \tfrac{1}{\mu_k}(r_k-r_{k+1})-{\tau(f)}r_k + \tfrac{2}{\tau(f)}C_1^2n^2\delta_k^4 \\ 
    && + {\mu_k}\left({C_2 n^2\delta_k^4} + C_3 n^2\delta_k^4+  2 n^2L_1(f)^2 r_k+ \tfrac{n^2}{\delta_k^2}\sigma_{\xi} \right)
\end{array}
\end{aligned}
\end{equation}
Where in the last step we used 
\begin{equation*}
    \delta_k^2 \|x_k-x^{\star}\|_2 \leq \tfrac{1}{2}\left(\tfrac{ C_3}{2 L_1(f)} \delta_k^4+ \tfrac{2 L_1(f)}{ C_3}\|x_k-x^{\star}\|_2^2\right) 
\end{equation*}
to rewrite~\eqref{equ:g:var}. 

Now we use the step- and smoothingsize for $k=1,\dots,K_0$, that is, $\mu_k= 1/(\tau(f) K)$, $\delta_k = \delta K^{-\tfrac{1}{6}}$, and observe that 

\begin{align*}
    r_{k+1}&\leq r_k-\tau(f)\mu_kr_k + \tfrac{2\mu_k}{\tau(f)}C_1^2n^2\delta_k^4 + \mu_k^2\left({C_2 n^2\delta_k^4} + C_3 n^2\delta_k^4+ 2n^2L_1(f)^2 r_k+ \tfrac{n^2}{\delta_k^2}\sigma_{\xi} \right)\\
    &= \left(1-\tfrac{1}{K}+\tfrac{2 L_1(f)^2}{(\tau(f)K)^2}n^2\right)r_k + \nu_K = a_K r_k + \nu_K
\end{align*}
for $a_k$ as between brackets and $\nu_K$ defined as 
\begin{align*}
    \nu_K =& \tfrac{2}{\tau(f)^2 K}C_1^2n^2\delta_k^4 + \tfrac{1}{(\tau(f) K)^2}\left({C_2 n^2\delta_k^4}+ C_3 n^2\delta_k^4+ \tfrac{n^2}{\delta_k^2}\sigma_{\xi} \right)\\
    \leq& \tfrac{1}{\tau(f)^2 K}\left( n^2 \delta^4 C_4 + \tfrac{n^2}{\delta^2}\sigma_{\xi}\right)K^{-\tfrac{2}{3}}.
\end{align*}
We now proceed with bounding $r_{K_0+1}$. As in~\cite{ref:akhavan2020exploiting}, set 
\begin{equation*}
    q_K = 1+\tfrac{2 L_1(f)^2}{(\tau(f)K)^2}n^2
\end{equation*}
by iterating over $r_k$ it follows from a geometric series argument that
\begin{align*}
    r_{K_0+1}\leq a_K^{K_0}r_1 + \textstyle\sum^{K_0-1}_{i=0}a_K^i \nu_K \leq \left( r_1 + \tfrac{(\tau(f)K)^2}{2 L_1(f)^2n^2}\nu_K \right)q_K^{K_0}. 
\end{align*}
Now for $\floor*{\cdot}$ being the floor function, let $K_0$ be as in the theorem. 
Then, as $\log(1+x)\leq x$, on $\mathbb{R}_{\geq 0}$, one has
\begin{align*}
    q_K^{K_0} =& \mathrm{exp}\left(K_0 \log\left(1+\tfrac{2L_1(f)^2}{(\tau(f)K)^2}n^2 \right) \right)
    \leq  \mathrm{exp}\left(\tfrac{8n^2L_1(f)^2}{\tau(f)^2} \log\left(1+\tfrac{2L_1(f)^2}{(\tau(f)K)^2}n^2 \right) \right)
    \leq  \mathrm{exp}\left(\tfrac{16 n^4 L_1(f)^4}{\tau(f)^4 K^2} \right).
\end{align*}
Fix any $\bar{C}\in (0,\tfrac{1}{32})$, when 
\begin{equation*}
    K= \sqrt{\tfrac{8 n^4 L_1(f)^4}{\tau(f)^4\bar{C}}}
\end{equation*}
then $K\geq 2K_0$ and $q^{K_0}_K\leq e^{\bar{C}}=C_5$. 
As such,
\begin{align*}
    r_{K_0+1}\leq& \left( r_1 + \tfrac{(\tau(f)K)^2}{2L_1(f)^2n^2}\nu_K \right)C_5\\ 
    \leq& \left( r_1 + \tfrac{(\tau(f)K)^2}{L_1(f)^2n^2}\tfrac{1}{\tau(f)^2 K}\left( n^2 \delta^4 C_4 + \tfrac{n^2}{\delta^2}\sigma_{\xi}\right)K^{-\tfrac{2}{3}} \right)C_5\\ 
    =& \left( r_1 + \tfrac{1}{L_1(f)^2n^2}\left( n^2 \delta^4C_4 + \tfrac{n^2}{\delta^2}\sigma_{\xi}\right)K^{\tfrac{1}{3}} \right)C_5.
\end{align*}
Now we return to our normal step- and smoothingsizes, that is $\mu_k=2/(\tau(f) k)$, $\delta_k =\delta  k^{-\tfrac{1}{6}}$, for $k\geq K_0+1$. By plugging this into~\eqref{equ:fk} we get
\begin{equation*}
\begin{aligned}
    (K-K_0)\mathbb{E}[f(\bar{x}_{K_0,K})-f(x^{\star})]
    &\leq \textstyle \sum^K_{k=K_0+1}\tfrac{\tau(f)k}{2}(r_k-r_{k+1})-{\tau(f)}r_k + \tfrac{4}{\tau(f)k}n^2L_1(f)^2 r_k \\ 
    & +\textstyle\sum^K_{k=K_0+1} \tfrac{2}{\tau(f)k}C_6 n^2\delta^4 k^{-\tfrac{2}{3}} \\
    & +\textstyle\sum^K_{k=K_0+1} \tfrac{2}{\tau(f)k}\tfrac{n^2}{\delta^2}k^{\tfrac{1}{3}}\sigma_{\xi} 
\end{aligned}
\end{equation*}
By construction of $K_0$ we have that for $k\geq K_0+1$, $\tau(f)/2\geq (4n^2L_1(f)^2)/(\tau(f) k)$. Hence 
\begin{equation*}
\begin{aligned}
    (K-K_0)\mathbb{E}[f(\bar{x}_{K_0,K})-f(x^{\star})]
    &\leq \tfrac{\tau(f)}{2}\left(\textstyle\sum^K_{k=K_0+1}k(r_k-r_{k+1})- r_k\right) + U_{K_0,K}.
\end{aligned}
\end{equation*}
where by Lemma~\ref{lem:fractional}
\begin{align*}
    U_{K_0,K} =& \textstyle\sum^K_{k=K_0+1} \tfrac{2}{\tau(f)k}C_6 n^2\delta^4 k^{-\tfrac{2}{3}}  + \tfrac{2}{\tau(f)k}\tfrac{n^2}{\delta^2}k^{\tfrac{1}{3}}\sigma_{\xi} \\
    \leq& \tfrac{2n^2}{\tau(f)}\left(C_6\delta^4+\tfrac{1}{\delta^2}\sigma_{\xi} \right)\textstyle\sum^K_{k=K_0+1}k^{-\tfrac{2}{3}}\\
    \leq& \tfrac{2n^2}{\tau(f)}\left(C_6\delta^4+\tfrac{1}{\delta^2}\sigma_{\xi} \right)3 K^{\tfrac{1}{3}}. 
\end{align*}
As demonstrated in~\cite{ref:akhavan2020exploiting} (below), one can now construct the bound $\textstyle\sum^K_{k=K_0+1} k(r_k - r_{k+1})-r_k \leq K_0 r_{K_0+1}$ where the last term is exactly the term we could bound before. In combination with the bound on $K_0$ itself, we find that 
\begin{equation*}
\begin{aligned}
    (K-K_0)\mathbb{E}[f(\bar{x}_{K_0,K})-f(x^{\star})]
    &\leq \tfrac{\tau(f)}{2}\tfrac{8n^2L_1(f)^2}{\tau(f)^2}\left( r_1 + \tfrac{1}{L_1(f)^2n^2}\left( n^2 \delta^4 C_4 + \tfrac{n^2}{\delta^2}\sigma_{\xi}\right)K^{\tfrac{1}{3}} \right)C_5\\
    &+ \tfrac{2n^2}{\tau(f)}\left(C_6\delta^4+\tfrac{1}{\delta^2}\sigma_{\xi} \right)3 K^{\tfrac{1}{3}}.
\end{aligned}
\end{equation*}
By our selection of $\bar{C}$ we have that $K\geq 2 K_0$ and as such 
\begin{equation*}
\begin{aligned}
    \mathbb{E}[f(\bar{x}_{K_0,K})-f(x^{\star})]
    &\leq \tfrac{8n^2L_1(f)^2}{\tau(f) K}\left( r_1 + \tfrac{1}{L_1(f)^2}\left( \delta^4 C_3 + \tfrac{1}{\delta^2}\sigma_{\xi}\right)K^{\tfrac{1}{3}} \right)C_5\\
    &+ \tfrac{4n^2}{\tau(f)K}\left(C_6\delta^4+\tfrac{1}{\delta^2}\sigma_{\xi} \right)3 K^{\tfrac{1}{3}}.
\end{aligned}
\end{equation*}
Now, reordering terms yields~\eqref{equ:unconstrained:rate}. 
\end{proof}

\begin{remark}[On unconstrained anytime algorithms]
\label{rem:anytime}
\upshape{
The unconstrained algorithms (Theorem~\ref{thm:strong:convex::unconstrained:noise} and Theorem~\ref{thm:strong:convex::unconstrained:noise:quadratic}) require the user to pre-define the full length of the algorithm as the stepsize depends on $K$. One can mitigate this by shifting the dependence on $k$, \textit{e.g.}, by using $\mu_k = 1/\big(\tau(f)(k+2K_0)\big)$. Although the rate~\eqref{equ:unconstrained:rate} remains unaffected, this does come at the cost of potentially sacrificing progress in the first $K_0$ steps of the algorithm. A detailled study is left for future work. 
}
\end{remark}

\subsection{Optimal convergence rates}
\label{sec:strong:convex:opt:opt}
Now we consider the special case of $f$ being quadratic.
Here we improve upon the previous section due to exploitation of the quadratic nature of $f$, that is, by using $\Im(f(x+i\delta u))=\delta \langle \nabla f(x), u \rangle$ for any $\delta>0$. 

 Better yet, we see that for quadratic functions we incur optimal regret. Optimality can be shown along the lines of~\cite{shamir2013complexity}, or along the lines of~\cite{ref:agarwal2009information} after observing that in the quadratic case the gradient estimator $g_{\delta}(x)$ becomes an unbiased estimator for $\nabla f(x)$. The test function used in~\cite{ref:akhavan2020exploiting} is smooth but unfortunately not analytic\footnote{Section~\ref{sec:smooth} highlights that this might not be an obstruction.}. We start by providing the bound from below. 

 \begin{theorem}[Bound from below]
\label{thm:quad:below}
Any possibly randomized zeroth-order algorithm of fixed length $K$, applying the estimator~\eqref{equ:grad:est:g} under Assumption~\ref{ass:complex:oracle}, cannot achieve a rate faster than
\begin{equation*}
    \Omega\left(\frac{n^2}{\tau(f) K} \right),
\end{equation*}
uniformly over all $\tau(f)$-strongly convex quadratic (real-analytic) functions. 
\end{theorem}
\begin{proof}
We largely follow~\cite[Theorem~3]{shamir2013complexity}, but for the sake of completeness we highlight the main arguments. 

Recall that based on $x_1,x_2,\dots,x_K$, in particular the function evaluations at those points, we compute some point ${x}'_K$ (this could be a non-uniform average estimator). In our case the function queries correspond to $v_k  = \Im (f(x_k+i\delta u))$ for some choice of $\delta>0$, $u\in \mathbb{S}^{n-1}$ and with the possibility of being corrupted by additive noise $\xi$.

Now, consider the $C^{\omega}$ function over $\mathbb{R}^n$
\begin{equation}
    \label{equ:fz}
    f_z(x) = \tfrac{\tau}{2}\|x\|_2^2 - \langle z, x \rangle. 
\end{equation}
The unique minimizer of $f_z(x)$ is given by $x^{\star}=\tfrac{1}{\tau}z$. 
Moreover, assume $z$ is drawn uniformly from $\{-\nu,\nu\}^n$ for some $\nu$ that will be specified later. 
It follows from the strong $\tau$-convexity of~\eqref{equ:fz} that $f_z(x)-f_z(x^{\star})\geq \tfrac{\tau}{2}\|x-\tfrac{1}{\tau}z\|_2^2$. As such, for any randomized strategy
\begin{equation*}
    \mathbb{E}[f_z(x'_K)-f_z(x^{\star})]\geq \tfrac{\tau}{2}\mathbb{E}[\|x'_K-\tfrac{1}{\tau}z\|_2^2]=\tfrac{\tau}{2}\mathbb{E}\left[\textstyle\sum^n_{i=1}(x'_{i,K}-\tfrac{1}{\tau}z_i)^2 \right]\geq\tfrac{ \nu^2}{2\tau}\mathbb{E}\left[\mathbbold{1}_{x'_{i,K}z_i<0} \right],  
\end{equation*}
where the expectation is taken over the quadratic functions of the form~\eqref{equ:fz}. 
This means that we can construct a bound from below if we can get a grip on the signs of each $z_i$. To that end, we follow the proof of~\cite[Theorem~3]{shamir2013complexity}. The idea is to consider deterministic strategies that have only access to a sequence of function evaluations. The KL-divergence will allow for relating these function evaluations and the sign of $z_i$. 

The key difference with respect to~\cite{shamir2013complexity}, however, is the estimator. Given some point $x_k$, our function evaluation $v_k$ is of the form $v_k = \Im(f(x_k+i\delta u))+\xi$ for some $\delta>0$, $u\in \mathbb{S}^{n-1}$ and noise realization $\xi$. Now observe that $\Im(f_z(x+i\delta u))=\delta (\tau \langle x,u\rangle - \langle z,u \rangle)$. Hence, conditioning on $z_i>0$ we get
\begin{subequations}
\begin{equation*}
    v_k = \delta\left(\tau \langle x_k, u \rangle - \textstyle\sum_{j\neq i}z_ju_j\right)-\nu u_i + \xi
\end{equation*}
whereas conditioning on $z_i<0$ yields
\begin{equation*}
    v_k = \delta\left(\tau \langle x_k, u \rangle - \textstyle\sum_{j\neq i}z_ju_j\right)+\nu u_i + \xi.
\end{equation*}
\end{subequations}
Under the assumption that the noise is Gaussian one can now bound the KL-divergence by $(2\nu u_i)^2/(2\sigma_{\xi})$,~\textit{e.g.}, see~\cite[Lemma~5]{shamir2013complexity}. 
Using the fact that $u\in \mathbb{S}^{n-1}$ one can now exploit~\cite[Lemma~4]{shamir2013complexity} and show that
\begin{equation*}
    \tfrac{\nu^2}{2\tau}\mathbb{E}[\mathbbold{1}_{x'_{i,K}z_i<0}]\geq \tfrac{ n\nu^2}{4\tau}\left(1-\sqrt{\tfrac{2\nu^2 K}{n\sigma_{\xi}}} \right). 
\end{equation*}
As such, selecting $\nu = \sqrt{({n\sigma_{\xi}})/({4K})}$
yields the desired result. 
\end{proof}
 In the light of Theorem~\ref{thm:quad:below} and Remark~\ref{rem:suffix:averaging}, the following algorithms are \textit{rate} optimal. More specifically, one can show that the dependence on $\sigma_{\xi}$ is also optimal. Note that for quadratic functions we should not simply appeal to Theorem~\ref{thm:strong:convex:noise} as that proof relies on $L_2(f)>0$. 
 
 \begin{theorem}[Convergence rate of Algorithm~\ref{alg:convex_unconstrained}~(b) with noise, $f$ being quadratic]
\label{thm:strong:convex:noise:quadratic}
Let $f\in C^{\omega}(\mathcal{D})$ be a $\tau(f)$-strongly convex function satisfying Assumption~\ref{ass:holomorphic} for some $\bar{\delta}\in (0,1)$ and let $\mathcal{K}\subset \mathcal{D}$ be a compact convex set. Suppose that $f$ has a constant Hessian over $\mathcal{K}$, that is,~\eqref{equ:Hess:Lip} holds with $L_2(f)=0$.   
 Let $\{x_k\}_{k\ge 1}$ be the sequence of iterates generated by Algorithm~\ref{alg:convex_unconstrained}~(b) with stepsize $ \mu_k  = {2}/({\tau(f) k})$ and the sequence of smoothing parameters defined for all $k\geq 1$ by $\delta_k = \delta $ with $\delta\in (0,\kappa\bar{\delta}]$ for some $\kappa\in (0,1)$. 
Then, if the oracle satisfies Assumption~\ref{ass:complex:oracle},  the uniformly-averaged iterate $\bar{x}_K=K^{-1}\sum^K_{k=1}x_k$ achieves the optimization error 
\begin{equation*}
    \mathbb{E}[f(\bar{x}_K)-f(x^{\star})]
    \leq  {O}\left(\frac{n}{\tau(f)} K^{-1}\right)+\widetilde{O}\left(\frac{n^2\sigma_{\xi}}{\tau(f)\delta^2}K^{-1}\right).
\end{equation*}
\end{theorem}
\begin{proof}
We can mainly follow the proof of Theorem~\ref{thm:strong:convex:noise}, which relies itself largely on~\cite{ref:akhavan2020exploiting}. To that end, let again $\sup_{x\in \mathcal{K}} \|\nabla f(x)\|_2\leq G$ and set $a_k=\|x_k-x^{\star}\|_2^2$ such that 
\begin{align*}
\everymath={\displaystyle}
\begin{array}{lcl}
    \mathbb{E}[f(x_k)-f(x^{\star})|x_k] &\leq & \left\|\mathbb{E}[g_{\delta_k}(x_k)|x_k]-\nabla f(x_k) \right\|_2\|x_k-x^{\star}\|_2 + \tfrac{1}{2\mu_k}\mathbb{E}[a_k-a_{k+1}|x_k]\\
    && \tfrac{\mu_k}{2}\mathbb{E}[\|g_{\delta_k}(x_k)\|_2^2|x_k]-\tfrac{\tau(f)}{2}\mathbb{E}[a_k|x_k]\\
    &\overset{\eqref{equ:grad:delta:approx:error}}{\leq} &  \tfrac{1}{2\mu_k}\mathbb{E}[a_k-a_{k+1}|x_k]+ \tfrac{\mu_k}{2}\mathbb{E}[\|g_{\delta_k}(x_k)\|_2^2|x_k]-\tfrac{\tau(f)}{2}\mathbb{E}[a_k|x_k].
    \end{array}
\end{align*}
Next, let $r_k=\mathbb{E}[a_k]$ such that we can write
\begin{equation}
\label{equ:fk:Q}
    \mathbb{E}[f(x_k)-f(x^{\star})]\leq \tfrac{1}{2\mu_k}(r_k-r_{k+1})-\tfrac{\tau(f)}{2}r_k  + \tfrac{\mu_k}{2}\mathbb{E}[\|g_{\delta_k}(x_k)\|_2^2]. 
\end{equation}
To allow for an identical stepsize as before, we replace $-\tau(f)/2$ with $-\tau(f)/4$.
Summing~\eqref{equ:fk:Q} over $k$ yields
\begin{align*}
    \textstyle\sum^{K}_{k=1}\mathbb{E}[f(x_k)-f(x^{\star})]\leq& \tfrac{1}{2}\textstyle\sum^{K}_{k=1}\left( \tfrac{1}{\mu_k}(r_k-r_{k+1})-\tfrac{\tau(f)}{2}r_k\right)+ \textstyle\sum^{K}_{k=1} \tfrac{\mu_k}{2}\mathbb{E}[\|g_{\delta_k}(x_k)\|_2^2]. 
\end{align*}
As we selected $\mu_k = 2/(\tau(f) k)$ we can simplify the above by using the same reasoning as in~\cite{ref:akhavan2020exploiting}, that is
\begin{align*}
    \textstyle\sum^K_{k=1}\left( \tfrac{1}{\mu_k}(r_k-r_{k+1})-\tfrac{\tau(f)}{2}r_k\right) \leq r_1 \left(\tfrac{1}{\mu_1}-\tfrac{\tau(f)}{2}\right) + \textstyle\sum^K_{k=2}r_k \left(\tfrac{1}{\mu_k}-\tfrac{1}{\mu_{k-1}}-\tfrac{\tau(f)}{2}\right) =0. 
\end{align*}
Indeed, without the scaling of $\tau(f)$ our stepsize would have been $\mu_k = 1/(\tau(f) k)$.
Note that we rely on the $\tau(f)$-strong convexity. Using the observation from above and plugging in the stepsize $\mu_k$ elsewhere yields
\begin{align*}
\everymath={\displaystyle}
\begin{array}{lcl}
    \textstyle\sum^{K}_{k=1}\mathbb{E}[f(x_k)-f(x^{\star})] &\leq&  \tfrac{1}{\tau(f)}\textstyle\sum^{K}_{k=1} \tfrac{1}{k}\mathbb{E}[\|g_{\delta_k}(x_k)\|_2^2]\\
    &\leq&  \tfrac{n^2}{\tau(f)}\textstyle\sum^{K}_{k=1}\tfrac{1}{k}\left[ \tfrac{1}{n}\|\nabla f(x_k)\|_2^2+ \frac{1}{\delta_k^2}\sigma_{\xi}\right].
    \end{array}
\end{align*}
Now, minimizing over $\{\delta_k\}_k$ clearly yields a desire to pick a larger and fixed $\delta$ \textit{cf}.~Theorem~\ref{thm:strong:convex:noise}. Combining this with the bound on $\nabla f(x)$ yields by~\eqref{equ:log:bound} 
\begin{align*}
\everymath={\displaystyle}
\begin{array}{lcl}
    \textstyle\sum^{K}_{k=1}\mathbb{E}[f(x_k)-f(x^{\star})]
    &\leq&  \tfrac{n^2}{\tau(f)}\left[ \tfrac{1}{n}G^2+ \tfrac{1}{\delta^2}\sigma_{\xi}\right](\log(K)+1),
    \end{array}
\end{align*}
as such we obtain the optimization error $\mathbb{E}[f(\bar{x}_k)-f(x^{\star})]
    \leq  \tfrac{n^2}{\tau(f)K}\left[ \tfrac{1}{n}G^2+ \tfrac{1}{\delta^2}\sigma_{\xi}\right](\log(K)+1)$.
\end{proof}

\begin{example}[Numerical strongly-convex optimization]
\label{ex:strongcvx}
\upshape{
Here we exemplify what can go wrong and how the proposed complex-step method handles this. 
Consider for $n=10^3$ the problem of solving 
\begin{equation*}
    \begin{aligned}
    \minimize_{x\in \mathbb{B}^n}\quad & \tfrac{1}{2}\langle x,x\rangle. 
    \end{aligned}
\end{equation*}
We let $\sigma_{\xi}=\epsilon_M^4$, with $\delta=1$ or $\delta=10^{-100}$ (two extremes) and compare Theorem~\ref{thm:strong:convex:noise:quadratic} (CS algorithm) against a state-of-the-art multi-point method~\cite[Theorem~5.1]{ref:akhavan2020exploiting} ($\beta$ algorithm)\footnote{The small $\sigma_{\xi}$ aides the exposition as a larger $\sigma_{\xi}$ would merely delay the effect.}. Their stepsize equals ours, but their smoothing parameter equals
\begin{equation*}
    \delta_k = \left(\frac{3n^2\sigma_{\xi}}{4k+9n^2} \right)^{\tfrac{1}{4}},\quad k=1,2,\dots,K. 
\end{equation*}
In Figure~\ref{fig:strongcvx} we show for 250 experiments ($x_1\overset{i.i.d.}{\sim} \mathcal{N}(0,I_n)$), the differences in convergence. Indeed, the proposed complex-step does not suffer from cancellation errors as can be seen in Figure~\ref{fig:sim:cont}.  
The reason why $\delta=10^{-100}$ works (unreasonably) well is due to the constraints and the averaging, each iteration lives on $\partial \mathbb{B}^n$ (recall error terms of the form $1/\delta^2$). Although the setting is somewhat esoteric, this example does show the possibility of catastrophic cancellation and how to resolve it.  
}
\end{example}

Similar to Theorem~\ref{thm:strong:convex::unconstrained:noise}, we analyze unconstrained zeroth-order optimization when $f$ is quadratic. 

\begin{theorem}[Convergence rate of Algorithm~\ref{alg:convex_unconstrained}~(a) with noise, $f$ being quadratic]
\label{thm:strong:convex::unconstrained:noise:quadratic}
Let $f\in C^{\omega}(\mathcal{D})$ be a $\tau(f)$-strongly convex quadratic function satisfying Assumption~\ref{ass:holomorphic} for some $\bar{\delta}\in (0,1)$ with $x^{\star}\in \mathrm{int}(\mathcal{D})$. Suppose that $f$ has a Lipschitz gradient and constant Hessian, that is, ~\eqref{equ:grad:Lipschitz} and~\eqref{equ:Hess:Lip} hold, for $L_1(f)>0$ and $L_2(f)=0$, respectively.   
 Let $\{x_k\}_{k\ge 1}$ be the sequence of iterates generated by Algorithm~\ref{alg:convex_unconstrained}~(a) for 
 \begin{align*}
\everymath={\displaystyle}
 \begin{array}{llll}
 &\mu_k = \tfrac{1}{\tau(f) K},\quad  &\delta_k = \delta,\quad &k=1,\dots,K_0,  \\
   &\mu_k = \tfrac{2}{\tau(f) k},\quad  &\delta_k = \delta,\quad &k=K_0+1,\dots,K,
 \end{array}
 \end{align*}
 with $K_0 = \floor*{\tfrac{4 nL_1(f)^2}{\tau(f)^2}}$ and $\delta\in (0,\kappa\bar{\delta}]$ for some $\kappa\in (0,1)$. Then, if the oracle satisfies Assumption~\ref{ass:complex:oracle} and $K\geq 2K_0$ we incur for $\bar{x}_{K_0,K}=\tfrac{1}{K-K_0}\sum^K_{k=K_0+1}x_k$
 the optimization error 
\begin{equation}
\label{equ:quad:unconstrained:rate}
\begin{aligned}
    \mathbb{E}[f(\bar{x}_{K_0,K})-f(x^{\star})]
    &\leq O\left(\frac{nL_1(f)^2}{\tau(f) }\|x_1-x^{\star}\|_2^2K^{-1}\right) + \widetilde{O}\left( \frac{n^2\sigma_{\xi}}{\tau(f)\delta^2}K^{-1}\right). 
\end{aligned}
\end{equation}
\end{theorem}
\begin{proof}
The proof will be a combination of Theorem~\ref{thm:strong:convex::unconstrained:noise} and Theorem~\ref{thm:strong:convex:noise:quadratic}.

Again, set $a_k=\|x_k-x^{\star}\|_2^2$, then
\begin{align*}
\everymath={\displaystyle}
\begin{array}{lcl}
    \mathbb{E}[f(x_k)-f(x^{\star})|x_k]
    &\overset{\eqref{equ:grad:delta:approx:error}}{\leq} & \tfrac{1}{2\mu_k}\mathbb{E}[a_k-a_{k+1}|x_k] + \tfrac{\mu_k}{2}\mathbb{E}[\|g_{\delta_k}(x_k)\|_2^2|x_k]-\tfrac{\tau(f)}{2}\mathbb{E}[a_k|x_k].
    \end{array}
\end{align*}
Next, let $r_k=\mathbb{E}[a_k]$ such that by $\|\nabla f(x_k)\|_2^2\leq L_1(f)^2 \|x_k-x^{\star}\|_2^2$ we can write
\begin{equation}
\label{equ:fk:quad:unconstrained}
\begin{aligned}
    \mathbb{E}[f(x_k)-f(x^{\star})]
    &\leq \tfrac{1}{\mu_k}(r_k-r_{k+1})-{\tau(f)}r_k  + {\mu_k}\left(  nL_1(f)^2 r_k+ \tfrac{n^2}{\delta_k^2}\sigma_{\xi} \right)
\end{aligned}
\end{equation}

Now we use the step- and smoothingsize for $k=1,\dots,K_0$, that is, $\mu_k= 1/(\tau(f) K)$, $\delta_k = \delta $, and observe that 

\begin{align*}
    r_{k+1}&\leq r_k-\tau(f)\mu_kr_k  + \mu_k^2\left( nL_1(f)^2 r_k+ \tfrac{n^2}{\delta_k^2}\sigma_{\xi} \right)= \left(1-\tfrac{1}{K}+\tfrac{L_1(f)^2}{(\tau(f)K)^2}n\right)r_k + \nu_K = a_K r_k + \nu_K
\end{align*}
for $a_k$ as between brackets and $\nu_K$ defined as 
\begin{align*}
    \nu_K =&  \tfrac{n^2}{(\delta \tau(f) K)^2}\sigma_{\xi}. 
\end{align*}
We now proceed with bounding $r_{K_0+1}$. As in~\cite{ref:akhavan2020exploiting}, set 
\begin{equation*}
    q_K = 1+\tfrac{L_1(f)^2}{(\tau(f)K)^2}n
\end{equation*}
by iterating over $r_k$ it follows that
\begin{align*}
    r_{K_0+1}\leq a_K^{K_0}r_1 + \textstyle\sum^{K_0-1}_{i=0}a_K^i \nu_K \leq \left( r_1 + \tfrac{(\tau(f)K)^2}{L_1(f)^2n}\nu_K \right)q_K^{K_0}. 
\end{align*}
Now assume that that $K_0$ is as in the theorem, 
then as $\log(1+x)\leq x$
\begin{align*}
    q_K^{K_0} =& \mathrm{exp}\left(K_0 \log\left(1+\tfrac{L_1(f)^2}{(\tau(f)K)^2}n \right) \right)
    \leq  \mathrm{exp}\left(\tfrac{4nL_1(f)^2}{\tau(f)^2} \log\left(1+\tfrac{L_1(f)^2}{(\tau(f)K)^2}n \right) \right)
    \leq  \mathrm{exp}\left(\tfrac{4 n^2 L_1(f)^4}{\tau(f)^4 K^2} \right).
\end{align*}
Fix any $\bar{C}\in (0,\tfrac{1}{16})$, when 
\begin{equation*}
    K\geq \sqrt{\tfrac{4 n^2 L_1(f)^4}{\tau(f)^4\bar{C}}}
\end{equation*}
then $K\geq 2K_0$ and $q^{K_0}_K\leq e^{\bar{C}}=C_4$. 
As such,
\begin{align*}
    r_{K_0+1}\leq& \left( r_1 + \tfrac{(\tau(f)K)^2}{L_1(f)^2n}\nu_K \right)C_4
    \leq \left( r_1 + \tfrac{(\tau(f)K)^2}{L_1(f)^2n}\tfrac{n^2}{(\delta \tau(f) K)^2}\sigma_{\xi} \right)C_4 
    = \left( r_1 + \tfrac{n}{\delta^2 L_1(f)^2}\sigma_{\xi} \right)C_4.
\end{align*}
Now we return to our normal step- and smoothingsizes, that is $\mu_k=2/(\tau(f) k)$, $\delta_k = \delta $, for $k\geq K_0+1$. By plugging this into~\eqref{equ:fk:quad:unconstrained} we get
\begin{equation*}
\begin{aligned}
    (K-K_0)\mathbb{E}[f(\bar{x}_{K_0,K})-f(x^{\star})]
    \leq &\textstyle \sum^K_{k=K_0+1}\tfrac{\tau(f)k}{2}(r_k-r_{k+1})-{\tau(f)}r_k + \tfrac{2}{\tau(f)k}nL_1(f)^2 r_k \\ 
     +&\textstyle\sum^K_{k=K_0+1} \tfrac{2}{\tau(f)k}\tfrac{n^2}{\delta^2}\sigma_{\xi} 
\end{aligned}
\end{equation*}
By construction of $K_0$ we have that for $k\geq K_0+1$, $\tau(f)/2\geq (2nL_1(f)^2)/(\tau(f) k)$. Hence 
\begin{equation*}
\begin{aligned}
    (K-K_0)\mathbb{E}[f(\bar{x}_{K_0,K})-f(x^{\star})]
    &\leq \tfrac{\tau(f)}{2}\textstyle\sum^K_{k=K_0+1}k(r_k-r_{k+1})- r_k + U_{K_0,K}.
\end{aligned}
\end{equation*}
where by~\eqref{equ:log:bound} 
\begin{align*}
    U_{K_0,K} = \textstyle\sum^K_{k=K_0+1}   \tfrac{2}{\tau(f)k}\tfrac{n^2}{\delta^2}\sigma_{\xi} 
    \leq \tfrac{2}{\tau(f)}\tfrac{n^2}{\delta^2}\sigma_{\xi}(\log(K)+1). 
\end{align*}
As demonstrated in~\cite{ref:akhavan2020exploiting}, one can now construct the bound $\textstyle\sum^K_{k=K_0+1} k(r_k - r_{k+1})-r_k \leq K_0 r_{K_0+1}$ where the last term is exactly the term we could bound before. In combination with the bound on $K_0$ itself, we find that 
\begin{equation*}
\begin{aligned}
    (K-K_0)\mathbb{E}[f(\bar{x}_{K_0,K})-f(x^{\star})]
    &\leq \tfrac{\tau(f)}{2}\tfrac{4nL_1(f)^2}{\tau(f)^2}\left( r_1 + \tfrac{n}{\delta^2 L_1(f)^2}\sigma_{\xi} \right)C_4+\tfrac{2}{\tau(f)}\tfrac{n^2}{\delta^2}\sigma_{\xi}(\log(K)+1).
\end{aligned}
\end{equation*}
By our selection of $\bar{C}$ we have that $K\geq 2 K_0$ and as such 
\begin{equation*}
\begin{aligned}
 \mathbb{E}[f(\bar{x}_{K_0,K})-f(x^{\star})]
    &\leq \tfrac{2nL_1(f)^2}{\tau(f) K}\left( r_1 + \tfrac{n}{\delta^2 L_1(f)^2}\sigma_{\xi} \right)C_4+\tfrac{2n^2}{\delta^2\tau(f)K}\sigma_{\xi}(\log(K)+1).
\end{aligned}
\end{equation*}
Now, reordering terms yields~\eqref{equ:quad:unconstrained:rate}. 
\end{proof}

It is important to highlight that the stepsizes for the quadratic cases are identical to the general case. As such, no knowledge of the quadratic nature is required, but if $f$ happens to be quadratic, the algorithm performs optimally. 

\subsection{Online optimization}
Online optimization shows up in settings where the objective might change due to the presence of more information, say, when more data becomes available. In the online case one is interested in bounding the regret of the form
\begin{equation}
    \tfrac{1}{K}\textstyle \sum^K_{k=1}\mathbb{E}[f_k(x_k)] - \inf_{x\in \mathcal{D}} \tfrac{1}{K}\sum^K_{k=1} f_k(x). 
\end{equation}
As in for example~\cite{ref:bach:smooth}, the proof techniques are largely the same as for the stochastic cases above. We consider the following setting to exemplify the possibilities. Note, here the algorithm proceeds as 
\begin{equation*}
   x_{k+1} = \Pi_{\mathcal{K}}\left( x_k - \mu_k\cdot \left(\tfrac{n}{\delta_k}\Im\left(f_{k+1}(x_k + i \delta_k u_k)\right)u_k+\tfrac{n}{\delta_k}\xi_k u_k \right)\right). 
\end{equation*}

\begin{theorem}[Online optimization, convergence rate of Algorithm~\ref{alg:convex_unconstrained}~(b) with noise, $f_k$ being quadratic]
\label{thm:online:strong:convex:noise:quadratic}
Let all $f_k\in C^{\omega}(\mathcal{D})$ be $\tau(f_k)$-strongly convex functions satisfying Assumption~\ref{ass:holomorphic} for some $\bar{\delta}\in (0,1)$ and let $\mathcal{K}\subset \mathcal{D}$ be a compact convex set. Suppose that all $f_k$ have a mutual Lipschitz gradient and a constant Hessian over $\mathcal{K}$, that is, ~\eqref{equ:grad:Lipschitz} and~\eqref{equ:Hess:Lip} hold, for some constant $L_1(f_k)>0$ and $L_2(f_k)=0$, respectively.   
 Set $\tau=\min_{k}\tau(f_k)$ and let $\{x_k\}_{k\ge 1}$ be the sequence of iterates generated by Algorithm~\ref{alg:convex_unconstrained}~(b) with stepsize $ \mu_k  = {2}/({\tau k})$ and the sequence of smoothing parameters defined for all $k\geq 1$ by $\delta_k = \delta $ with $\delta\in (0,\kappa\bar{\delta}]$ for some $\kappa\in (0,1)$. Then, if the oracle satisfies Assumption~\ref{ass:complex:oracle} we incur the regret 
\begin{equation}
\label{equ:online:quad:regret}
    \tfrac{1}{K}\textstyle\sum^{K}_{k=1}\mathbb{E}[f_k(x_k)]-\inf_{x\in \mathcal{K}}\tfrac{1}{K}\sum^K_{k=1}f_k(x)
    \leq  \widetilde{O}\left(\frac{n^2}{\tau }K^{-1}\right)
\end{equation}
\end{theorem}
\begin{proof}
As~\eqref{equ:f:ineq:tau:th1} holds for any $x$ and not just $x^{\star}$, the proof is effectively identical to that of Theorem~\ref{thm:strong:convex:noise:quadratic}. By appealing to that proof we have the following bound immediately
\begin{align*}
\everymath={\displaystyle}
\begin{array}{lcl}
    \textstyle\sum^{K}_{k=1}\mathbb{E}[f_k(x_k)-f_k(x)]
    &\leq&  \tfrac{n^2}{\tau }\left[ \tfrac{1}{n}G^2+ \tfrac{1}{\delta^2}\sigma_{\xi}\right](\log(K)+1)
    \end{array}
\end{align*}
and as such we obtain the regret bound~\eqref{equ:online:quad:regret}.
\end{proof}


\subsection{Numerical estimation of $\tau(f)$}
\label{sec:est:tau}
As most regret bounds and stepsizes contain terms of the form $1/\tau(f)$ one should take care in estimating the strong convexity parameter $\tau(f)$. An arbitrarily small $\tau(f)$ complies with the definition but could lead for instance to numerical overflow due to large stepsizes. In fact, it is known that either under- or overestimating $\tau(f)$ can have detrimental effects on convergence properties, especially in accelerated schemes~\cite{ref:oc2015adaptive}. When one has access to gradients, line-search-like schemes are possible to estimate both $\tau(f)$ and $L_1(f)$~\cite{ref:nesterov2013gradient}. However, when the gradient direction is random, this is less straight-forward. 

Fitting a quadratic model using a (recursive) least squares approach can grossly overestimate $\tau(f)$. For example, consider the function $f(x)=x^4+x^2+\tfrac{1}{2}\lambda x^2$. One might have access to $\lambda>0$,~\textit{e.g.}, by means of being a regularization parameter. Then, fitting a quadratic model to this function yields (asymptotically) a strong convexity estimate of $6+\lambda$ instead of $2+\lambda$.  

We propose simple routines to estimate the largest $\tau(f)$ satisfying~\eqref{equ:strong:convex}, denoted $\bar{\tau}(f)$. Here we exploit the fact that we have a sequence of function evaluations, which remain commonly and unfortunately unused in this line of zeroth-order optimization schemes. We also assume to have knowledge of some lower bound $\tau_0>0$ such that $\tau_0\leq \bar{\tau}(f)$, which is frequently available due to regularization. In terms of the dimension $n$ we identify two regimes, small (medium) scale $n\leq 10^3$ and large scale $n>10^3$.  

\begin{enumerate}[(i)]
    \item (Small scale): Using the data at hand we can construct an explicit quadratic model in $(P,q,r)$ that bounds $f(x)$ from below. Due to the inherent randomness, $f\in C^{\omega}$, and the possibility of selecting $\delta_k$ close to $0$, one has (for $\delta_k$ sufficiently small) a sufficiently accurate quadratic model by using $N(n)\geq \tfrac{1}{2}(n+1)(n+2)$ data-points in the following \textbf{\textit{semidefinite program}} (SDP)   
\begin{equation}
\label{equ:SDP:tau}
    \begin{aligned}
    \minimize_{(P,q,r)} \quad & \textstyle\sum^{K'+N(n)-1}_{k=K'}\Re(f(x_k+i\delta_k u_k))-\langle \tfrac{1}{2}Px_k+q,x_k \rangle -r + \delta_k^2 \langle \tfrac{1}{2}Pu_k, u_k \rangle\\
    \subjectto \quad & P\in \mathcal{S}^n_{\succ 0},\,q\in \mathbb{R}^n,\,r\in \mathbb{R},\,P\succeq \tau_0 I_n,\\
    \quad & \Re(f(x_k+i\delta_k u_k))\geq \langle \tfrac{1}{2}Px_k+q,x_k \rangle + r - \delta_k^2 \langle \tfrac{1}{2}Pu_k, u_k \rangle,\\ 
    \quad & \text{for }k=K',\dots,K'+N(n)-1, 
    \end{aligned}
\end{equation}
for some $K'\geq 1$. Now, an approximation of $\bar{\tau}(f)$ follows by setting $\widehat{\tau}(f)=\min_i\{\lambda_i(P^{\star})\}$. Indeed, models as such can now also be used to further fine-tune the proposed algorithms. 

Let us elaborate on the aforementioned claims in the \textit{unconstrained} case. The constrained case is less predictable. First of all, we want to have a tight quadratic model that approximates the function $f$ from below. To use our data economically, we have to do with samples of the form $\Re(f(x_k+i\delta_k u_k))$ instead of $f(x_k)$. Note, these samples might be corrupted by noise. Then, the objective in combination with the inequality constraints in~\eqref{equ:SDP:tau} enforce that $(P,q,r)$ parametrizes a quadratic model, approximately from below, that is as close as possible to the available data. To parametrize this model one needs at most $n(n+1)/2+n+1=\tfrac{1}{2}(n+1)(n+2)$ data-points indeed. Here, a data-point compromises the $4$-tuple $(x_k,u_k,\delta_k,f(x_k+i\delta_k u_k))$. Now, as we sample $u_1,u_2,\dots$ uniformly and independently from $\mathbb{S}^{n-1}$, the set $\{u_1,u_2,\dots,u_n\}$ will $\mathbb{P}$-\textit{a.s.} span $\mathbb{R}^n$. Then, as the noise terms $\xi_1,\xi_2,\dots$ that potentially enter the oracle are independent of $u_1,u_2,\dots$ and $f\neq 0$ can only vanish on sets of measure $0$ by the real-analytic assumption, we must have that $x_k$ is $\mathbb{P}$-\textit{a.s.} not parallel to $x_{k+1}$,~\textit{cf.}~Algorithm~\ref{alg:convex_unconstrained}.   

    \item (Large scale): When $n$ is large, we follow the ideas as set forth in~\cite{ref:ahmadi2019dsos}. Denote by $\mathsf{dd}^n$ the \textit{\textbf{diagonally dominant}} matrices in $\mathbb{R}^{n\times n}$. That is, $A\in \mathsf{dd}^n$ when $a_{ii}\geq \sum_{j\neq i}|a_{ij}|$ for all $i\in [n]$. This allows for a polytopic representation of the constraint $P-\tau_0 I_n\succeq 0$.
    Now we transform~\eqref{equ:SDP:tau} in the diagonally dominant program (DDP) by identifying $p\in \mathbb{R}^{n(n+1)/2}$ with $\mathrm{svec}(P)$, that is, $P$ is not an additional decision variable but merely an auxiliary variable to simplify notation  
\begin{equation}
\label{equ:DDP:tau}
    \begin{aligned}
    \minimize_{(p,q,r)} \quad & \textstyle\sum^{K'+N(n)-1}_{k=K'}\Re(f(x_k+i\delta_k u_k))-\tfrac{1}{2}\langle p,x_k \otimes_s x_k-\delta_k^2 u_k \otimes_s u_k \rangle-\langle q,x_k \rangle -r \\
    \subjectto \quad & p\in \mathbb{R}^{n(n+1)/2},\,q\in \mathbb{R}^n,\,r\in \mathbb{R},\,P=\mathrm{smat}(p),\\
    \quad & P_{ii}-\tau_0\geq \textstyle\sum_{j\neq i}|P_{ij}|,\,i=1,\dots n,\\
    \quad & \Re(f(x_k+i\delta_k u_k))\geq \tfrac{1}{2}\langle p,x_k \otimes_s x_k-\delta_k^2 u_k \otimes_s u_k \rangle + \langle q,x_k \rangle   + r,\\ 
    \quad & \text{for }k=K',\dots,K'+N(n)-1. 
    \end{aligned}
\end{equation}
Here,~$\otimes_s$ denotes the symmetric Kronecker product. See~\cite{ref:majumdar2020SDP} for a recent survey on large-scale SDPs and Section~\ref{sec:DDP} for more on the DDP-based approximation of $\bar{\tau}(f)$. Specifically, we can iteratively improve the basis in~\eqref{equ:DDP:tau}, such that $p^{\star}$ with respect to~\eqref{equ:DDP:tau} converges weakly to $P^{\star}$ with respect to~\eqref{equ:SDP:tau}.
\end{enumerate}

\begin{example}[Numerical performance of $\tau(f)$ estimation]
\label{ex:reg}
\upshape{
To show how the proposed estimation scheme for $\tau(f)$ can be beneficial we look at a transparent (closed-form solutions are available) example. 
Consider the $\ell_2$-regularized least-squares problem
\begin{equation*}
    \minimize_{x\in r \mathbb{B}^n}\quad \tfrac{1}{2}\langle Ax-b, Ax-b \rangle  + \tfrac{1}{2}\lambda \langle x,x \rangle
\end{equation*}
for $r>0$ such that $x^{\star}\in \mathrm{int}( r \mathbb{B}^n)$. In many problems one might have knowledge of the regularization parameter $\lambda>0$ but not of the remaining objective terms. As such, we start with $\tau_0= \lambda$ and use the SDP formulation~\eqref{equ:SDP:tau} to approximate $\bar{\tau}(f)$ from below by $\widehat{\tau}(f)$. We do 250 experiments ($x_1\overset{i.i.d.}{\sim} \mathcal{N}(0,I_n)$, $\mathrm{vec}(A)\overset{i.i.d.}{\sim}  \mathcal{N}(0,I_{mn})$, $b\overset{i.i.d.}{\sim}  \mathcal{N}(0,I_m)$) for $n=10$, $m=20$ and $\lambda=10^{-4}$. We plug the estimation scheme into Algorithm~\ref{alg:convex_unconstrained}~(b) (for $\delta_k=\epsilon_M$ and $\sigma_{\xi}=\epsilon_M^4$), that is, compute $\widehat{\tau}(f)$ once at $K=\tfrac{1}{2}(n+1)(n+2)$, and show the results in Figure~\ref{fig:cvxreg}. The approximation clearly speeds up the convergence and closely resembles that under $\bar{\tau}(f)$. Section~\ref{sec:DDP} (Appendix) presents a similar example for~\eqref{equ:DDP:tau}.     
}
\end{example}

The take away of this section is not only a routine to estimate $\tau(f)$, but also the observation that this can be done directly using the complex function evaluations of the form $f(x_k+i\delta_k u_k)$.


\section{Outlook: nonconvex zeroth-order optimization}
\label{sec:non:convex}

At last we consider a critical point in a possibly non-convex program. Note, we do not assume that our function is \textit{locally} convex. We exploit that the gradient of $f\in C^{\omega}(\mathcal{D})$ is uniformly bounded over any compact subset of $\mathcal{D}$. 

\begin{theorem}[Convergence rate of Algorithm~\ref{alg:convex_unconstrained}~(b) to a critical point]
\label{thm:ncvx}
Let $f\in C^{\omega}(\mathcal{D})$ be a | not necessarily convex | function that satisfies Assumption~\ref{ass:holomorphic} for some $\bar{\delta}\in (0,1)$. Suppose that $f$ has a Lipschitz gradient and Hessian on $\mathcal{K}\subset\mathcal{D}$, that is, ~\eqref{equ:grad:Lipschitz} and~\eqref{equ:Hess:Lip} hold, for some constants $L_1(f)>0$ and $L_2(f)\geq 0$, respectively.   
Let $\{x_k\}_{k\ge 1}$ be the sequence of iterates generated by Algorithm~\ref{alg:convex_unconstrained}~(b) with stepsize $\mu_k = {1}/({n L_1(f)k^{2/3}})$ and the sequence of smoothing parameters defined for all $k\geq 1$ by $\delta_k = \delta k^{-1/6}$ with $\delta\in (0,\kappa\bar{\delta}]$ for some $\kappa\in (0,1)$. Let $x^{\star}$ be a global minimum of $f$, then,  
\begin{equation}
\label{equ:non:convex:rate}
   \begin{aligned}
    \min_{k\geq 1} \mathbb{E}[\|\nabla f(x_k)\|^2_2] \leq O\left({nL_1(f)(f(x_1)-f(x^{\star}))}K^{-\tfrac{1}{3}} \right)+\widetilde{O}\left(n{(\delta^2  + \frac{\sigma_{\xi}}{\delta^2})}K^{-\tfrac{1}{3}} \right).
\end{aligned} 
\end{equation}
\end{theorem}
\begin{proof}
Our proof will be similar to constructions as set forth in~\cite{nesterov2003introductory}.
As $f\in C^{1,1}_{L_1(f)}$ one has 
\begin{align*}
\everymath={\displaystyle}
\begin{array}{lcl}
    f(x_{k+1}) & \overset{\eqref{equ:L1:upper}}{\leq}& f(x_k) - \mu_k \langle \nabla f(x_k), g_{\delta_k}(x_k) \rangle + \tfrac{1}{2}\mu_k^2 L_1(f) \|g_{\delta_k}(x_k)\|_2^2\\
    &=& f(x_k) - \mu_k\|\nabla f(x_k)\|_2^2 - \mu_k \langle \nabla f(x_k),g_{\delta_k}(x_k)-\nabla f(x_k)\rangle + \tfrac{1}{2}\mu_k^2 L_1(f) \|g_{\delta_k}(x_k)\|_2^2.
\end{array}
\end{align*}
Now taking expectation, applying the Cauchy-Schwarz inequality and using both~\eqref{equ:grad:delta:approx:error} and~\eqref{equ:g:var} results in 
\begin{align*}
    \mathbb{E}_{u_k\sim \mathbb{S}^{n-1}}\left[f(x_{k+1})|x_k \right] \leq& f(x_k) - \mu_k \|\nabla f(x_k)\|_2^2 + \mu_k {C_1 n\delta_k^2}\|\nabla f(x_k)\|_2 \\
    & + \tfrac{1}{2}\mu_k^2 L_1(f)\left( n\|\nabla f(x_k)\|_2^2 + {C_2 n^2 \delta_k^4 }+C_3n^2\delta_k^2\|\nabla f(x_k)\|_2+\tfrac{n^2}{\delta_k^2}\sigma_{\xi}\right). 
\end{align*}
Then, taking expectation over $u_1,\dots,u_{k-1}$, plugging in our stepsize $\mu_k = 1/(n L_1(f)k^{\tfrac{2}{3}})$ applying Jensen's inequality and rearranging yields 
\begin{equation}
\label{equ:ncvx:grad:bound}
    \begin{aligned}
    k^{-\tfrac{2}{3}}\tfrac{1}{2nL_1(f)}\mathbb{E}\left[\|\nabla f(x_k)\|_2^2\right]
    \leq& \mathbb{E}[f(x_k)-f(x_{k+1})]\\
    &+k^{-\tfrac{2}{3}}\tfrac{C_4 \delta_k^2 }{L_1(f) }\mathbb{E}\left[\|\nabla f(x_k)\|_2\right] +k^{-\tfrac{4}{3}}\tfrac{C_2  \delta_k^4 }{2 L_1(f)}+k^{-\tfrac{4}{3}}\tfrac{1}{2L_1(f)\delta_k^2 }\sigma_{\xi}. 
\end{aligned}
\end{equation}
As we consider a global minimum we have that $f(x_k)\geq f(x^{\star})$. Now, define $\phi_k=(1/k^{\tfrac{2}{3}})\mathbb{E}[\|\nabla f(x_k)\|^2_2]$ and $\sup_{x\in \mathcal{K}}\|\nabla f(x)\|_2=G<+\infty$, then, a telescoping argument yields 
\begin{align*}
    \textstyle\sum^K_{k=1}\phi_K \leq 2nL_1(f)(f(x_1)-f(x^{\star})) + 2n C_4G\sum^K_{k=1}k^{-\tfrac{2}{3}}\delta_k^2  + n C_2\sum^K_{k=1}k^{-\tfrac{4}{3}}\delta_k^4+\sum^K_{k=1}k^{-\tfrac{4}{3}}\tfrac{n}{\delta_k^2}\sigma_{\xi}. 
\end{align*}
Now plug in $\delta_k=\delta k^{-\tfrac{1}{6}}$ and get
\begin{align*}
    \textstyle\sum^K_{k=1}\phi_K &\leq \textstyle 2nL_1(f)(f(x_1)-f(x^{\star})) + {2n C_4 G}\sum^K_{k=1}k^{-1}\delta^2 + {n C_2 }\sum^K_{k=1}k^{-2}\delta^4+\sum^K_{k=1}k^{-1}\tfrac{n}{\delta^2}\sigma_{\xi}\\
    &\leq 2nL_1(f)(f(x_1)-f(x^{\star})) + C_5(\log(K)+1)(n G \delta^2 + \tfrac{n}{\delta^2}\sigma_{\xi})+C_6 n \delta^4.
\end{align*}
As such $\sum^K_{k=1}\phi_K\leq h(K)$, for $h(K)$ corresponding to the right-most term above. This implies that $\min_k \phi_k \leq {h(K)/K}$. By definition of $\phi_k$ we have $K^{-\tfrac{2}{3}}\min_k \mathbb{E}[\|\nabla f(x_k)\|^2_2]\leq \min_k \phi_k$. Combing these observations yields
\begin{equation*}
    \min_{k\geq 1} \mathbb{E}[\|\nabla f(x_k)\|^2_2] \leq K^{-\tfrac{1}{3}}\left(2nL_1(f)(f(x_1)-f(x^{\star})) + C_5 n(\log(K)+1)(G \delta^2 + \tfrac{\sigma_{\xi}}{\delta^2})+C_6 n \delta^4\right). 
\end{equation*}
\end{proof}

Although the rate is relatively slow~\textit{cf.}~\cite{ref:ghadimi2013stochastic}, the approach appears to be scalable, in contrast to common Monte Carlo methods~\cite{ref:polyak2017MC}. Sharpening and further generalizing Theorem~\ref{thm:ncvx} is left for future work.

Now we provide an numerical experiment, showing that Theorem~\ref{thm:ncvx} can handle the noise, in contrast to the nonconvex algorithm as proposed in~\cite{ref:JongeneelYueKuhnZO2021} for the deterministic setting. 
    \begin{example}[Himmelblau function]
\label{ex:nonconvex}
\upshape{
Consider optimizing a Himmelblau function over a closed ball centred at $0$, in particular, consider 
\begin{equation}
    \label{equ:himmelblau:constrained}
    \minimize_{x\in 6\mathbb{B}^2}\quad  \left(\big(x^{(1)}\big)^2+x^{(2)}-11\right)^2 + \left(x^{(1)}+\big(x^{(2)}\big)^2-7\right)^2. 
\end{equation}
The minimum value of~\eqref{equ:himmelblau:constrained} is $f^{\star}=0$. We will compare~\cite[Algorithm~1]{ref:JongeneelYueKuhnZO2021} with stepsize $\mu_k= \mu = 1/(nL_1(f))$ against Theorem~\ref{thm:ncvx}. That is, we compare a plain nonconvex algorithm (with $\delta_k=10^{-6}/k)$ against its counterpart that \textit{is} designed to handle noise (with $\delta_k = 10^{-6}k^{-1/6}$). We consider $8$ initial conditions (circles) and show the results in Figure~\ref{fig:noncvx}. The dark stars indicate minima of $f$, whereas the light star is merely a local minima. We see that the algorithm adapted to the noise can handle the perturbations well whereas the other algorithm diverges. Note that formalizing these observations is left for future work. 
\begin{figure*}[t!]
    \centering
    \begin{subfigure}[b]{0.3\textwidth}
        \includegraphics[width=\textwidth]{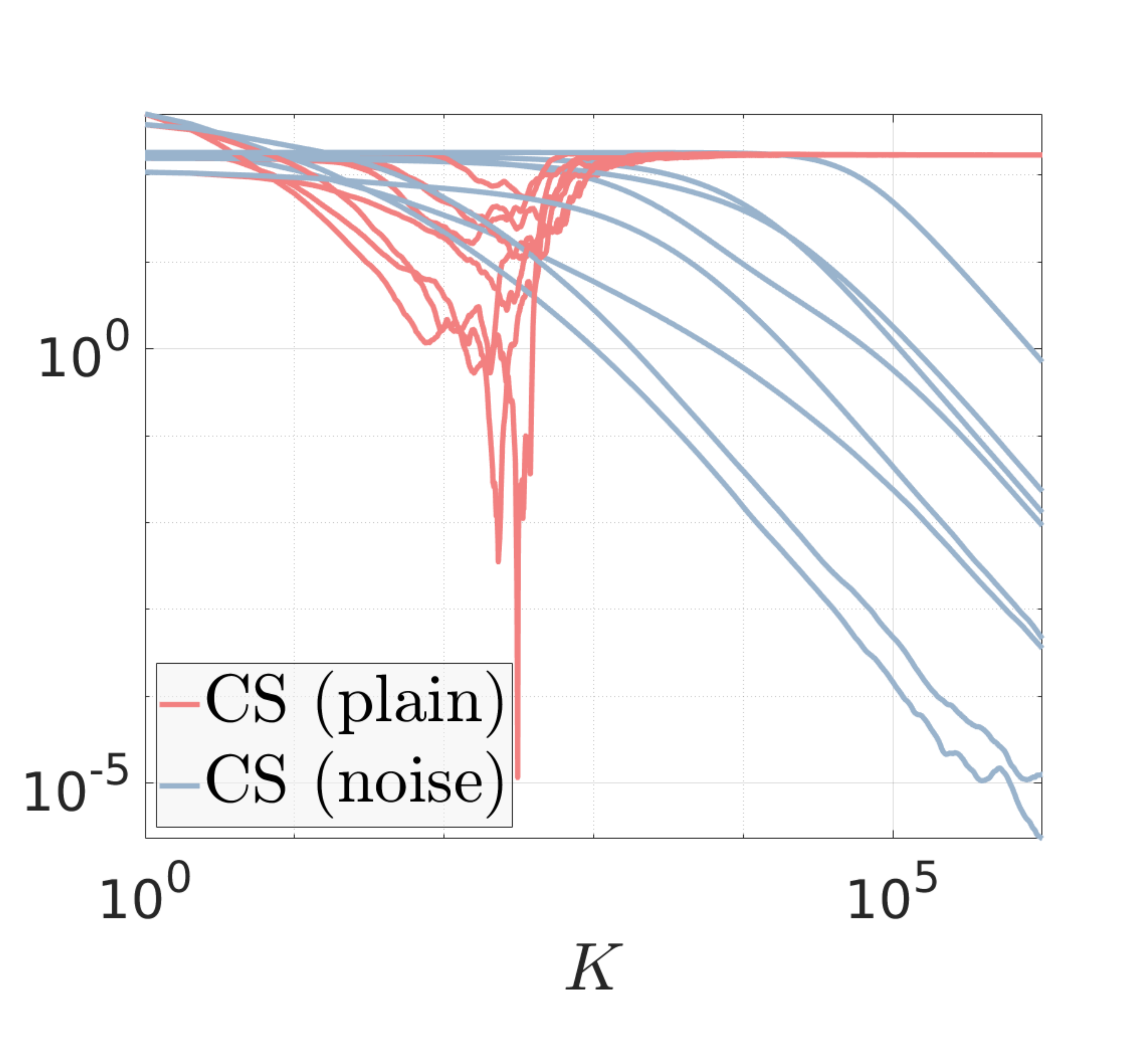}
        \caption{The error $f(\bar{x}_K)-f^{\star}$}
        \label{fig:nonconvexK}
    \end{subfigure}\quad
    \begin{subfigure}[b]{0.3\textwidth}
        \includegraphics[width=\textwidth]{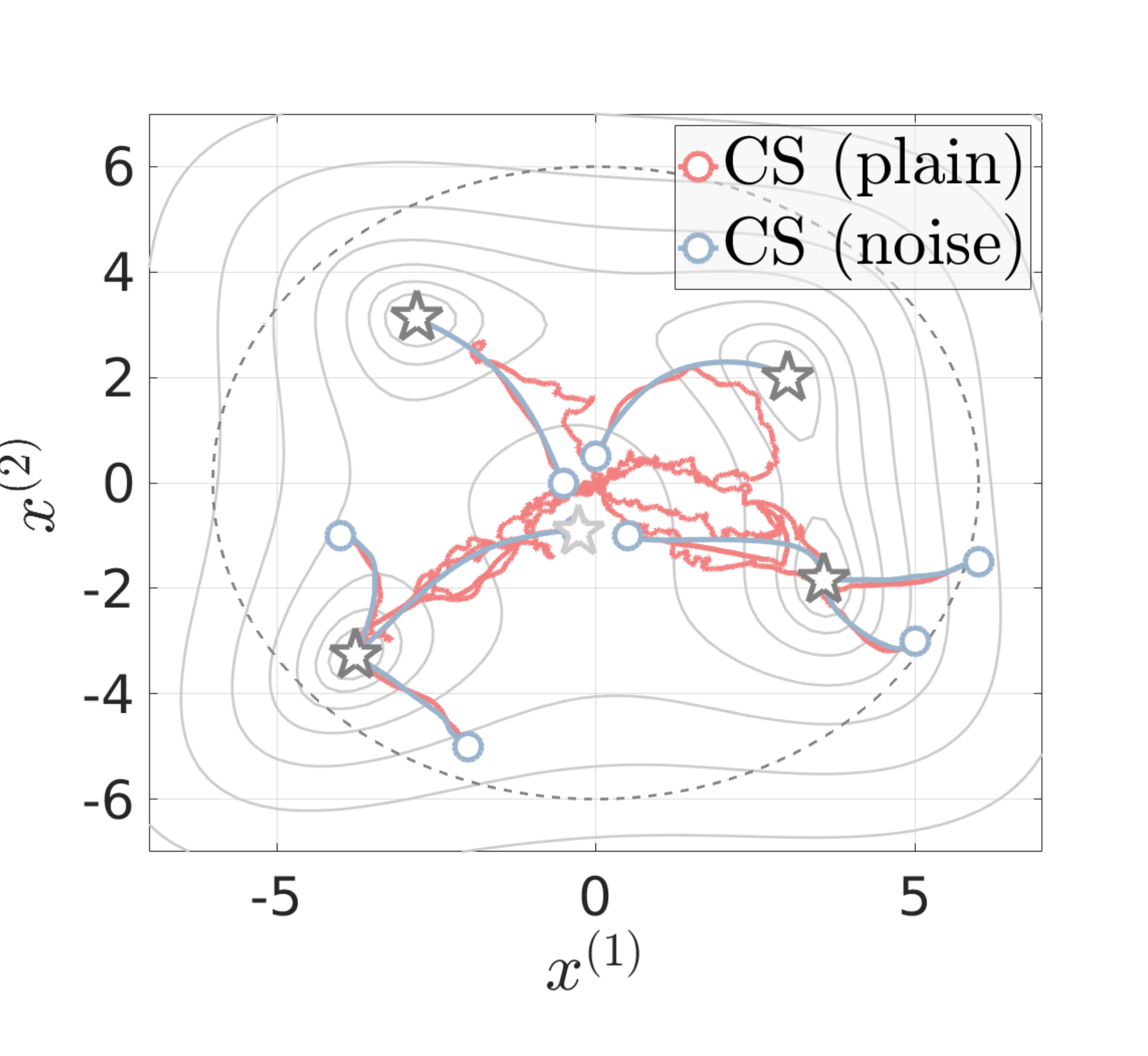}
        \caption{The trajectories $\bar{x}_1,\dots,\bar{x}_K$.}
        \label{fig:nonconvextrajec}
    \end{subfigure}
    \caption[]{Numerical outcomes of Example~\ref{ex:nonconvex}.}
    \label{fig:noncvx}
\end{figure*}
}
\end{example}

\begin{figure*}[t!]
    \centering
    \begin{subfigure}[b]{0.3\textwidth}
        \includegraphics[width=\textwidth]{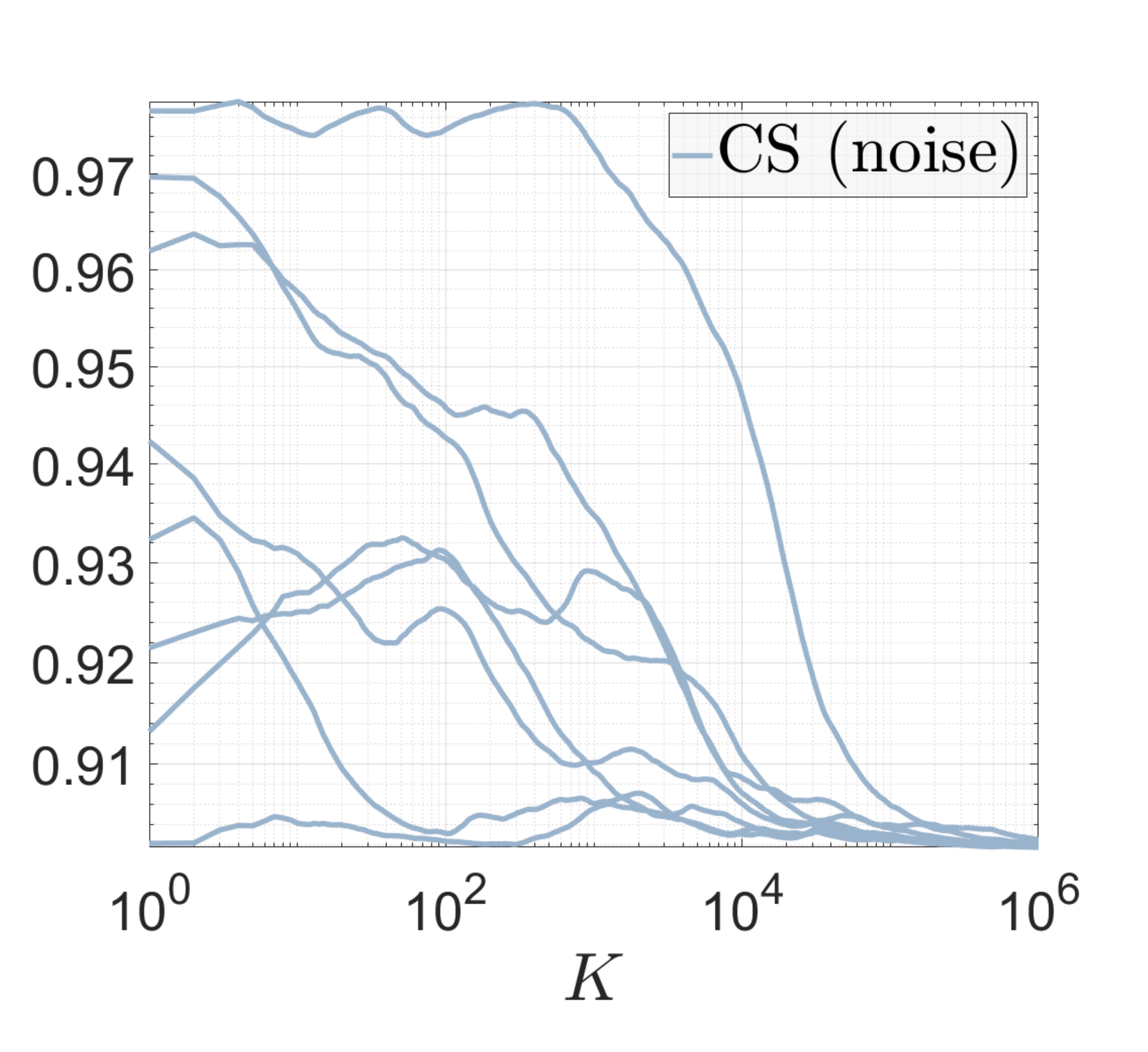}
        \caption{Example~\ref{ex:PDE}, convergence of the cost (averaged estimator).}
        \label{fig:PDEcost}
    \end{subfigure}\quad
    \begin{subfigure}[b]{0.3\textwidth}
        \includegraphics[width=\textwidth]{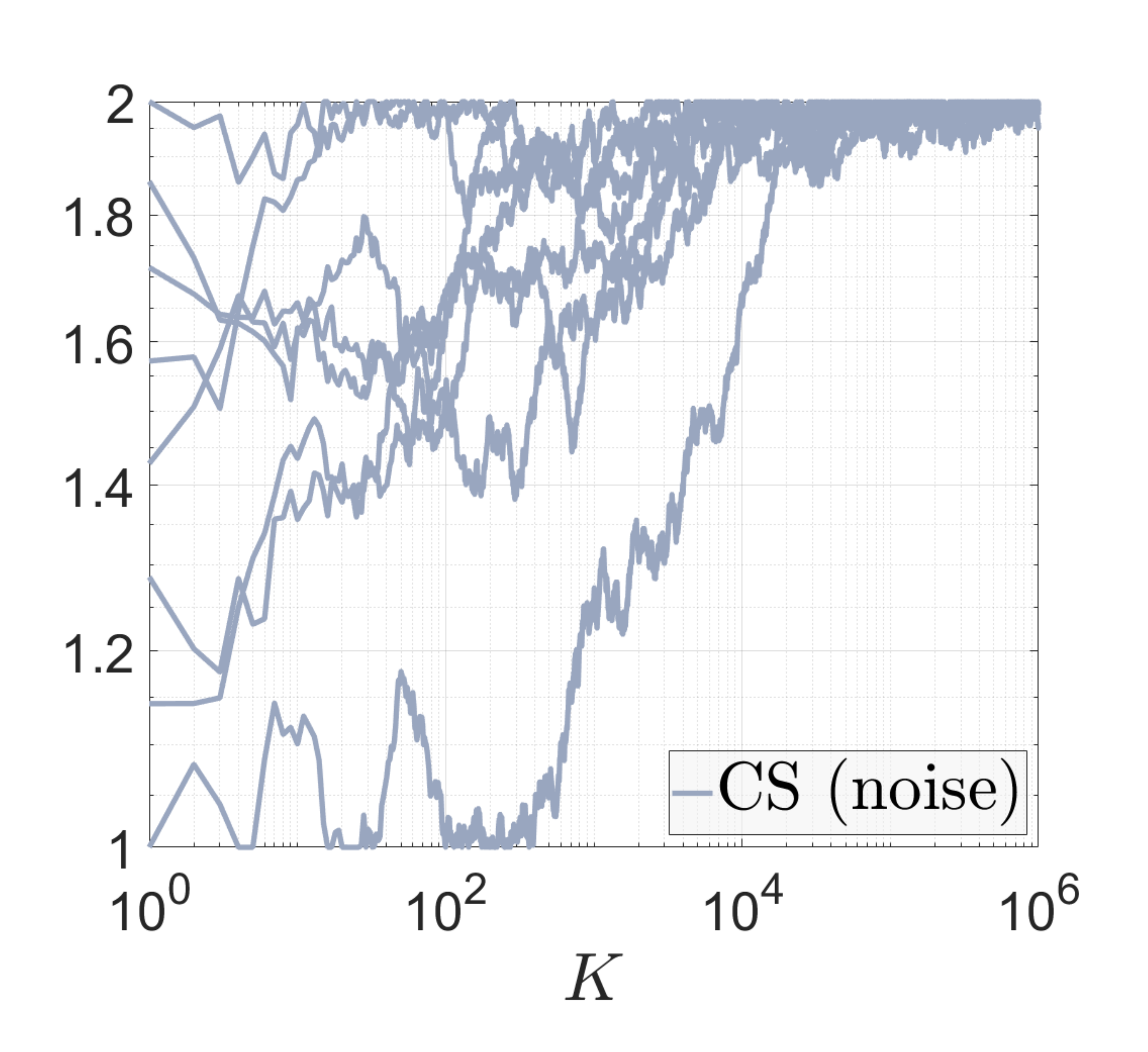}
        \caption{Example~\ref{ex:PDE}, convergence of $r_0\in [1,2]$ (non-averaged).}
        \label{fig:PDExK}
    \end{subfigure}\quad 
    \begin{subfigure}[b]{0.3\textwidth}
        \includegraphics[width=\textwidth]{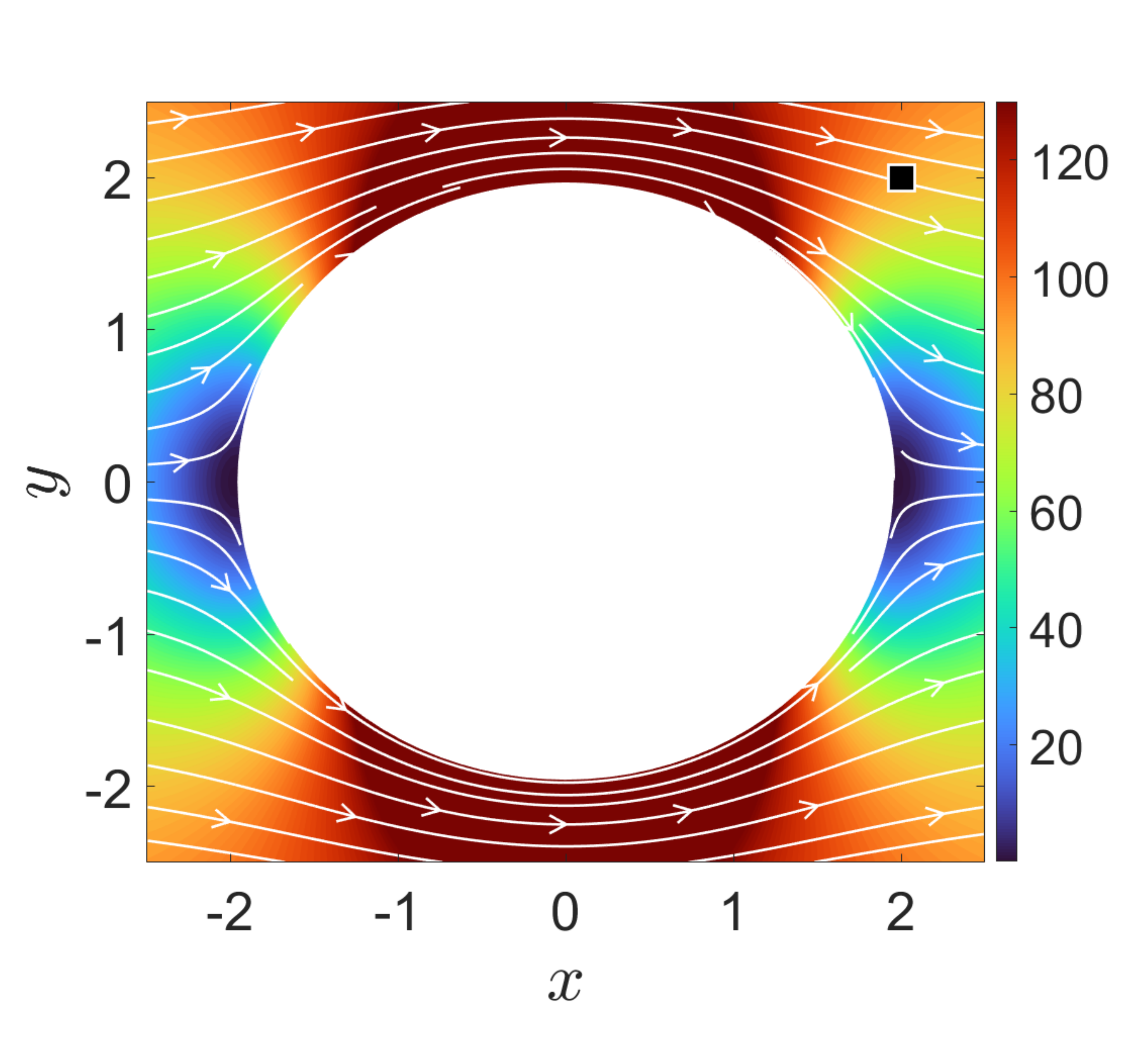}
        \caption{Example~\ref{ex:PDE}, typical result for $K=10^6$ (non-averaged).}
        \label{fig:PDE1}
    \end{subfigure}
    \caption[]{PDE-constrained optimization. Each figure displays all available data.}
    \label{fig:PDE}
\end{figure*}

We end with an example pertaining to partial differential equations (PDEs). PDEs are relevant as on the one hand, closed-form solutions are rare and numerical solutions (approximations) are often a necessity, on the other hand, analyticity of solutions has been studied since the early 50s, see for example~\cite{ref:morrey1958analyticity,ref:morrey1958analyticity2}.
\begin{example}[PDE-constrained optimization]
\label{ex:PDE}
\upshape{
PDEs can rarely be solved in closed-form and one commonly resorts to numerical schemes, however, schemes that often lend themselves to the complex-lifting as set forth in this article. In this example we show that there are already examples that meet the conditions of Theorem~\ref{thm:ncvx}. In particular, let $u$ be a velocity field on $\mathbb{R}^2$, with abuse of notation, $(x,y)$ denote the usual coordinates on $\mathbb{R}^2$. This velocity field is induced by a solid sphere in $\mathbb{R}^2$ that moves in the negative $x$-direction with a velocity $V$. We are interested in finding the optimal radius $r$ of this sphere such that norm of the velocity field at the point $(2,2)\in \mathbb{R}^2$ is minimized. When constraining the radius to the interval $[1,2]$, then under idealized conditions (incompressibility and irrotationality), we can consider the following PDE-constrained optimization problem 
\begin{equation}
    \label{equ:PDE:opt}
    \begin{aligned}
    \minimize_{r,u} \quad &\|u(2,2)\|_2^2\\
    \subjectto \quad &\mathrm{div}(u)=0,\,\mathrm{curl}(u)=0,\\
    &\langle u,n_x \rangle = \langle V , n_x \rangle, \, \forall x\in r\mathbb{S}^1,\,n_x\in (T_x r\mathbb{S}^1)^{\perp},\\
    &r\in [1,2],\,u\in C^1([-3,3]). 
    \end{aligned}
\end{equation}
As the PDE in~\eqref{equ:PDE:opt} admits a closed-form solution parametric in $r$\footnote{See for example Section~4.5.1 of the lectures notes by Dr. Evy Kersal\'e~\url{http://www1.maths.leeds.ac.uk/~kersale/2620/Notes/chapter_4.pdf}.}, one can easily bound $L_1(f)$,~\textit{e.g.}, we use $L_1(f)=10$. Moreover, we set $\sigma_{\xi}=10^{-12}$ to simulate numerical noise, set $\delta=10^{-6}$ and perform the constrained optimization by means of Algorithm~\ref{alg:convex_unconstrained}~(b) and by using the potential function one can find for~\eqref{equ:PDE:opt}, that is, a function $\varphi$ such that $u=\nabla\varphi$. Note that using our scheme and some numerical PDE-solver as an inner-loop (instead of the closed-form solution) is also possible,~\textit{e.g.}, one needs to solve a linear system, not over $\mathbb{R}$, but over $\mathbb{C}$. We select $8$ initial conditions uniformly from $[1,8]$ and show the convergence in Figure~\ref{fig:PDE}\footnote{See \url{http://wjongeneel.nl/PDE.gif} for an animated version of Figure~\ref{fig:PDE1}.}. Note in particular that the non-averaged iterates perform similar to their averaged counterparts. 
}
\end{example}

\section{Discussion}
\label{sec:discussion}
\subsection{On the necessity of leaving the real numbers}

Given the results from the previous section, one might wonder if this ``\textit{complex-lifting}'' is needed. Real single-point gradient estimators evidently exist, \textit{cf.}~\cite{ref:Flaxman}, but with problematic variance bounds for $\delta \downarrow 0$. The common solution is to bring back some relation with the (directional) derivative~\cite{ref:agarwal2010optimal,ref:nesterov2017random}. Hence, one might wonder if there is a purely real analogue to~\eqref{equ:cs:deriv}. The next proposition strongly hints at a negative answer.

\begin{proposition}[On the necessity of leaving the real numbers]
\label{prop:real}
Consider some open, convex set $\mathcal{D}\subseteq \mathbb{R}^n$ with $\mathrm{int}(\mathcal{D})\neq \emptyset$. Then, there does not exist a continuous map $G:\mathbb{R}\to \mathbb{R}$ such that for all real-analytic functions $f:\mathcal{D}\to \mathbb{R}$ 
\begin{equation*}
    \lim_{\delta\downarrow 0}G\left(\tfrac{1}{\delta}f(x+\delta y)\right)=Df(x)[y]=\langle \nabla f(x), y \rangle \quad \forall x\in \mathcal{D},\,y\in \mathbb{S}^{n-1}. 
\end{equation*}
\end{proposition}
\begin{proof}
As $f\in C^{\omega}(\mathcal{D})$ we can construct for sufficiently small $\delta$ and any $y\in \mathbb{S}^{n-1}$ the convergent Taylor series of $f(x+\delta y)$ around $x$ and as $G\in C^0(\mathbb{R})$ we can consider the limit in $\delta$ with respect to the argument of $G$, hence, we have $ \lim_{\delta\downarrow 0}G\left(\tfrac{1}{\delta} f(x+\delta y)\right)=G( \langle  \nabla f(x), y \rangle +\lim_{\delta\downarrow 0}\tfrac{1}{\delta}f(x))$. When $f(x)=0$ we end up with the fixed-point problem $G(\langle \nabla f(x), y \rangle)=\langle \nabla f(x), y \rangle$. As $\mathcal{D}$ is open, convex and with $\mathrm{int}(\mathcal{D})\neq \emptyset$, then for any $\alpha \in \mathbb{R}$ one can always find a pair $(f,y)$ such that $\langle \nabla f(x), y \rangle = \alpha$, \textit{e.g.}, construct a linear function over $\mathcal{D}$. Therefore, $G$ is forced to be the identity map on $\mathbb{R}$. Thereby, obstructing the case $f(x)\neq 0$. 
\end{proof}
Observe from the proof of Proposition~\ref{prop:real} that if we would generalize $G(\frac{1}{\delta} f(x+\delta y))$ to $\bar{G}(f(x+\delta y),\delta)$ with $\bar{G}$ continuous in $\mathbb{R}\times\mathbb{R}_{>0}$, the conclusion would not change.

\subsection{$C^{\infty}$-smooth imaginary zeroth-order optimization}
\label{sec:smooth}

Consider the smooth function $\psi:\mathbb{R}\to \mathbb{R}$ defined by $\psi(x)=|x|^2$. When evaluating $\psi$ at some complex point $z=x+iy\in \mathbb{C}$ one finds that $\psi(z)=x^2+y^2$, as such, $\psi$ does not satisfy the Cauchy-Riemann equations and is \textit{nowhere} (complex) analytic. This, however, means that one cannot appeal to the complex-step framework from~\cite{ref:JongeneelYueKuhnZO2021}~\textit{cf.}~Section~\ref{sec:im:grad}. Next, consider the prototypical smooth, yet non-analytic, function $\varphi:\mathbb{R}\to \mathbb{R}$ defined by
\begin{equation*}
    \varphi(x) = \begin{cases}
    \mathrm{exp}\left(\tfrac{-1}{x}\right)& \text{if }x>0\\
    0 & \text{otherwise}
    \end{cases}.
\end{equation*}
This function only fails to be analytic at $0$ and, interestingly, by the following expansions of $\mathrm{exp}(-1/z)$
\begin{equation}
    \label{equ:expz}
  \begin{aligned}
    \mathrm{exp}\left(\tfrac{-1}{x+iy} \right)
    = \mathrm{exp}\left(\tfrac{-x}{x^2+y^2}+i\tfrac{y}{x^2+y^2} \right)
= \mathrm{exp}\left(\tfrac{-x}{x^2+y^2} \right)\left(\cos\left(\tfrac{y}{x^2+y^2}\right)+i\sin\left(\tfrac{y}{x^2+y^2}\right)\right) 
\end{aligned}  
\end{equation}
one can readily show that $\varphi$ \textit{does} satisfy the Cauchy-Riemann equations. Indeed, recall~\eqref{equ:cs:deriv} and consider now the imaginary part of~\eqref{equ:expz}, then by the series expansion of $\mathrm{exp}(\cdot)$ and $\sin(\cdot)$ 
one observes that
\begin{equation*}
\partial_x\, \mathrm{exp}\left(\tfrac{-1}{x}\right)=\tfrac{1}{x^2}\mathrm{exp}\left(\tfrac{-1}{x}\right) = 
\lim_{\delta \downarrow 0}\,    \tfrac{1}{\delta}\mathrm{exp}\left(\tfrac{-x}{x^2+\delta^2} \right)\sin\left(\tfrac{\delta}{x^2+\delta^2}\right).
\end{equation*}
Hence, although $\varphi\in C^{\infty}\setminus C^{\omega}$, the complex-step framework is not obstructed. 

It turns out that from a topological point of view, the function $\varphi$ is somewhat of a special case. Let $\mathsf{X}$ be a topological space. Then the set $M\subset \mathsf{X}$ is of the \textit{\textbf{first category}}, in the sense of Baire, when $M$ is a countable union of nowhere dense sets in $\mathsf{X}$. A set $A\subseteq \mathsf{X}$ is said to be \textbf{\textit{nowhere dense}} when $\mathrm{cl}(A)^c$ is dense in $\mathsf{X}$, or equivalently, when $\mathrm{int}(\mathrm{cl}(A))=\emptyset$. Now one can show that under the sup-norm, the complement to the space of nowhere differentiable functions in $C^0([0,1])$ is of the first category~\cite[Chapter~5]{ref:Folland}. Differently put, almost every continuous function on $[0,1]$ is nowhere differentiable. A similar topological statement can be made about nowhere analytic functions in the space of smooth functions $C^{\infty}([0,1])$ under a sup-metric,~\textit{e.g.}, see\footnote{See in particular this post~\url{https://web.archive.org/web/20161009194815/mathforum.org/kb/message.jspa?messageID=387148} by Dave L. Renfro for more context.}~\cite{ref:darst1973most,ref:cater1984differentiable}. Again, bluntly put, almost every smooth function is nowhere analytic. An important question that comes with such an observation is where in the space of smooth functions optimization takes place?


\subsection{Future work}
This work exploits smoothness to be able to appeal to the Cauchy-Riemann equations. Other work, like~\cite{polyak1990optimal,ref:bach:smooth,ref:akhavan2020exploiting,ref:novitskii2021improved} exploit the knowledge of smoothness and construct kernels to (optimally) filter out (all) low-order errors. For increasing smoothness, however, we observe numerical instability in this approach, that is, the kernels become ill-defined. It would be worthwhile to further study how to exploit smoothness while taking the implementation into consideration. Given Proposition~\ref{prop:real}, it would also be interesting to explore the possibility of applying generalized versions of the complex-step approach,~\textit{e.g.}, using \textit{hyper-dual} numbers to extract second-order information~\cite{ref:fike2011}.

This work is mostly positioned within the scope of randomized methods via Lemma~\ref{lem:CR:grad}. Recent work indicated that in fact \textit{non}-randomized methods can outperform their randomized/smoothed counterparts~\cite{ref:berahas2021theoretical,ref:scheinberg2022finite}. This provides for interesting future work, especially in the presence of noise. Estimating the noise statistics itself also provides for relevant future work as it allows for a more appropriately scaled sequence of smoothing parameters. 



\subsection{Conclusion}
\label{sec:conclusion}

We have presented a line of algorithms that can theoretically \textit{and} practically deal with any suitable sequence $\delta_k \to 0$ (conditioned on appropriate stepsizes $\{\mu_k\}_k$). In contrast to~\cite{ref:JongeneelYueKuhnZO2021} we can also deal with computational noise and demand less prior knowledge of problem parameters. 

Only if we understand all the vulnerabilities of our algorithms | as esoteric as they are | we can safely implement them. With that in mind, we hope this work provides for more future work on numerical optimization.


\appendix
\thispagestyle{empty}
\section*{Appendix}
This appendix contains auxiliary results related to the work above.

\section{Auxiliary results}
The following results are well-known.  
\begin{lemma}[Logarithm bound]
\label{lem:log}
For any $J\in \mathbb{N}_{\geq 1}$ one has
\begin{equation}
\label{equ:log:bound}
    \textstyle\sum^J_{j=1}\tfrac{1}{j}\leq \log(J)+1.
\end{equation}
\end{lemma}

\begin{lemma}[Fractional bound]
\label{lem:fractional}
For any $\beta\geq 1$ one has
\begin{equation}
    \textstyle\sum^J_{j=1}j^{-1+1/\beta}\leq \beta J^{1/\beta}. 
\end{equation}
\end{lemma}

\section{Estimation of $\bar{\tau}(f)$ via diagonally dominant programming}
\label{sec:DDP}

We highlight the \textit{basis pursuit} approach as proposed in~\cite{ref:ahmadi2017seq}. 
A constraint of the form $P\in \mathsf{dd}^n$ translates to a set of linear constraints.
The same is true for $P_z\in \mathsf{dd}^n(U_z)$ with
\begin{equation*}
    \mathsf{dd}^n(U_z)=\{M\in \mathcal{S}^n: M=U_z^{\mathsf{T}}QU_z,\,Q\in \mathsf{dd}^n\}
\end{equation*}
for some \textit{basis} matrix $U_z$.
Now to iteratively change the basis $U_z$ one can use 
\begin{equation*}
    U_{z+1}=\mathsf{chol}(P_z),\quad U_0=I_n,
\end{equation*}
for $P_z$ the solution of the $z^{\mathrm{th}}$ program. By construction one has $P_z\in \mathsf{dd}^n(U_{z+1})$ such that each new iteration is at least as good as the previous one. As by the construction in Section~\ref{sec:est:tau} we demand that $P_z\succeq \tau_0 I_n\succ 0$, then, by~\cite[Theorem~3.1]{ref:ahmadi2017seq} $P_z\to P^{\star}$ (weakly) for $z\to +\infty$ and $P^{\star}$ being the (a) solution of the original problem. In practice, one could terminate the algorithm when $Q_z$ is sufficiently close to $I_n$ and set $\widehat{\tau}(f)=\min_i\{\lambda_i(P_z)\}$, which can be found using a dedicated large-scale algorithm.

\begin{figure*}[t!]
    \centering
    \begin{subfigure}[b]{0.3\textwidth}
        \includegraphics[width=\textwidth]{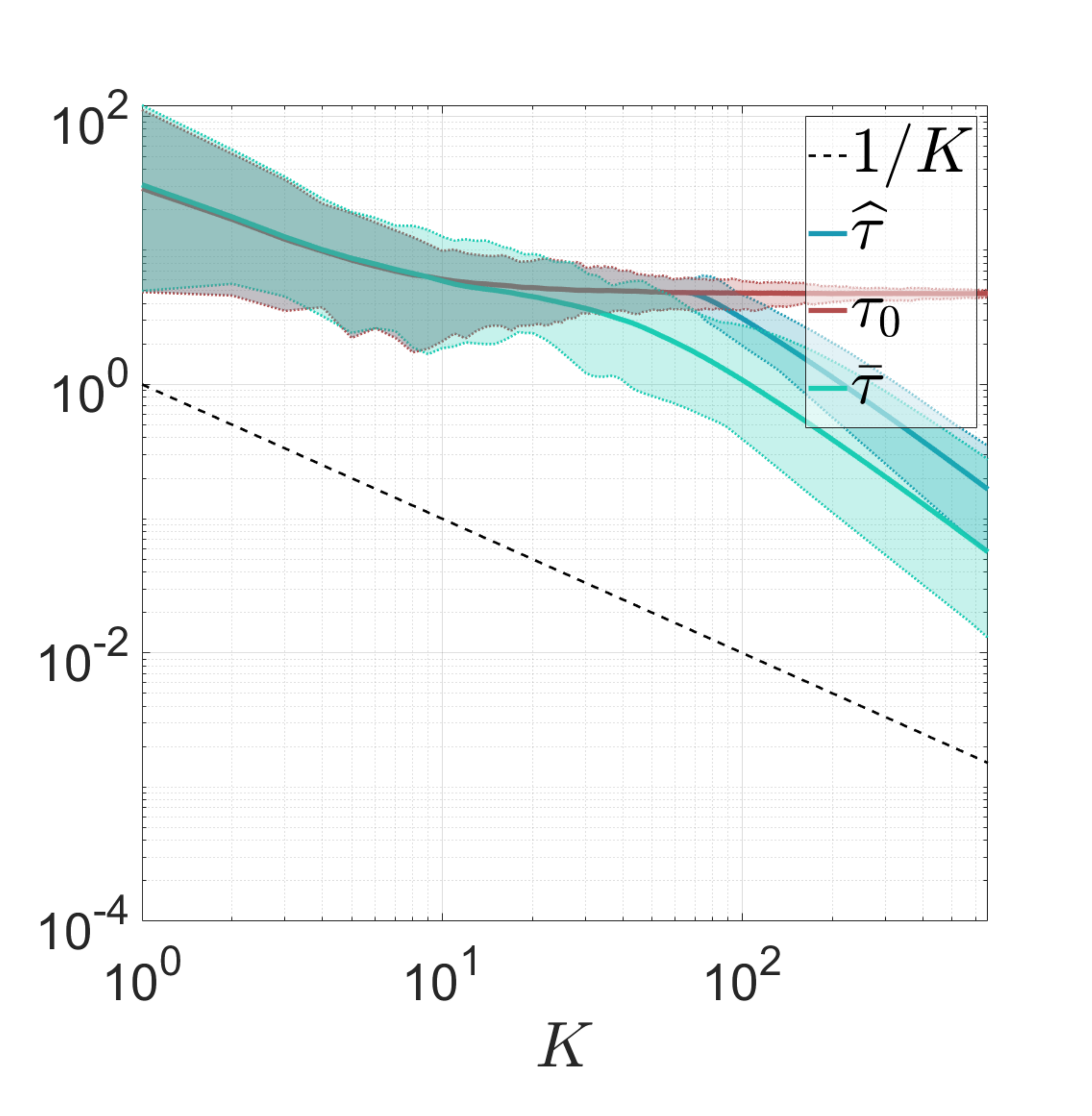}
        \caption{Algorithm~\ref{alg:convex_unconstrained}~(b) in combination with DDP-based estimation of $\bar{\tau}$.}
        \label{fig:DDPcost}
    \end{subfigure}\quad
    \begin{subfigure}[b]{0.3\textwidth}
        \includegraphics[width=\textwidth]{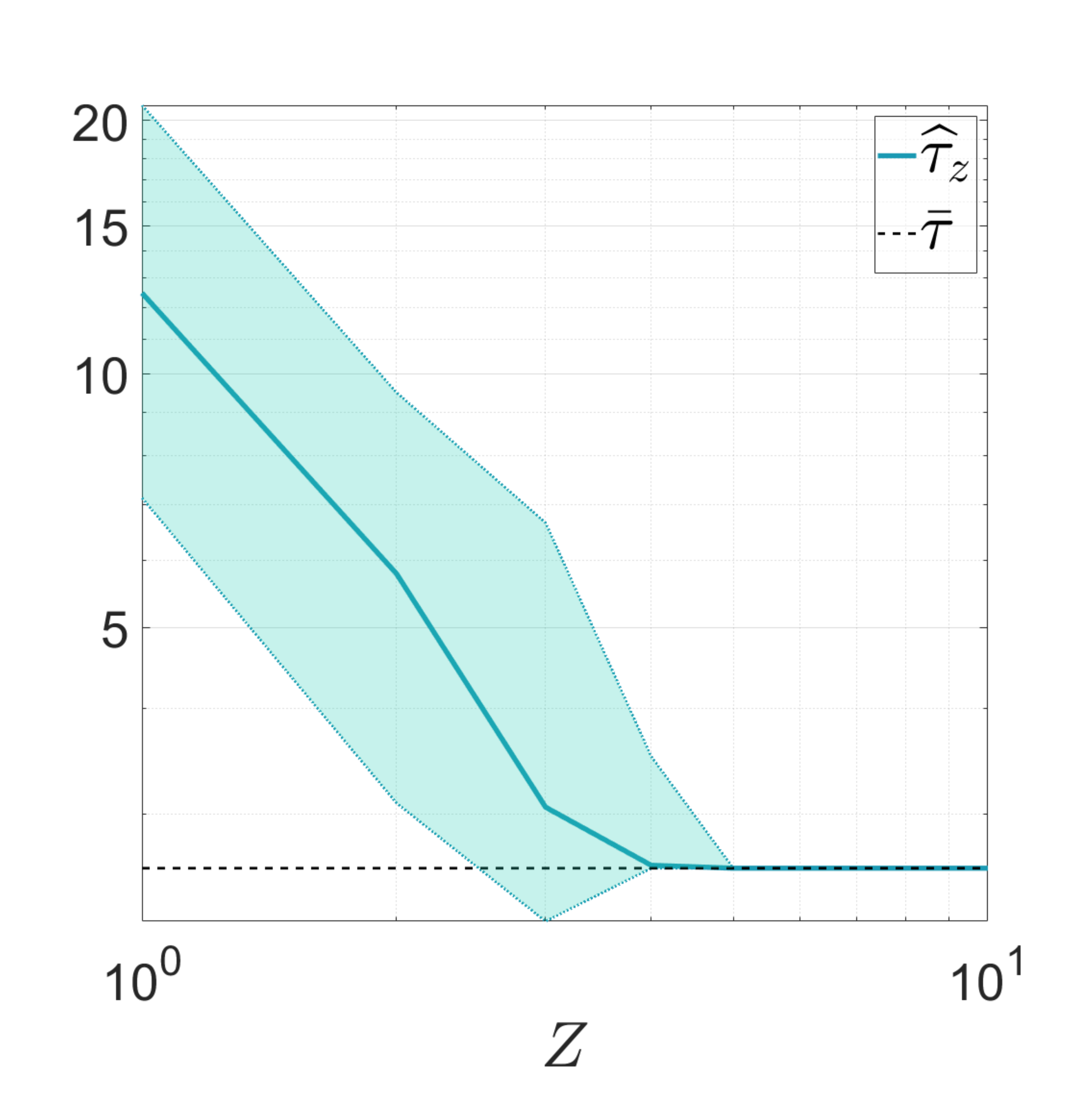}
        \caption{Basis pursuit convergence for the program~\eqref{equ:DDP:tau}.\\}
        \label{fig:tauconv}
    \end{subfigure}
    \caption[]{Numerical outcomes of the DDP example (Section~\ref{sec:DDP}).}
    \label{fig:DDP}
\end{figure*}

To showcase the approach we redo Example~\ref{ex:num:est:stab}, but by using~\eqref{equ:DDP:tau}. Here, we fix a random pair $(A,b)$ and show for $100$ initial conditions $x_1\overset{i.i.d.}{\sim}\mathcal{N}(0,I_n)$ the effect of an improved estimate of $\bar{\tau}$. Here, we apply the basis pursuit approach as sketched above for $n$ iterations and set $\widehat{\tau}= \widehat{\tau}_n$. 
The results are shown in Figure~\ref{fig:DDP}. Again, we observe the benefit of estimating $\bar{\tau}$, plus, we see that the inner-routine convergences quickly, yet, usually from above. Quantifying the behaviour as seen in Figure~\ref{fig:tauconv} would be interesting and is left for future work.  

\section{Lipschitz inequalities}
In this section we gather a variety of inequalities which come in useful later. Note that convexity of $f$ is usually not a necessary assumption. 
If $f$ is convex then, by~\cite[Theorem~2.1.5]{nesterov2003introductory}~\eqref{equ:grad:Lipschitz} implies that 
\begin{equation}
    \label{equ:grad:Lipschitz:bound2}
    f(x)\geq f(y) + \langle \nabla f(y) , x -y \rangle + \tfrac{1}{2 L_1(f)}\|\nabla f(x) - \nabla f(y) \|_2^2,\quad \forall x,y \in \mathcal{D} 
\end{equation}
and thus for any (local) minimum $x^{\star}$ such that $\nabla f(x^{\star})=0$ one has $2L_1(f)\left(f(x)-f(x^{\star})\right)\geq \|\nabla f(x)\|_2^2$. 
Also, as~\cite[Equation (6)]{ref:Nesterov:2011}, for $f\in C^{1,1}_{L_1(f)}(\mathcal{D})$ one has
\begin{equation}
    \label{equ:L1:upper}
    |f(y)-f(x) - \langle \nabla f(x),y-x \rangle | \leq \tfrac{1}{2}L_1(f) \|x-y\|_2^2, \quad \forall x,y\in \mathcal{D}.  
\end{equation}
It follows from~\cite[Lemma 1.2.4]{nesterov2003introductory} that if $f\in C^2(\mathcal{D})$ then
\begin{equation}
\label{equ:L2:bound}
    |f(y)-f(x) - \langle \nabla f(x), y-x\rangle - \tfrac{1}{2}\langle \nabla^2 f(x)(y-x), y-x \rangle |\leq \tfrac{1}{6}L_2(f) \|x-y\|_2^3,\quad \forall x,y \in \mathcal{D}.
\end{equation}
See that~\eqref{equ:Hess:Lip} is equivalent to
\begin{equation}
    \label{equ:thirdorder:smooth}
    |\langle \nabla^2 f(x)u, u \rangle - \langle \nabla^2 f(y)u, \rangle | \leq L_2(f)\|x-y\|_2 \quad \forall x,y \in \mathcal{D},\;u\in \mathbb{S}^{n-1}, 
\end{equation}
which is commonly referred to as $f$ being $3$rd-order smooth, \textit{cf.}~\cite[Section~1.1]{ref:bach:smooth}. Now it follows directly from~\eqref{equ:L2:bound} and the definition of a derivative that $f\in C^{3,2}_{L_2(f)}(\mathcal{D})$ implies that for all $x\in \mathcal{D}$ one has 
\begin{equation}
    \label{equ:third:deriv:bound}
    \left|\partial^3_t f(x+tu)|_{t=0}\right|\leq L_2(f),\quad \forall u\in \mathbb{S}^{n-1}.
\end{equation}

\section{Further numerical comments}
\label{sec:app:num}

\begin{figure*}[t!]
    \centering
    \begin{subfigure}[b]{0.3\textwidth}
        \includegraphics[width=\textwidth]{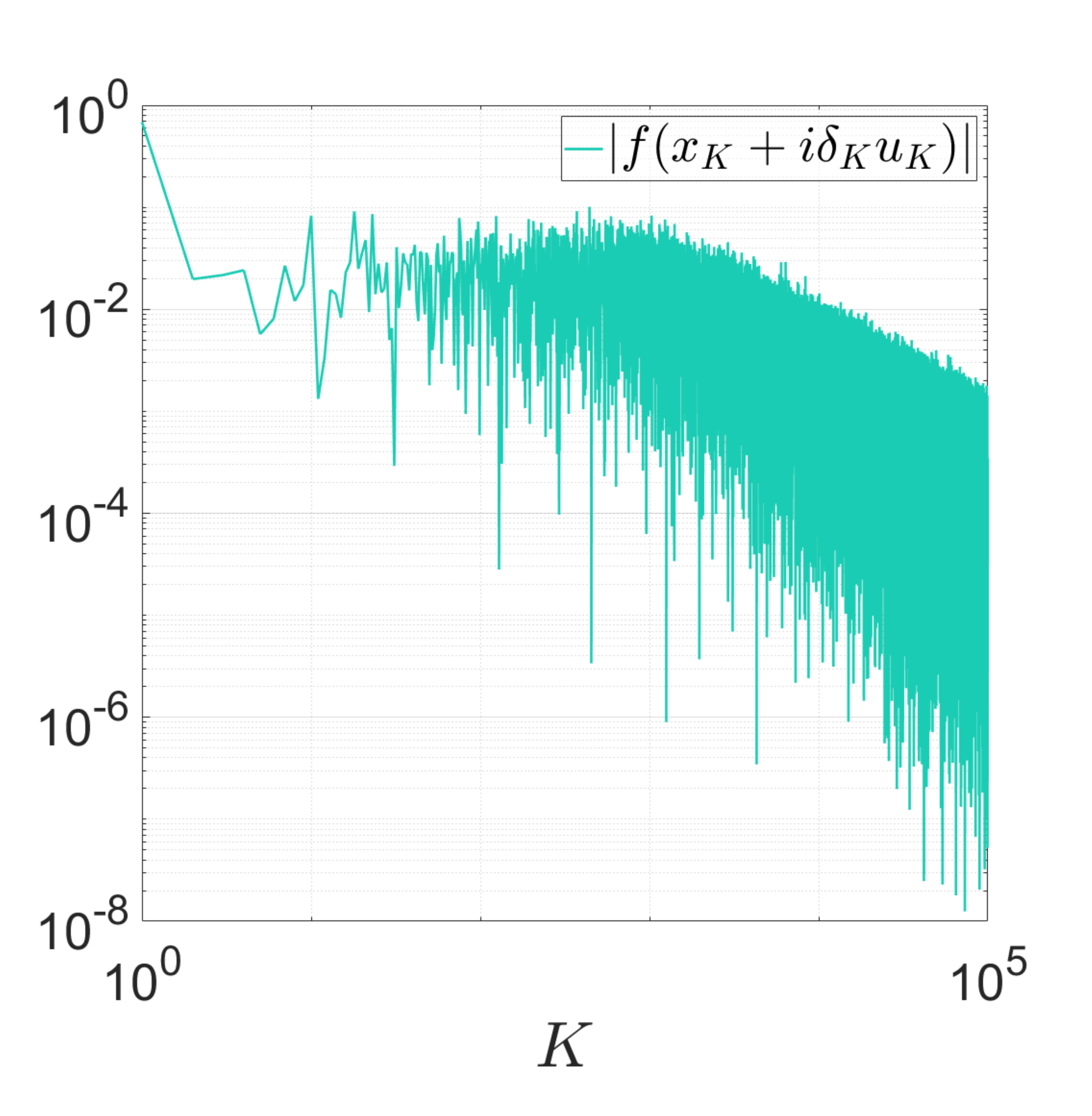}
        \caption{CS ($\delta=1$) estimator}
        \label{fig:strongcvxCS0f}
    \end{subfigure}\quad
    \begin{subfigure}[b]{0.3\textwidth}
        \includegraphics[width=\textwidth]{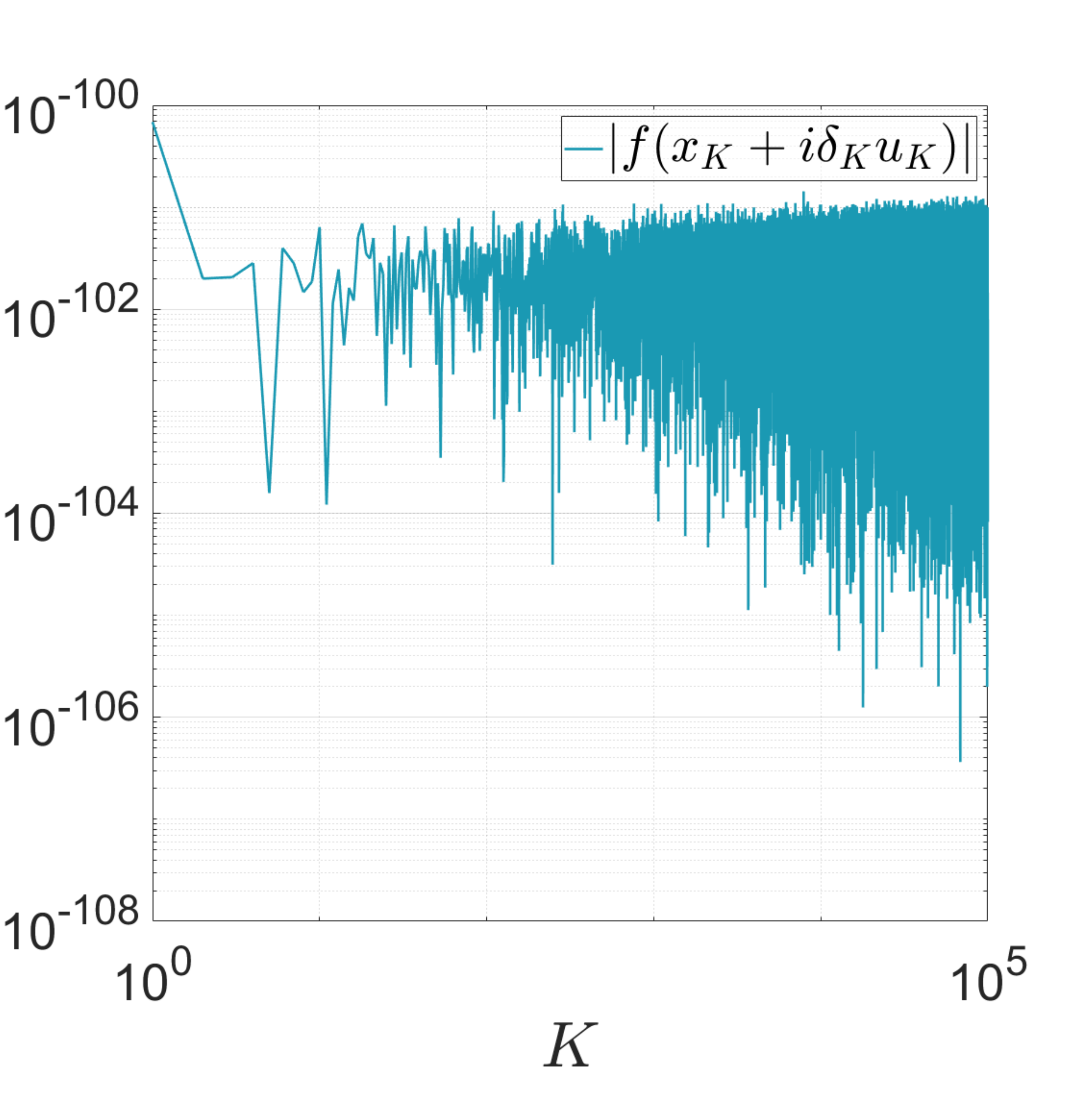}
        \caption{CS ($\delta=10^{-100}$) estimator}
        \label{fig:strongcvxCSf}
    \end{subfigure}\quad 
    \begin{subfigure}[b]{0.3\textwidth}
        \includegraphics[width=\textwidth]{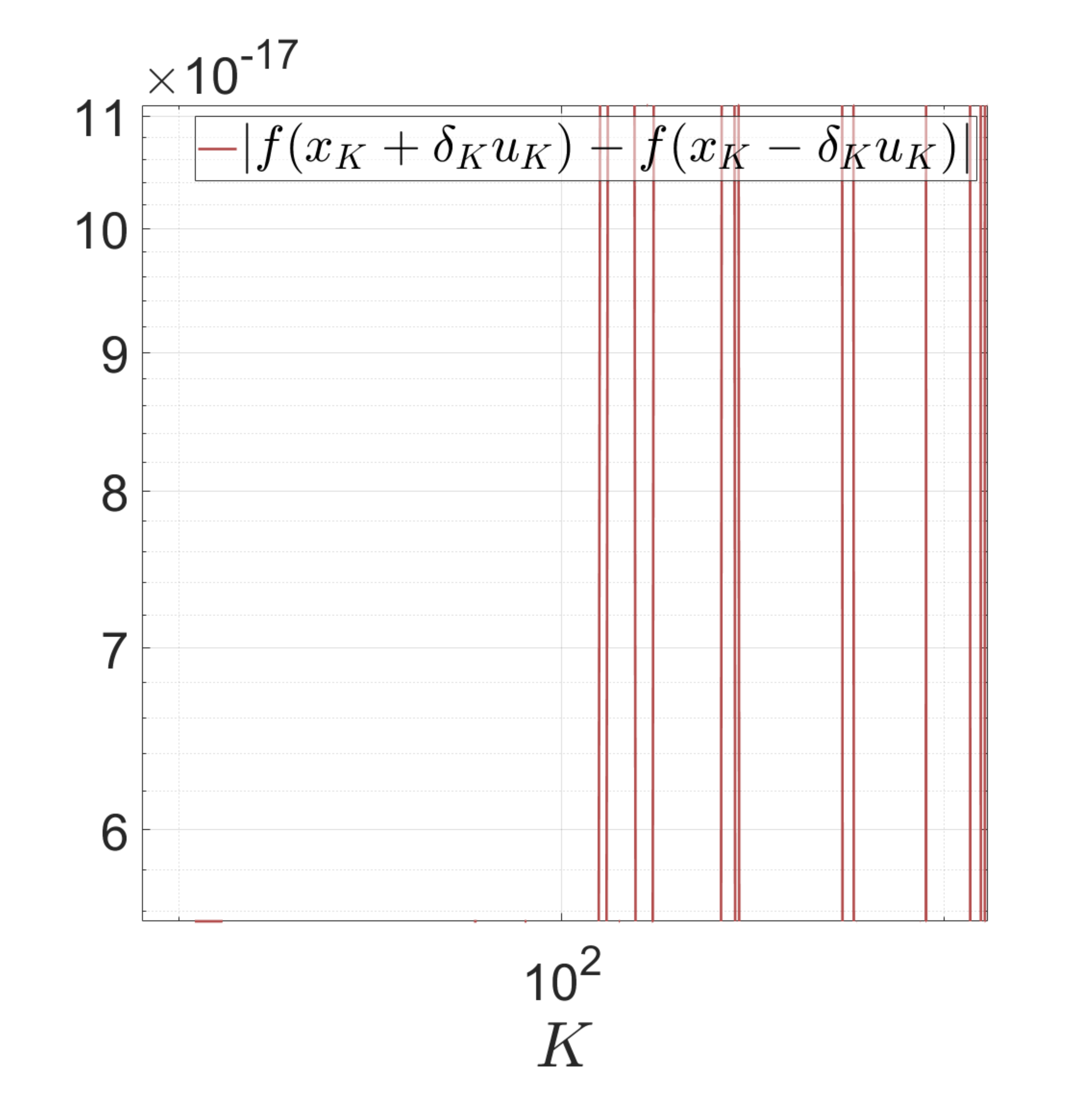}
        \caption{$\beta$ estimator}
        \label{fig:strongcvxMf}
    \end{subfigure}
    \caption[]{Further remarks on Example~\ref{ex:strongcvx}.}
    \label{fig:sim:cont}
\end{figure*}

Example~\ref{ex:strongcvx} continued. In Figure~\ref{fig:strongcvx} we see a clear difference in behaviour. This can be explained by looking at the corresponding estimators. We see that for the estimator as proposed in this work no cancellation occurs, while for the frequently employed central-difference scheme as used in~\cite{ref:akhavan2020exploiting} the two function evaluations can cancel catastrophically. See Figure~\ref{fig:strongcvxCS0f}-\ref{fig:strongcvxCSf} and Figure~\ref{fig:strongcvxMf}. We like to remark, in line with the analysis, that the scheme for $\delta=1$ is better conditioned. 

All numerical experiments are carried out in MATLAB using the SDPT3 solver~\cite{ref:toh1999sdpt3}. 


\paragraph*{Data availability statement}
All data generated or analysed during this study are included in this article.

\paragraph*{Conflict of interest}
The author has no competing interests to declare that are relevant to the content of this article.

\pagestyle{basicstyle}
\addcontentsline{toc}{section}{Bibliography}
\subsection*{Bibliography}
\printbibliography[heading=none]


\end{document}